\documentclass{conm-p-l}

   \usepackage{amsthm}
   \usepackage{amsmath}
   \usepackage{amssymb}
   \usepackage{amsfonts}
   \usepackage{amscd}
   \usepackage[mathscr]{eucal}
   \usepackage{verbatim}
   \usepackage{graphics}

\input xy
\xyoption{all}

\newtheorem{athm}{Theorem}[subsection]
\newtheorem{alem}[athm]{Lemma}
\newtheorem{aprop}[athm]{Proposition}

\theoremstyle{definition}
 \newtheorem{arem}[athm]{Remark}
 \newtheorem{arems}[athm]{Remarks}
 \newtheorem{adefi}[athm]{Definition}

\numberwithin{equation}{athm}

\newcommand{\ssss}{\stepcounter{athm}\subsubsection{}}

\newcommand{\EQAL}[1]%
{\,\begin{picture}(#1,0)%
\put(0,3){\line(1,0){#1}}%
\put(0,1){\line(1,0){#1}}%
\end{picture}\,}%

\newcommand{\vlto}[1]%
{\,\begin{picture}(#1,3)%
\put(0,2){\vector(1,0){#1}}%
\end{picture}\,}%

\newcommand{\vllarrow}[1]%
{\,\begin{picture}(#1,3)%
\put(#1,2){\vector(-1,0){#1}}%
\end{picture}\,}%

\newcommand{\dirlm}[1]%
  {
     {\lim\hskip-1.58em\lower.65ex
       \hbox{$
                {}_{\stackrel{\lower1ex\hbox
                                        {$\scriptstyle -\!\!\<\longrightarrow$}
                             }{ ^{#1} }
                   }
            $}
     }
\:}

\newcommand{\subdirlm}[1]%
  {
     {\lim\hskip-1.5em\lower.6ex
       \hbox{$
                   {}_{\stackrel{\lower1ex\hbox
                                           {$\scriptstyle\longrightarrow$}
                                }{ ^{#1} }
                      }
             $}
     }
\:}

\newcommand{\inlm}[1]%
   {
      {\lim\hskip-1.58em\lower.65ex
        \hbox{$
                 {}_{\stackrel{\lower1ex\hbox
                                        {$\scriptstyle \longleftarrow\!\!\<-$}
                              }{ ^{#1} }
                    }
             $}
      }
\:}

\def\>{\mspace {1mu}}
\def\<{\mspace{-1mu}}
\def\({{\textup(}}
\def\){{\textup)}}

\newcommand{\X}{{\mathscr X}}

\newcommand{\Y}{{\mathscr Y}}
\newcommand{\Z}{{\mathscr Z}}
\newcommand{\V}{{\mathscr V}}

\newcommand{\W}{{\mathscr W}}
\newcommand{\I}{{\mathscr I}}

\newcommand{\sS}{{\mathscr S}}
\newcommand{\sT}{{\mathscr T}}
\newcommand{\sU}{{\mathscr U}}
\newcommand{\sV}{{\mathscr V}}

\newcommand{\eC}{{\mathscr C}}

\newcommand{\sfC}{{\mathsf C}}
\newcommand{\sfI}{{\mathsf I}}
\newcommand{\sfF}{{\mathsf F}}
\newcommand{\sfO}{{\mathsf O}}
\newcommand{\sfP}{{\mathsf P}}
\newcommand{\sfQ}{{\mathsf Q}}
\newcommand{\sfR}{{\mathsf R}}

\newcommand{\A}{{\mathcal A}}
\newcommand{\B}{{\mathcal B}}
\newcommand{\C}{{\mathcal C}}
\newcommand{\cD}{{\mathcal D}}

\newcommand{\cS}{{\mathcal S}}

\newcommand{\cFb}{{\mathcal F}^{\bullet}}

\newcommand{\cO}{{\mathcal O}}

\newcommand{\rh}{{\mathrm {h}}}
\newcommand{\rv}{{\mathrm {v}}}

\newcommand{\bC}{{\mathbf C}}
\newcommand{\D}{{\mathbf D}}

\newcommand{\Dqc}{\D_{\mkern-1.5mu\mathrm {qc}}}

\newcommand{\Dqct}{\D_{\mkern-1.5mu\mathrm{qct}}}

\newcommand{\bt}{{\bf t}}
\newcommand{\bs}{{\bf s}}

\newcommand{\boldC}{\boldsymbol{C}}
\newcommand{\bbeta}{\boldsymbol{\beta}}
\newcommand{\bphi}{\boldsymbol{\phi}}

\newcommand{\cfr}{\mathfrak c}

\newcommand{\lfr}{\mathfrak l}

\newcommand{\pfr}{\mathfrak p}
\newcommand{\qfr}{\mathfrak q}
\newcommand{\rfr}{\mathfrak r}
\newcommand{\sfr}{\mathfrak s}

\newcommand{\R}{{\mathbf R}}

\newcommand{\Hom}{{\mathrm {Hom}}}

\newcommand{\iGp}[1]{{\varGamma_{\<\!#1}'}}
\newcommand{\iG}[1]{{\varGamma_{\<\!#1}^{\phantom\prime}}}

\newcommand{\set}{\!:=}

\newcommand{\Lra}{\Longrightarrow}

\newcommand{\xto}[1]{\xrightarrow{#1}}

\newcommand{\ssbox}{{\scriptscriptstyle{\Box}}}
\newcommand{\bxt}{\scriptscriptstyle{\boxtimes}}
\newcommand{\sbsq}{{\scriptstyle{\blacksquare}}}

\newcommand{\wt}[1]{{\widetilde{#1}}}
\newcommand{\ov}[1]{{\overline{#1}}}

\newcommand{\wtcD}{\wt{\cD}}

\newcommand{\wtDqc}{\wt{\D}_{\mkern-1.5mu\mathrm {qc}}}

\newcommand{\sbeta}{\textup{\ss}}

\newcommand{\dg}[1]{\ensuremath{\dagger_{\<\raisebox{-.2ex}{\scriptsize{#1}}}}}
\newcommand{\oneD}[1]{\smash{\mathbf{1}_{\<\<%
\raisebox{-.1ex}{$\scriptstyle\cD_{\<\<%
\raisebox{-.1ex}{$\scriptscriptstyle{#1}$}}$}}}}
\newcommand{\xyoneD}[1]{\;\;\mathbf{1}_{\<\<%
\raisebox{-.1ex}{$\scriptstyle\cD_{\<\<%
\raisebox{-.1ex}{$\scriptscriptstyle{#1}$}}$}}}

\newcommand{\iso}%
{{\mkern8mu\longrightarrow \mkern-25.5mu{}^\sim\mkern17mu}}
\newcommand{\osi}%
{{\mkern8mu\longleftarrow \mkern-24.5mu{}^\sim\mkern16mu}}
\def\Otimes{\underset
  {\vbox to 0pt {\vskip-1ex\hbox{$\scriptscriptstyle=$}\vss}}
    {\otimes}\vadjust{\kern.4pt}}

\newcommand{\smcirc}%
  {{\raise.15ex\hbox to.7em{$\hss \scriptstyle\circ\hss$}}}

\title{Pasting pseudofunctors}
\author[ S.\:Nayak]{Suresh Nayak}
\address{Chennai Mathematical Institute \\
   92, G.\,N.\,Chetty Road\\
   Chennai-600017, INDIA}
\email {snayak@cmi.ac.in}

\thanks{The author thanks the Mathematisches 
Forschungsinstitut Oberwolfach, the Banff International Research
Station and the Institute of Mathematical Sciences at Chennai
for providing access to their conducive environment in which much of
this research was carried out. The author was funded by the 
National Board of Higher Mathematics.}

\begin{document}

\begin{abstract}
We give an abstract criterion for pasting pseudofunctors
on two subcategories of a category into a pseudofunctor 
on the whole category. As an application we extend the
variance theory of the twisted inverse image $(-)^!$ over
schemes to that over formal schemes.
Thus we show that over composites of compactifiable maps of noetherian 
formal schemes, there is a pseudofunctor $(-)^!$ 
such that if $f$ is a pseudoproper map, then $f^!$ is the right adjoint to the
derived direct image $\R f_*$ and if $f$ is \'{e}tale, then $f^!$ is the
inverse image $f^*$. We also show that $(-)^!$ is compatible with 
flat base change.
\end{abstract}

\maketitle

\tableofcontents

\renewcommand{\theenumi}{\roman{enumi}}

\section{Introduction}

\subsection{Background}
\label{subsec:background}
The purpose of this paper is to give an abstract criterion 
for pasting pseudofunctors on two subcategories of a category 
into a pseudofunctor on the whole category. The inspiration for 
the criterion comes from applications lying in the theory of
Grothendieck duality. 

For the sake of this introduction,
let $\mathbb S$ be the category of separated finite-type maps
of noetherian schemes and $\mathbb F$ the category of
separated pseudofinite-type maps of noetherian formal schemes. 
Here $\mathbb S$ can be thought of as the full subcategory 
of~$\mathbb F$ consisting of formal schemes whose structure sheaf of rings 
has discrete topology.

One of the main applications of our pasting result 
concerns extending the variance theory of the twisted-inverse-image
functor $(-)^!$ to the setup of formal schemes. 
Recall that, in \cite[Theorem 2]{AJL2} it has been shown that for any map 
$f \colon \X \to \Y$ in $\mathbb F$, the functor
$\R f_* \colon \Dqct^+(\X) \to \Dqct^+(\Y)$ has a right 
adjoint $f^{\times}$, where $\Dqct^+(-)$ stands for the derived category of 
bounded-below complexes whose homology modules are quasi-coherent 
and torsion.
Let $(-)^{\times}$ be the resulting pseudofunctor on $\mathbb F$.
Following the case of ordinary schemes, a basic 
goal in the abstract side of duality theory 
over formal schemes is of
constructing a pseudofunctor
$(-)^!$, possibly over the whole of $\mathbb F$,
such that the restriction of~$(-)^!$ to the subcategory of
pseudoproper maps is isomorphic to $(-)^{\times}$, and
the restriction to open immersions is the 
inverse-image pseudofunctor $(-)^*$. Thus the problem of constructing 
$(-)^!$ is a problem of pasting pseudofunctors.

It is known that a solution for this 
pasting problem exists over $\mathbb S$, henceforth to be also called as the
classical case, though 
a complete presentation of details
has not appeared in published form yet. But there have been
serious hurdles in carrying over the approach of 
the classical case to that over $\mathbb F$.
Let us first briefly recall the main inputs  
used in the classical case.

Over $\mathbb S$, the main ingredients used for 
pasting are (i) Nagata's compactification
theorem (\cite{Lu}, \cite{Nag}) which says that every map $f$ in $\mathbb S$
factors as $f = pi$ where $p$ is a proper map and $i$ an open immersion;
(ii) The open-base-change theorem due to Deligne and more generally,
the flat-base-change theorem of Verdier (\cite{Ve}).
While~(ii) has been generalized to the formal-scheme setup 
(\cite[Theorem 3]{AJL2}),
we do not know if the analog of (i) for formal schemes holds. 

The nonavailability of (i) has been a critical hurdle for pasting 
over $\mathbb F$. In fact, pasting
even within the limited scope of
compactifiable maps in $\mathbb F$ had not been achieved
since the classical proofs, in their intermediate stages,
rely on some form of~(i).
For instance, in the situation of a map~$f$ in~$\mathbb F$
admitting two factorizations $f= p_1i_1 = p_2i_2$  
where $p_j$ are pseudoproper 
and $i_j$ open immersions, a canonical choice for an isomorphism 
$i_1^*p_1^{\times} \iso i_2^*p_2^{\times}$ had not 
been demonstrated till now. (see~\cite{AJL2a})

The functorial details that need to be addressed in these
pasting problems are best tackled in an abstract framework. 
In \cite{De} Deligne gave an abstract criterion for pasting 
pseudofunctors. But its requirement that factorizations 
such as in (i) above hold in the working category is a serious drawback.
Moreover, even in some contexts where (i) does hold,
satisfying Deligne's input conditions seems to require 
considerable functorial manipulations which could possibly be 
carried out at an abstract level itself.

Against this backdrop we offer a way of pasting that does not 
rely on the existence of factorizations as in (i). Instead, 
our approach exploits compatibilities concerning what we call
fundamental isomorphisms.

In the context of pasting $(-)^{\times}$ with $(-)^*$, 
a fundamental isomorphism is one
associated to every factorization 
of an identity map $1_{\X} = gf$ with~$g$ open and~$f$ pseudoproper,
and is of the form 
$f^{\times}g^* \iso \mathbf{1}_{\Dqct^+(\X)}$. 
Such an isomorphism is easily constructed: 
in this case $f,g$ are necessarily 
isomorphisms and hence we may use $g^* = g^{\times}$ and 
$f^{\times}g^{\times} \iso (1_X)^{\times}$.
Nevertheless, these fundamental isomorphisms play a very
useful role in pasting.

Our principal results are abstract in nature in the spirit of Deligne's
pasting result. Before summarizing, we point out the main application
to duality that comes out of them.

Loosely speaking, we now have the following result (see Theorems 
\ref{thm:twisted4}, \ref{thm:twisted5} and~\ref{thm:twisted6}). 

{\it{Over the subcategory of composites of  
compactifiable maps in $\mathbb F$ (equivalently, the subcategory 
of composites of open immersions and pseudoproper maps in $\mathbb F$)
there is a pseudofunctor $(-)^!$ whose restriction to the subcategory 
of pseudoproper maps is isomorphic to $(-)^{\times}$  and whose
restriction to the 
subcategory of open immersions is $(-)^*$, and which 
furthermore, is uniquely determined via these isomorphisms by its 
compatibility with open-base-change isomorphisms. Moreover $(-)^!$ 
satisfies compatibility with flat base change and can
also be extended to the subcategory of composites 
of \'{e}tale maps and pseudoproper maps in $\mathbb F$.}}

Note however that since we do not know whether the 
subcategory of composites of  
compactifiable maps in $\mathbb F$ (or even composites of \'{e}tale 
and pseudoproper) is the whole of $\mathbb F$,
therefore we do not know if $(-)^!$ as constructed by us is 
defined over the whole of $\mathbb F$. 

Finally, here is a rough summary of
our abstract pasting result. (Theorem~\ref{thm:output4})

For its input, we assume that there is a category $\sfC$ and 
there are pseudofunctors~$(-)^{\times}$
and~$(-)^{\ssbox}$ defined over two subcategories $\sfP$ and $\sfO$ 
respectively. Furthermore, we assume that to every fibered square $\sfr$ 
resulting from the fibered product of a $\sfP$-map with an $\sfO$-map there
is an associated ``base-change'' isomorphism $\beta_{\sfr}$ and that to every 
factorization~$\sigma$ of an identity map into a~$\sfP$-map 
followed by an $\sfO$-map,
there is an associated ``fundamental'' isomorphism $\phi_{\sigma}$.
The comparison maps of $(-)^{\times}, (-)^{\ssbox}$, the 
base-change isomorphisms $\beta_{-}$ and the fundamental 
isomorphisms~$\phi_{-}$ are assumed to be compatible in certain ways. 

For the output we obtain a pseudofunctor $(-)^!$ on the smallest 
subcategory~$\sfQ$ of~$\sfC$ containing~$\sfP$ and~$\sfO$, 
together with isomorphisms
between the restrictions of~$(-)^!$ to~$\sfP,\sfO$ and 
$(-)^{\times}, (-)^{\ssbox}$ respectively such that $(-)^!$
is compatible with~$\beta_{-}$ and~$\phi_{-}$.

We also have an abstract analog for flat base change. 
(Theorem~\ref{thm:flat5})

\subsection{Localness on source of torsion twisted inverse image}
\label{subsec:local}
In order to illustrate some of our main arguments and concerns, we take
up here the following task. 
Let $f \colon \X \to \Y$ be a map in $\mathbb F$ and
let $\X \xto{i} \Z_1 \xto{p} \Y$  and 
$\X \xto{j} \Z_2 \xto{q} \Y$ be factorizations
of $f$ such that $p,q$ are pseudoproper and $i,j$ are open immersions.
Let us call these factorizations as $\sigma_1 \set (i,p)$ and 
$\sigma_2 \set (j,q)$.
Our aim is to give a canonical isomorphism 
$c(\sigma_2, \sigma_1) \colon i^{*}p^{\times} \iso j^{*}q^{\times}$.

Consider the following diagram where the squares involved are 
fiber-product diagrams.
\[
\begin{CD}
\X @>{\Delta}>> \W_{11} @>{i_2}>> \W_{12} @>{p_2}>> \X \\
@. @V{j_2}VV @V{j_1}VV @VV{j}V \\
@. \W_{21} @>{i_1}>> \W_{22} @>{p_1}>> \Z_2 \\
@. @V{q_2}VV @V{q_1}VV @VV{q}V \\
@. \X @>>{i}> \Z_1 @>>{p}> \Y
\end{CD}
\]
We claim that there are natural isomorphisms
\begin{equation*}
\tag{\dag}
\Delta^{\!\times}j_2^*q_2^{\times} \iso \mathbf{1} 
\osi \Delta^{\!\times}i_2^*p_2^{\times}
\end{equation*} 
where~$\mathbf{1}$ is the identity on $\Dqct^+(\X)$. Assuming the claim 
we proceed as follows. Consider isomorphisms $\beta_i$ for $i=1,2,3,4$ 
as shown below where $\beta_1$ and $\beta_4$ are the 
obvious pseudofunctorial isomorphisms, and $\beta_2, \beta_3$ result 
from the base-change isomorphism \cite[Theorem 3]{AJL2}
corresponding to the 
northeastern and southwestern squares above respectively.
\[
j_2^*i_1^* \xto{\beta_1} i_2^*j_1^* \qquad   
j_1^*p_1^{\times} \xto{\beta_2} p_2^{\times}j^* \qquad   
q_2^{\times}i^* \xto{\beta_3} i_1^*q_1^{\times} \qquad
q_1^{\times}p^{\times} \xto{\beta_4} p_1^{\times}q^{\times}
\]
Then the desired isomorphism $c(\sigma_2,\sigma_1)$
is obtained via the following ones 
\begin{align}
i^{*}p^{\times} \osi \Delta^{\!\times}j_2^*q_2^{\times}i^{*}p^{\times}
\xto{\text{via }\beta_3} \Delta^{\!\times}j_2^*i_1^*q_1^{\times}p^{\times}
&\xto{\text{via }\beta_1} \Delta^{\!\times}i_2^*j_1^*q_1^{\times}p^{\times} 
 \notag \\
&\xto{\text{via }\beta_4} \Delta^{\!\times}i_2^*j_1^*p_1^{\times}q^{\times}
 \notag \\
& \xto{\text{via }\beta_2} \Delta^{\!\times}i_2^*p_2^{\times}j^*q^{\times}
 \iso j^*q^{\times}. \notag
\end{align}

For the isomorphisms in (\dag), let us begin with 
the first one, viz.,
$\Delta^{\!\times}j_2^*q_2^{\times} \iso \mathbf{1}$.
Consider the following diagram described below.
\[
\begin{CD}
@. \X \\
@. @V{a_1}VV \\
@. \V @>{a_2}>> \X \\
@. @V{a_3}VV @VV{a_4}V \\
\X @>>{\Delta}> \W_{11} @>>{j_2}> \W_{21} @>>{q_2}> \X
\end{CD}
\]
Here $a_4 \set j_2\Delta$, the square is  
obtained by taking fibered products and $a_1$ is the diagonal map. 
Thus $a_1,a_3,a_4$ are all
pseudoproper, in fact closed immersions.
In particular, via pseudofunctoriality we obtain 
isomorphisms
\[
\Delta^{\!\times} \iso a_1^{\times}a_3^{\times}, \qquad \qquad
a_4^{\times}q_2^{\times} \iso \mathbf{1}. 
\]
Since $a_2a_1 = 1_{\X}$, it follows that $a_2$ is surjective. Since $a_2$ 
is also an open immersion, therefore it is an isomorphism. In particular,
$a_2^*$ and $a_2{}_* = \R a_2{}_*$ are both left-adjoint and right-adjoint
to each other. Hence we may use $a_2^* = a_2^{\times}$. 
There results a natural isomorphism 
\[
a_1^{\times}a_2^* = a_1^{\times}a_2^{\times} \iso (1_{\X})^{\times} = 
\mathbf{1}.
\]   

We define the first isomorphism in (\dag) via the 
following isomorphisms where 
$\beta_5 \colon a_3^{\times}j_2^* \iso a_2^*a_4^{\times}$ 
is the inverse of the base-change isomorphism associated to 
the fibered square in the preceding diagram.
\[
\Delta^{\!\times}j_2^*q_2^{\times} \iso 
a_1^{\times}a_3^{\times}j_2^*q_2^{\times} \xto{\text{via }\beta_5}
a_1^{\times}a_2^*a_4^{\times}q_2^{\times} \iso
a_4^{\times}q_2^{\times} \iso \mathbf{1}.
\]
The second isomorphism in (\dag) is obtained through a 
similar procedure. 

Thus we have demonstrated an isomorphism 
$c(\sigma_2,\sigma_1) \colon i^{*}p^{\times} \iso j^{*}q^{\times}$. 
In what sense is it canonical?
We discuss some points concerning this.

To begin with, note the usage of fiber-product diagrams 
in the construction. Since fibered products are unique only up to isomorphism,
the objects $\W_{ij}, \V$ used above and the arrows coming in 
and out of them are not unique as functions of the maps that we
started with. Therefore we must 
verify that $c(\sigma_2,\sigma_1)$ is independent 
of the choices of these objects.
This can be done by comparing two different choices of 
fiber-product diagrams via the unique isomorphisms relating them. 

Clearly $c(\sigma_2,\sigma_1)$ is a natural transformation, even one of
triangulated functors.

Finally, if there is a third factorization $\sigma_3$
of $f$, then the cocycle condition, viz., 
$c(\sigma_1,\sigma_2) \circ c(\sigma_2,\sigma_3) = c(\sigma_1,\sigma_3)$
holds. Proving this requires working with triple-fiber-product
diagrams involving the three factorizations.

Since we achieve canonicity 
at the level of complexes, our result is
an improvement upon the one in
\cite[Prop.~2]{AJL2a}, where the isomorphism obtained is canonical
only at the level of the homology modules.

Much of our work is concentrated around generalizing the isomorphism
$c(\sigma_1,\sigma_2)$ and its cocycle property to the situation of 
factorizations of arbitrary
length, i.e., the case when~$\sigma_1$ and~$\sigma_2$ are arbitrarily
long sequences of open and pseudoproper maps. Addressing the concerns
raised above then requires considerable effort. The whole task is
carried out in an abstract setting.  

\[
* \hspace{8em} * \hspace{8em} * 
\]

In \S\ref{sec:abs} we state the abstract results without proof.
The proofs span \S\ref{sec:proofsI}--\S\ref{sec:proofsIV}.
In~\S\ref{sec:apps} we state some applications, most of which   
concern Grothendieck duality.

\subsection{Conventions}
\label{subsec:conv}
\hfill
\begin{enumerate}
\item \label{conv1}
We interchangeably use the terms cartesian square and fibered square.
We call a diagram composed of cartesian squares as a cartesian diagram.
\item \label{conv2}
Reduced notation : Given functors 
\[
\Gamma_1 \colon \B \to \A, \quad \Gamma_2 \colon \C \to \B, \quad 
\Gamma_2'\colon \C \to \B, \quad \Gamma_3 \colon \cD \to \C,
\] 
and a natural transformation 
$\theta \colon \Gamma_2 \to \Gamma_2'$, we call a natural transformation
$\delta \colon \Gamma_1\Gamma_2\Gamma_3 \to \Gamma_1\Gamma_2'\Gamma_3$ 
as ``$\theta$ under reduced notation'' if 
$\delta$ is the obvious transformation induced by $\theta$
and acting as identity on $\Gamma_1,\Gamma_3$.
Here $\Gamma_1,\Gamma_2,\Gamma_2',\Gamma_3$ may occur as
composites of other functors too.
\end{enumerate}

\section{The abstract pasting results}
\label{sec:abs}
Here we present the abstract pasting results without proofs. 
The basic input data
for pasting is stated in \S\ref{subsec:input} below while the output
occurs in Theorem \ref{thm:output4}. In \S\ref{subsec:flat} we discuss
the abstract form of flat-base-change results. The proofs of all these 
results span \S\ref{sec:proofsI}--\S\ref{sec:proofsIV}.

We begin with some preliminaries.

Recall that a \emph{contravariant normalized pseudofunctor $(-)^{\#}$} on a 
category~$\sfC$ assigns to every object $X$ in~$\sfC$ a category $X^{\#}$,
to every map $f\colon X \to Y$ in~$\sfC$ 
a functor $f^{\#} \colon Y^{\#} \to X^{\#}$, and to every pair 
of maps $X \xto{\;f\;} Y \xto{\;g\;} Z$ in~$\sfC$
an isomorphism $C^{\#}_{f,g} \colon f^{\#}g^{\#} \iso (gf)^{\#}$
such that the following conditions hold:
\begin{itemize}
\item $C^{\#}_{-,-}$ is associative vis-\`{a}-vis triple compositions.
\item $(1_Z)^{\#} = \mathbf{1}_{Z^{\#}}$ for all objects $Z$ in $\sfC$.
\item For any $f \colon X \to Y$ in~$\sfC$ the following 
isomorphisms are identity
\[
f^{\#} = 1_X^{\#}f^{\#} \xto{\; C^{\#}_{1_{\<X},f} \;} f^{\#} ,
\qquad
f^{\#} = f^{\#}1_Y^{\#} \xto{\; C^{\#}_{f,1_Y} \;} f^{\#}.
\]
\end{itemize}
Unless mentioned otherwise, any pseudofunctor 
occurring in this paper shall be assumed to be 
contravariant and normalized.

A morphism of pseudofunctors $(-)^{\#} \to (-)^!$ on $\sfC$ consists of the 
following data:
\begin{itemize}
\item For every object $X$ in $\sfC$, there is a functor 
$S_X \colon X^{\#} \to X^!$;
\item For every morphism  $f \colon X \to Y$ in $\sfC$ there is a natural 
transformation $S_Xf^{\#} \to f^!S_Y$;
\end{itemize}
such that for any pair of maps $X \xto{\;f\;} Y \xto{\;g\;} Z$,
the following diagram (of obvious natural maps) commutes,
\[
\begin{CD}
S_Xf^{\#}g^{\#} @>>> f^!S_Yg^{\#} @>>> f^!g^!S_Z \\
@VVV @. @VVV \\
S_X(gf)^{\#} @. \makebox[0pt]{$\xto{\hspace{8em}}$} @. (gf)^!S_Z
\end{CD}
\]
and for any object $X$, the natural map 
$S_X = S_X1_X^{\#} \to 1_X^{!}S_X = S_X$
is the identity. 

The composition of two maps of pseudofunctors $(-)^{\#} \to (-)^! \to (-)^*$
is defined in the obvious way. Note that if $(-)^{\#} \to (-)^!$
is an isomorphism, then the associated functors $S_X$ and 
natural transformations $S_Xf^{\#} \to f^!S_Y$ are isomorphisms. 

\subsection{The input data for gluing}
\label{subsec:input}

Here then is the input data for gluing, consisting of [A]-[D] below. 

\medskip
{\bf [A].} There is a category $\sfC$ 
and there are subcategories $\sfO, \sfP$ that are stable 
under base change by maps in $\sfC$, i.e., 
for any map $f \colon X \to Y$ in $\sfP$ (resp.~$\sfO$)
and any map $g \colon Z \to Y$ in $\sfC$, the fibered product
of $f$ with $g$ exists and 
the induced map $f' \colon X' \to Z$ is also in $\sfP$ (resp.~$\sfO$).  
Moreover we require that the following holds.
\begin{enumerate}
\item The category $\sfI$ of isomorphisms in $\sfC$ is contained 
in~$\sfO$ and in~$\sfP$. (In particular, every object of~$\sfC$
is also in $\sfO$ and $\sfP$.)
\item Let $X \xto{\;f\;} Y \xto{\;g\;} Z$ be morphisms
in~$\sfC$ such that $gf$ is in~$\sfP$. If $g$
is in~$\sfO$ or in~$\sfP$, then $f$ is in~$\sfP$. 
(Therefore, the conclusion that $f$ is in~$\sfP$
also holds if~$g$ is assumed to be a composite of maps 
in $\sfO$ and $\sfP$.)
\end{enumerate}
\medskip 

{\bf [B].} There is a pseudofunctor $(-)^{\times}$ on $\sfP$ 
and a pseudofunctor $(-)^{\ssbox}$ on~$\sfO$ such that 
$X^{\times} = X^{\ssbox}$ for 
any object $X$ in $\sfC$. Henceforth we use $\cD_{X}$ for 
$X^{\times}$ or $X^{\ssbox}$. 


\medskip

{\bf [C].} For any cartesian square $\sfr$ in $\sfC$ as follows
such that $f$ is in $\sfP$ and $i$ is in $\sfO$ 
(so that~$f'$ is in~$\sfP$ and~$i'$ in~$\sfO$), 
\[
\begin{CD}
U @>{i'}>> X \\
@V{f'}VV @VV{f}V \\
V @>{i}>> Y
\end{CD}
\]
there is a ``base-change'' isomorphism
$\beta_{\sfr} \colon i^{'\ssbox}f^{\times} \iso f^{'\times}i^{\ssbox}$,
also sometimes referred to as $\beta_{f,i}$ is there is no cause for 
confusion. Moreover $\beta_{-}$ obeys the transitivity rules~(i) 
and~(ii) stated below.

Consider the following extensions of $\sfr$ where $j \in \sfO$
and $g \in \sfP$. In each case the appended
square is also cartesian and is called $\sfr_1$, while the 
composite cartesian square is called $\cfr$. 
\[
\begin{CD}
U_1 @>{j'}>> U @>{i'}>> X \\
@V{f''}VV @V{f'}VV @V{f}VV \\
V_1 @>{j}>> V @>{i}>> Y
\end{CD}
\hspace{8em}
\begin{CD}
U_1 @>{i''}>> X_1 \\ 
@V{g'}VV @V{g}VV \\
U @>{i'}>> X \\
@V{f'}VV @V{f}VV \\
V @>{i}>> Y
\end{CD}
\] 
\begin{enumerate}
\item Horizontal transitivity: 
The following diagram of isomorphisms commutes. 
\[
\begin{CD}
j'{}^{\ssbox}i'{}^{\ssbox}{f}^{\times} 
 @>{\quad j'{}^{\ssbox}(\beta_{\sfr})\quad}>> 
 j'{}^{\ssbox}f'{}^{\times}{i}^{\ssbox} 
 @>{\quad \beta_{\sfr_1}({i}^{\ssbox}) \quad}>> 
 f''{}^{\times}{j}^{\ssbox}{i}^{\ssbox} \\
@V{C^{\ssbox}_{j',i'}(f^{\times})}VV  @. 
 @VV{f''{}^{\times}(C^{\ssbox}_{j,i})}V \\
(i'j')^{\ssbox}f^{\times} @.
 \makebox[0pt]{$\xto{\hspace{6em} \beta_{\cfr} \hspace{6em}}$}
 @. f''{}^{\times}(ij)^{\ssbox}
\end{CD}
\]
\item Vertical transitivity: The following diagram of isomorphisms commutes. 
\[
\begin{CD}
i''{}^{\ssbox}g^{\times}f^{\times} 
 @>{\quad \beta_{\sfr_1}(f^{\times}) \quad}>> 
 g'{}^{\times}i'{}^{\ssbox}f^{\times} 
 @>{\quad g'{}^{\times}(\beta_{\sfr})\quad}>> 
 g'{}^{\times}f'{}^{\times}i^{\ssbox} \\
@V{i''{}^{\ssbox}(C^{\times}_{g,f})}VV @. 
 @VV{C^{\times}_{g',f'}(i^{\ssbox})}V \\
i''{}^{\ssbox}(gf)^{\times} @.
 \makebox[0pt]{$\xto{\hspace{6em} \beta_{\cfr} \hspace{6em}}$}
 @.  {i}^{\ssbox}(g'f')^{\times} 
\end{CD}
\]
\end{enumerate}

\medskip

{\bf [D].} For any object $X$ in $\sfC$ and for any 
sequence $X \xto{\;f\;} Y \xto{\;g\;} X$ 
factoring the identity map on $X$ such that 
$g \in \sfO$ (so that $f$ is necessarily in~$\sfP$),   
there is a ``fundamental'' isomorphism 
$\phi_{f,g}\colon f^{\times}g^{\ssbox} \iso \oneD{X}$. 
Moreover $\phi_{f,g}$ is compatible with isomorphisms
defined earlier in the following ways (i) and (ii). 
\begin{enumerate}
\item 
Let $h \colon X' \to X$ be a map in~$\sfC$.
Consider the following induced diagram where $g'f' = 1_{X'}$ 
and each square is cartesian,
with the one on the left being called $\lfr$ and the one on right $\rfr$.
\[
\begin{CD}
X' @>{f'}>> Y' @>{g'}>> X' \\
@V{h}VV @V{h'}VV @V{h}VV \\
X @>{f}>> Y @>{g}>> X
\end{CD}
\]
Then the following conditions hold.
\begin{enumerate}
\item If $h$ is in~$\sfP$ 
then the following diagram of isomorphisms commutes.
\[
\begin{CD}
f'{}^{\times}g'{}^{\ssbox}h^{\times} @>{\text{via }\beta_{\rfr}}>> 
 f'{}^{\times}h'{}^{\times}g^{\ssbox} @>{\text{via }(-)^{\times}}>> 
 h^{\times}f^{\times}g^{\ssbox} \\
@V{\phi_{f',g'}(h^{\times})}VV @. @VV{h^{\times}(\phi_{f,g})}V  \\
\oneD{X'}h^{\times} @= h^{\times} @= h^{\times}\oneD{X} 
\end{CD}
\]
\item If $h$ is in~$\sfO$ 
then the following diagram of isomorphisms commutes
where $\lfr^{\text{t}}$ is the transpose of $\lfr$.
\[
\begin{CD}
f'{}^{\times}g'{}^{\ssbox}h^{\ssbox} @>{\text{via }(-)^{\ssbox}}>> 
 f'{}^{\times}h'{}^{\ssbox}g^{\ssbox} 
 @>{\text{via }\beta_{\lfr^{\text{t}}}^{-1}}>> 
 h^{\ssbox}f^{\times}g^{\ssbox} \\
@V{\phi_{f',g'}(h^{\ssbox})}VV @. @VV{h^{\ssbox}(\phi_{f,g})}V  \\
\oneD{X'}h^{\ssbox} @= h^{\ssbox} @= h^{\ssbox}\oneD{X} 
\end{CD}
\]
\end{enumerate}
\item For any diagram as follows where the square $\sfr$ is cartesian,
$f,h,i \in \sfP$, $g,j,k \in \sfO$ and $gf = 1_X = ki$,
\[
\begin{CD}
X \\
@VfVV \\
Y @>g>> X \\
@VhV{\hspace{2em}\sfr}V @VViV \\
Z @>>j> W @>>k> X
\end{CD}
\]
the following diagram of isomorphisms commutes.
\[
\begin{CD}
f^{\times}h^{\times}j^{\ssbox}k^{\ssbox} 
 @>{\qquad\text{via }\beta_{\sfr}\qquad}>> 
 f^{\times}g^{\ssbox}i^{\times}k^{\ssbox} \\
@V{\text{via $C^{\times}_{f,h}$ and $C^{\ssbox}_{j,k}$}}VV 
 @VV{\text{via $\phi_{f,g}$ and $\phi_{i,k}$}}V \\
(hf)^{\times}(kj)^{\ssbox} @>{\qquad\phi_{hf,kj}\qquad}>> \oneD{X} 
\end{CD}
\]
\end{enumerate}

\subsection{The output}
\label{subsec:output}
It will be convenient to develop some terminology before stating 
the output corresponding to the input conditions of~\S\ref{subsec:input}.

Let $\sfQ = \ov{\{\sfO,\sfP\}}$ be the smallest 
subcategory of $\sfC$ containing~$\sfO$ 
and~$\sfP$. Then~$\sfQ$ is also the category, whose objects are all the
objects of~$\sfC$ and whose morphisms are composites of maps 
in~$\sfO$ and~$\sfP$.

\begin{adefi}
\label{def:output1}
A \emph{pseudofunctorial cover} 
$\eC$ (with respect to data [A]--[D]) is a triple
$((-)^{!}, \mu^{!}_{\times},\mu^{!}_{\ssbox})$ consisting of 
\begin{enumerate}
\item a pseudofunctor $(-)^!$ on $\sfQ$;
\item a pseudofunctorial morphism 
$\mu^{!}_{\times} \colon (-)^!\big|_{\sfP} \xto{\quad} (-)^{\times}$ 
on $\sfP$;
\item a pseudofunctorial morphism 
$\mu^{!}_{\ssbox} \colon (-)^!\big|_{\sfO} \xto{\quad} (-)^{\ssbox}$ 
on $\sfO$;
\end{enumerate}
subject to the following conditions (a), (b) and (c).
\begin{enumerate}
\renewcommand{\theenumi}{\alph{enumi}}
\item For any object $X \in \sfC$, the functors 
$X^! \to X^{\times} = \cD_X$ and $X^! \to X^{\ssbox} = \cD_X$
induced by $\mu^{!}_{\times}$ and $\mu^{!}_{\ssbox}$ respectively,
are the same. Let us denote this functor $X^! \to \cD_X$ by~$S_X$.
\item For every cartesian square
$\sfr$ in $\sfC$ as follows
such that $f \in \sfP$ and $i \in \sfO$, 
\[
\begin{CD}
U @>{i'}>> X \\
@V{f'}VV @VV{f}V \\
V @>{i}>> Y
\end{CD}
\]
the following diagram commutes.
\[
\begin{CD}
S_Ui'{}^{!}f^{!} @>{\quad\text{via } (-)^!\quad}>> S_Uf'{}^{!}i^{!} \\
@V{\text{via (iii)}}VV @VV{\text{via (ii)}}V \\
i'{}^{\ssbox}S_Xf^{!} @. f'{}^{\times}S_Vi^{!} \\
@V{\text{via (ii)}}VV @VV{\text{via (iii)}}V \\
i'{}^{\ssbox}f^{\times}S_Y @>{\quad\text{via }\beta_{\sfr}\quad}>> 
 f'{}^{\times}i^{\ssbox}S_Y 
\end{CD}
\]
\item For every sequence $X \xto{i} Y \xto{h} X$ such that 
$hi = 1_X, i \in \sfP, h \in \sfO$, the following diagram
commutes.
\[
\begin{CD}
S_Xi^!h^! @>{\quad\text{via }(-)^!\quad}>> S_X\mathbf{1}_{X^!} \\
@V{\text{via (ii)}}VV @| \\
i^{\times}S_Yh^! @. S_X \\
@V{\text{via (iii)}}VV @| \\
i^{\times}h^{\ssbox}S_X @>{\quad\text{via }\phi_{i,h}\quad}>> 
 \mathbf{1}_{\cD_X}S_X
\end{CD}
\]
\end{enumerate}
\end{adefi}

\ssss
\label{sssec:output2}
A morphism between two pseudofunctorial covers is defined as follows. 
Let $\eC_1 = ((-)^{!}, \mu^{!}_{\times},\mu^{!}_{\ssbox}),
\eC_2 = ((-)^{\#}, \mu^{\#}_{\times},\mu^{\#}_{\ssbox})$ be 
two pseudofunctorial covers. Then a morphism $\eC_2 \to \eC_1$ 
of pseudofunctorial covers is a morphism of 
pseudofunctors $\epsilon \colon (-)^{\#} \to (-)^{!}$
such that the following two diagrams commute.
\[
\xymatrix{
(-)^{\#}\big|_{\sfP} \ar[rd] \ar[rr]^{\epsilon|_{\sfP}} 
 & & (-)^!\big|_{\sfP} \ar[ld]^{\cong} \\
& (-)^{\times}
}
\qquad\qquad
\xymatrix{
(-)^{\#}\big|_{\sfO} \ar[rd] \ar[rr]^{\epsilon|_{\sfP}}  
 & & (-)^!\big|_{\sfO} \ar[ld]^{\cong} \\
& (-)^{\ssbox}
}
\]

\begin{adefi}
\label{def:output3} 
A pseudofunctorial cover 
$\eC = ((-)^{!}, \mu^{!}_{\times},\mu^{!}_{\ssbox})$ is 
called perfect if $\mu^{!}_{\times}$, $\mu^{!}_{\ssbox}$ 
are isomorphisms.
\end{adefi}

The output for the input data of \S\ref{subsec:input} may now be 
stated as follows. For convenience we call a pseudofunctorial cover
simply as a cover.

\begin{athm}
\label{thm:output4}
Under input conditions \textup{[A]--[D]} of \S\textup{\ref{subsec:input}}
the following hold.
\begin{enumerate}
\item There exists a perfect cover.
\item Any perfect cover is final in the category of 
all covers, i.e., for any perfect 
cover $\eC$ and any cover $\eC'$ there exists a 
unique map of covers $\eC' \to \eC$. In particular, 
any two perfect covers are isomorphic via a unique isomorphism.
\end{enumerate}
\end{athm}

A proof of this theorem is given in \S\ref{subsec:Thm1} and is 
based on the results of \S\ref{sec:proofsI}--\ref{sec:proofsIII}.

\begin{arem}
\label{rem:triangle}
Suppose we further assume for the input conditions that 
for any object $X$ in $\sfC$, $\cD_X$ is a triangulated category
and that the functors $f^{\times}, f^{\ssbox}$ for~$f$ in~$\sfP$ or~$\sfO$
are triangulated functors that ``commute'' with translation and that 
$C^{\times}_{-,-}$, $C^{\ssbox}_{-,-}$, $\beta_{-}$ and $\phi_{-}$
are morphisms of triangulated functors. Then the perfect 
cover of Theorem~\ref{thm:output4}, say
$\eC = ((-)^{!}, \mu^{!}_{\times},\mu^{!}_{\ssbox})$,
can be chosen such that for any~$f$ in~$\sfQ$, $f^!$ is a triangulated
functor that commutes with translation and such that the isomorphisms
$C^{!}_{-,-},\mu^{!}_{\times},\mu^{!}_{\ssbox}$ respect the 
appropriate triangulated
structures. This will follow from the proof of~\ref{thm:output4}. 
\end{arem}

\begin{arem}
\label{rem:intersection}
A curious consequence of Theorem \ref{thm:output4}(i) is that on the 
intersection $\sfR$ of $\sfP$ and $\sfO$, the restrictions 
of $(-)^{\times}$ and $(-)^{\ssbox}$ to $\sfR$ are isomorphic
as pseudofunctors. This is clearly a necessary condition for pasting
but we do not find it useful to have it as part of the input data.
In particular, it is not a substitute for the more powerful condition [D].
(cf.~\ref{thm:twisted6} where [D] is easy to attain.)

For now, we give a quick indication of how to construct, 
for any map $f\colon X \to Y$ in $\sfR$, a natural
isomorphism $f^{\ssbox} \iso f^{\times}$.
Consider the following diagram where the square $\sfr$ is cartesian
and $\Delta$ is the diagonal map, which, by \S\ref{subsec:input}[A](i),(ii),
is in $\sfP$.
\[
\begin{CD}  
X @>{\Delta \in \sfP}>> X^2 @>{\pi_2 \in \sfO}>> X \\
@. @V{\pi_1 \in \sfP}V{\hspace{2em} \sfr}V @VV{f \in \sfP}V \\
@. X @>>{f \in \sfO}> Y
\end{CD}
\]
Then the desired isomorphism is obtained as
\[
f^{\ssbox} \xto{\text{via }(1_X)^{\times} \cong \Delta^{\times}\pi_1^{\times}}
\Delta^{\times}\pi_1^{\times}f^{\ssbox} \xto{\text{via }\beta_{\sfr}} 
\Delta^{\times}\pi_2^{\ssbox}f^{\times} \xto{\text{via }\phi_{\Delta,\pi_2}}
f^{\times}.
\] 
Thus for obtaining this isomorphism alone, the compatibilities 
in [C], [D] of the input data are not needed. This isomorphism 
is generalized further in \ref{sssec:psi7} below. 
Verdier's proof of the dualizing
property of differentials (\cite[Theorem 3]{Ve}) can be seen 
to be an instance of this by using $\sfO =$ smooth maps, etc..
\end{arem}
 
\subsection{Flat base change}
\label{subsec:flat}
Now we consider change of base by maps in a new 
subcategory of $\sfC$ and a corresponding new class of 
base-change isomorphisms.
We distinguish it from the base-change isomorphisms 
of~\S\ref{subsec:input}[C], by calling it \emph{flat} base change.
 
\ssss
\label{sssec:flat1}
Consider the following addition [E1]--[E3] to the input data [A]--[D] 
of~\S\ref{subsec:input}.

\medskip
{\bf [E1].} There is a subcategory $\sfF$ of $\sfC$ that is stable under 
base change by maps in~$\sfO$ or $\sfP$ and contains the 
subcategory $\sfI$ of isomorphisms in $\sfC$. 

{\bf [E2].} On $\sfF$ there is a pseudofunctor $(-)^{\flat}$ such that 
$X^{\flat} = \cD_X$ for any object~$X$ in~$\sfC$. 

{\bf [E3].}
Let $\sfr$ be a cartesian square as follows
where $f$ is in $\sfO$ or $\sfP$ and $u$ is in $\sfF$. 

\[
\begin{CD}
W @>{u'}>> X \\
@V{f'}VV @VV{f}V \\
Z @>{u}>> Y
\end{CD}
\]
Then the following conditions (i) and (ii) hold subject to (a) and (b) below.
\begin{enumerate}
\item If $f$ is in $\sfP$, then there is an isomorphism 
$\sbeta_{\sfr}^{\times} \colon 
u'{}^{\flat}f^{\times} \iso f'{}^{\times}u^{\flat}$,
also sometimes written as $\sbeta_{f,u}^{\times}$ if there is no cause
for confusion. Moreover $\sbeta_{-}^{\times}$ 
is transitive vis-\`{a}-vis extensions of $\sfr$ via 
$\sfF$-maps horizontally or $\sfP$-maps vertically.  
\item If $f$ is in $\sfO$, then there is an isomorphism 
$\sbeta_{\sfr}^{\ssbox} \colon  
u'{}^{\flat}f^{\ssbox} \iso f'{}^{\ssbox}u^{\flat}$, 
also sometimes written as $\sbeta_{f,u}^{\ssbox}$ if there is no cause
for confusion. Moreover $\sbeta_{-}^{\ssbox}$ is transitive 
vis-\`{a}-vis extensions of $\sfr$ via 
$\sfF$-maps horizontally or $\sfO$-maps vertically.  
\end{enumerate}

The ``flat-base-change'' 
isomorphisms $\sbeta^{\times}_{-}, \sbeta^{\ssbox}_{-}$
are compatible with the isomorphisms $\beta_{-}$ and $\phi_{-,-}$ of
\S\ref{subsec:input}[C] and \S\ref{subsec:input}[D] 
as follows.
\begin{enumerate}
\renewcommand{\theenumi}{\alph{enumi}}
\item For any cartesian ``cube'' as follows, (i.e., each ``face''
of the cube is a fibered product diagram)
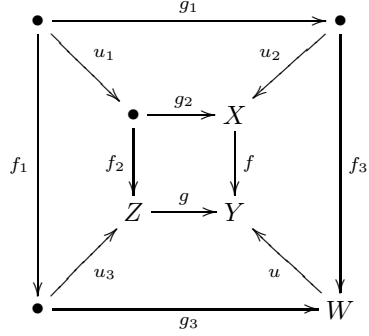
\[
\xymatrix{
\bullet \ar[ddd]_{f_1} \ar[dr]^{u_1} \ar[rrr]^{g_1} 
 & & & \bullet \ar[ddd]^{f_3} \ar[dl]_{u_2} \\
& \bullet \ar[d]_{f_2} \ar[r]^{g_2} & X \ar[d]^{f} & \\
& Z \ar[r]^{g} & Y  & \\
\bullet \ar[ur]_{u_3} \ar[rrr]_{g_3} & & & W  \ar[ul]^{u}
}
\]
if $f \in \sfP, g \in \sfO, u \in \sfF$, then 
the following hexagon of isomorphisms commutes.
\[
\xymatrix{
& u_1^{\flat}g_2^{\ssbox}f^{\times} \ar[rr]^{\text{via }}_{\beta_{f,g}} 
 & & u_1^{\flat}f_2^{\times}g^{\ssbox}
 \ar[dr]^{\text{via }\sbeta_{f_2,u_3}^{\times}}  \\
g_1^{\ssbox}u_2^{\flat}f^{\times} 
 \ar[ur]^{\text{via }(\sbeta_{g_2,u_2}^{\ssbox})^{-1}} 
 \ar[dr]_{\text{via }\sbeta_{f,u}^{\times}} 
 & & & & f_1^{\times}u_3^{\flat}g^{\ssbox} \\
& g_1^{\ssbox}f_3^{\times}u^{\flat} 
 \ar[rr]^{\text{via }}_{\beta_{f_3,g_3}} & & f_1^{\times}g_3^{\ssbox}u^{\flat} 
 \ar[ur]_{\text{via }(\sbeta_{g,u}^{\ssbox})^{-1}}
}
\]
\item For any cartesian diagram as follows where $gf=1_X$, $f \in \sfP, 
g \in \sfO, u \in \sfF$,
\[
\begin{CD}
X' @>{f'}>> Y' @>{g'}>> X' \\
@V{u}VV @V{u'}VV @V{u}VV \\
X @>{f}>> Y @>{g}>> X
\end{CD}
\]
the following diagram of isomorphisms commutes.
\[
\begin{CD}
f'{}^{\times}g'{}^{\ssbox}u^{\flat} 
 @>{\text{via }}>{(\sbeta_{g,u}^{\ssbox})^{-1}}> 
 f'{}^{\times}u'{}^{\flat}g^{\ssbox}
 @>{\text{via }}>{(\sbeta_{f,u'}^{\ssbox})^{-1}}> 
 u^{\flat}f^{\times}g^{\ssbox} \\
@V{\text{via }\phi_{f',g'}}VV @. @VV{\text{via }\phi_{f,g}}V \\
\mathbf{1}_{\cD_{X'}}u^{\flat} @= u^{\flat} @= u^{\flat}\mathbf{1}_{\cD_{X}}
\end{CD}
\]
\end{enumerate}

Under the conditions [A]--[D], [E1]--[E3], we obtain the following output.

\begin{athm}
\label{thm:flat5}
Let $\eC = ((-)^{!}, \mu^{!}_{\times},\mu^{!}_{\ssbox})$ be a 
perfect pseudofunctorial
cover. Then there is a unique family of flat-base-change isomorphisms 
$\sbeta^!_{-}$, (i.e.,
for any cartesian square $\sfr$ as follows 
where $f \in \sfQ = \ov{\{\sfO,\sfP\}}$ and $u \in \sfF$, 
\[
\begin{CD}
W @>{u'}>> X \\
@V{f'}VV @VV{f}V \\
Z @>{u}>> Y
\end{CD}
\]
there is an isomorphism
$\sbeta^!_{\sfr} \colon u'{}^{\flat}f^! \iso f'{}^!u^{\flat}$)
such that the following hold.
\begin{enumerate}
\item $\sbeta^!_{-}$ is vertically and horizontally transitive. 
\item For $\sfr$ as above, if $f$ is in $\sfP$ (resp.~in $\sfO$),
then the following diagram on the left (resp.~on the right) commutes.
\[
\begin{CD}
u'{}^{\flat}f^! @>{\sbeta_{\sfr}^!}>> f'{}^!u^{\flat} \\
@V{\textup{via }\mu^!_{\times}}VV @VV{\textup{via }\mu^!_{\times}}V \\
u'{}^{\flat}f^{\times} @>>{\sbeta_{\sfr}^{\times}}> f'{}^{\times}u^{\flat} 
\end{CD}
\qquad \qquad \qquad 
\begin{CD}
u'{}^{\flat}f^! @>{\sbeta_{\sfr}^!}>> f'{}^!u^{\flat} \\
@V{\textup{via }\mu^!_{\ssbox}}VV @VV{\textup{via }\mu^!_{\ssbox}}V \\
u'{}^{\flat}f^{\ssbox} @>>{\sbeta_{\sfr}^{\ssbox}}> f'{}^{\ssbox}u^{\flat} 
\end{CD}
\] 
\end{enumerate}
In particular, if $\eC' = ((-)^{\#}, \mu^{\#}_{\times},\mu^{\#}_{\ssbox})$ 
is another perfect cover with $\sbeta^{\#}_{-}$ the corresponding family of
flat-base-change isomorphisms, 
then for~$\sfr$ as above, via the unique isomorphism of covers $\eC' \to \eC$,
the following diagram commutes.
\[
\begin{CD}
u'{}^{\flat}f^{\#} @>{\sbeta^{\#}_{\sfr}}>> f'{}^{\#}u^{\flat} \\
@VVV @VVV \\
u'{}^{\flat}f^{!} @>{\sbeta^{!}_{\sfr}}>> f'{}^{!}u^{\flat} 
\end{CD}
\]
\end{athm}

A proof of this theorem is given in \S\ref{subsec:Thm2}.

\subsection{A variant}
\label{subsec:variant}
We give an alternate criterion for pasting pseudofunctors 
and the related base-change result. 

\ssss
\label{sssec:variant1}
For input data, we assume that the following conditions (1)--(5) hold.
\begin{enumerate}
\renewcommand{\theenumi}{\arabic{enumi}}
\item Conditions [A], [B] of \S\ref{subsec:input} 
and [E1], [E2], [E3](i) of \S\ref{subsec:flat} hold.
\item The subcategory $\sfO$ is contained in~$\sfF$ and $(-)^{\ssbox}$
is the restriction of~$(-)^{\flat}$ to~$\sfO$. For a fibered square $\sfr$
resulting from the change of base of a $\sfP$-map 
by an $\sfO$-map, 
we set $\beta_{\sfr} \set \sbeta_{\sfr}^{\times}$.
\item If $i$ is an isomorphism in $\sfC$, then $i^{\ssbox} = i^{\times}$. 
\item If $X \xto{f} Y \xto{g} X$ are $\sfC$-maps such that 
$gf=1_X, f \in \sfP, g \in \sfO$, then $f,g$ are isomorphisms.
\item Let $\sfr$ be a fibered square as follows where $f,g \in \sfP$ and
$i,j \in \sfF$.
\[
\begin{CD}
U @>{j}>> X \\
@V{g}VV @VV{f}V \\
V @>{i}>> Y
\end{CD}
\] 
If $i,j$ are isomorphisms, then among the following diagrams,
the left one commutes 
while if $f,g$ are isomorphisms, then the right one commutes.
\[
\begin{CD}
j^{\ssbox}f^{\times} @>{\beta_{\sfr}}>> g^{\times}i^{\ssbox} \\
@| @| \\
j^{\times}f^{\times} @>{\textup{via } (-)^{\times}}>> g^{\times}i^{\times} 
\end{CD}
\qquad
\qquad
\begin{CD}
j^{\flat}f^{\times} @>{\sbeta_{\sfr}^{\times}}>> g^{\times}i^{\flat} \\
@| @| \\
j^{\flat}f^{\flat} @>{\textup{via } (-)^{\flat}}>> g^{\flat}i^{\flat} 
\end{CD}
\]
\end{enumerate}

\ssss 
\label{sssec:variant2}
We define a pseudofunctorial cover $\eC$ 
(with respect to data (1)--(5))
the same way as in \ref{def:output1},
except for the following modification.
We do not require \ref{def:output1}(c) to hold, but instead require that
for every isomorphism $i \colon X \to Y$ in $\sfC$, the natural maps
$S_Xi^! \to i^{\times}S_Y$ and $S_Xi^! \to i^{\ssbox}S_Y$ are the 
same via $i^{\times}=i^{\ssbox}$.

The notion of a morphism of covers and that of a
perfect cover is the same as before.

For a fibered square $\sfr$
resulting from the change of base of an $\sfO$-map 
by an $\sfF$-map, 
we set $\sbeta_{\sfr}^{\ssbox}$ to be the obvious map 
resulting from pseudofunctoriality of $(-)^{\flat}$.

\begin{athm}
\label{thm:variant3}
Under the input conditions \textup{(1)--(5)}, the following hold.
\begin{enumerate}
\item There exists a perfect cover. Any such cover is final 
in the category of all covers.
\item For any perfect cover, there is a unique family of flat-base-change
isomorphisms such that the conclusion of \textup{\ref{thm:flat5}} holds.
\end{enumerate}
\end{athm}

We give a proof of this theorem in \S\ref{subsec:Thm3}.

\section{Proofs I (generalized isomorphisms in the labeled setup)}
\label{sec:proofsI}
The proofs of the output statements of \S\ref{subsec:output}, 
\S\ref{subsec:flat} and~\S\ref{subsec:variant} 
are based on the results of the next few sections.
The general theme underlying these proofs 
is of working with sequences of (composable) $\sfC$-maps where each map 
is in~$\sfO$ or in~$\sfP$. After proving results at the level 
of such sequences, one then descends to statements concerning maps
in $\sfQ = \ov{\{\sfO,\sfP\}}$. 

In this section we introduce the setup of labeled maps and sequences, 
and corresponding diagrams. This is the setup in which the results 
of \S\ref{sec:proofsII} and \S\ref{sec:proofsIII} are stated and proved.
The notion of labeled maps provides us with a uniform way of 
treating similar results through common notation, thus making the transition
to working with sequences easier.
Most of our work in this section goes into 
showing that the base-change isomorphisms and the fundamental isomorphisms  
of the input data of \S\ref{subsec:input} 
can be upgraded to the level of labeled sequences. 


\subsection{Labeled maps and sequences}
\label{subsec:lamase}

For a map $f$ which is both, in $\sfO$ and in~$\sfP$, the input conditions
of \S\ref{subsec:input} provide us with two immediate choices for a 
functor, viz., $f^{\ssbox}$ and~$f^{\times}$. This motivates the 
following definitions as a way of keeping track of the choice 
whenever necessary.
 
\begin{adefi}
\label{def:lamase1}
A \emph{labeled map} in $\sfC$ is a pair $F \set (f,\lambda)$
where $f$ is a map in $\sfO$ or~$\sfP$ and $\lambda$, the label, 
is an element of the set $\{ \sfO,\sfP \}$ such
that $f \in \lambda$.
\end{adefi}

\begin{adefi}
\label{def:lamase2}
For a labeled map $X \xto{\;F \;=\; (f,\lambda)\;} Y$ in $\sfC$ 
we define a functor from $\cD_Y$ to $\cD_X$ by 
\[
F^{\bxt} = (f,\lambda)^{\bxt} \set 
 \begin{cases} 
 f^{\ssbox}, &\text{if $\lambda = \sfO$;}  \\
 f^{\times}, &\text{if $\lambda = \sfP$.} 
 \end{cases}
\]
\end{adefi}

\begin{adefi}
\label{def:lamase2a}
By a labeled sequence $\sigma \set F_1,\ldots,F_n$ 
we mean a sequence of labeled maps
\[
X_1 \xto{\;F_1\;} X_2 \xto{\;F_2\;} X_3 
\xto{\quad} \cdots \xto{\quad}
X_{n} \xto{\;F_{n}\;} X_{n+1}.
\]
Suppose $F_i = (f_i,\lambda_i)$.
We define $|\sigma|$ to be the composite map 
$f_{n} \cdots f_1 \colon X_1 \to X_{n+1}$. 
We call $X_1$ the source of $\sigma$, $X_{n+1}$ its target
and $n$ its length. We frequently also use double-arrow
notation such as $\sigma \colon X_1 \Lra X_{n+1}$ to 
denote a labeled sequence.
\end{adefi}

\begin{adefi}
\label{def:lamase3}
For any labeled sequence $\sigma = F_1,\ldots,F_n$ 
we define the functor ${\sigma}^{\bxt} \colon \cD_{X_{n+1}} \to \cD_{X_1}$ by 
${\sigma}^{\bxt} \set F_1^{\bxt}F_2^{\bxt} \cdots F_n^{\bxt}$.
\end{adefi}

One of the main constructions used in the proofs of the output
statements in~\S\ref{subsec:output} and~\S\ref{subsec:flat}
is that of a canonical isomorphism $\Psi_{\sigma_1, \sigma_2} 
\colon \sigma_2^{\bxt} \iso \sigma_1^{\bxt}$ associated to any 
pair of labeled sequences $\sigma_1, \sigma_2$ such that 
$|\sigma_1| = |\sigma_2|$. This construction 
is carried out in \S\ref{sec:proofsII}. For the rest of this section
we concentrate on developing the tools used in its construction. 

\subsection{Compositions and fibered products of labeled maps}
\label{subsec:cafpols}

Here we give a collection of basic notions
involving labeled maps and sequences. Mainly, we
extend the isomorphisms of the input data in~\S\ref{subsec:input} 
and their compatibilities to the \emph{labeled} setup. Thus the isomorphisms
defined here are not really new, but are simply repackaged versions
of the old ones designed to work in the new setup.

\ssss
\label{sssec:cafpols-1}
Let $\sigma_1 = F_1, \ldots, F_{n}$
and $\sigma_2 = G_1, \ldots, G_{m}$ be labeled sequences such 
that the target of $\sigma_1$ equals the source of $\sigma_2$.
We may then concatenate
the two sequences in the obvious manner resulting in a labeled sequence 
\[
\sigma_1 \star \sigma_2 \set F_1, \ldots, F_{n},G_1, \ldots, G_{m}.
\] 
Note that $|\sigma_1 \star \sigma_2| = |\sigma_2||\sigma_1|$
and $(\sigma_1 \star \sigma_2)^{\bxt} = \sigma_2^{\bxt}\sigma_1^{\bxt}$.

\medskip
\ssss
\label{sssec:cafpols0}
For \emph{equi}-labeled maps $X \xto{(f_1,\lambda)} Y \xto{(f_2,\lambda)} Z$, 
with $f_3 \set f_2f_1$ and $F_i \set (f_i,\lambda)$ we define a 
comparison isomorphism 
$\boldC_{F_1, F_2} \colon F_1^{\bxt}F_2^{\bxt} \iso F_3^{\bxt}$
by
\[
\boldC_{F_1, F_2} = 
\begin{cases}
C^{\times}_{f_1,f_2} \vspace{2pt} &\text{if $\lambda_i = \sfP$;} \\  
C^{\ssbox}_{f_1,f_2} &\text{if $\lambda_i = \sfO$.}   
\end{cases}
\]

\medskip
\ssss
\label{sssec:cafpols3z}
Let $X \xto{(f_1,\lambda_1)} Y \xto{(f_2,\lambda_2)} X$ be 
labeled maps such that $\lambda_1 = \sfP$ and $f_2f_1 = 1_X$.
Set $F_i \set (f_i,\lambda_i)$. 
Then we define a fundamental isomorphism
\begin{equation}
\label{eq:cafpols3}
\boldsymbol{\phi}_{F_1, F_2} \colon 
F_1^{\bxt}F_2^{\bxt} \iso \oneD{X}
\end{equation}
by
\[
\boldsymbol{\phi}_{F_1, F_2} = 
\begin{cases}
C^{\times}_{f_1,f_2}
 &\text{if $\lambda_2 = \sfP$;} \\  
\phi_{f_1, f_2} 
 &\text{if $\lambda_2 = \sfO$.}  
\end{cases}
\]

\medskip
\ssss
\label{sssec:cafpols1a}
A \emph{labeled cartesian square} is a 
quadruplet $\sfr = (F_1, F_2, F_1',F_2')$
of labeled maps with
\[
F_1 = (f_1,\lambda_1), \quad F_2 = (f_2,\lambda_2), \quad
F_1' = (f_1',\lambda_1'), \quad F_2' = (f_2',\lambda_2'), 
\]
that fit into a diagram as follows
\begin{equation}
\label{eq:cafpols1}
\begin{CD}
U @>{F_2'}>> X \\ 
@V{F_1'}VV @VV{F_1}V \\ 
V @>{F_2}>> Y
\end{CD}
\end{equation}
such that the underlying diagram of $\sfC$-maps is a usual
cartesian square and $\lambda_i' = \lambda_i$.
We shall use the same name (such as $\sfr$) to denote a labeled
cartesian square and the underlying square of unlabeled maps.
We also say that $\sfr$ is a labeled cartesian square on $(F_1,F_2)$.
By default, the first and third element of the quadruplet underlying $\sfr$
are drawn vertically and the other two horizontally. 
Via this convention, the diagram uniquely determines the quadruplet. 

The \emph{transpose} of a labeled 
cartesian square $\sfr = (F_1, F_2, F_1',F_2')$
is the labeled cartesian square $\sfr^{\text{t}} \set (F_2, F_1, F_2',F_1')$.

Any two labeled maps $F_1,F_2$ having the same target, give rise to a 
labeled cartesian square on $(F_1,F_2)$, in the obvious way. Of course, 
the new ingredients of the square are not uniquely determined 
though they are unique up to a unique isomorphism.

\emph{From now on, we shall implicitly assume all maps and cartesian
squares to be labeled ones.}

For the cartesian square $\sfr$ in \eqref{eq:cafpols1},
we define a base-change isomorphism
\begin{equation}
\label{eq:cafpols2}
\boldsymbol{\beta}_{\sfr} \colon 
F_2^{'\bxt}F_1^{\bxt} \iso F_1^{'\bxt}F_2^{\bxt}
\end{equation}
by
\[
\boldsymbol{\beta}_{\sfr} = 
\begin{cases}
(C_{f_1',f_2}^{\times})^{-1}C_{f_2',f_1}^{\times},\vspace{3pt}  
 &\text{if $\lambda_1= \lambda_2 = \sfP$;} \\
(C_{f_1',f_2}^{\ssbox})^{-1}C_{f_2',f_1}^{\ssbox},\vspace{3pt}  
 &\text{if $\lambda_1= \lambda_2 = \sfO$;} \\
\beta_{\sfr}, 
 &\text{if $\lambda_1 = \sfP, \lambda_2 = \sfO$;} \\
\beta_{\sfr^{\text{t}}}^{-1}, 
 &\text{if $\lambda_1 = \sfO, \lambda_2 = \sfP$.} 
\end{cases}
\]
Note that $\bbeta_{-}$ is symmetric, i.e., 
$\bbeta_{\sfr} = \bbeta_{\sfr^{\text{t}}}^{-1}$.

The isomorphisms $\boldC_{-,-}$, $\bphi_{-,-}$ and $\bbeta_{-}$
defined above  
satisfy many compatibilities including those analogous to the ones 
in \S\ref{subsec:input}. We record some of these.

\begin{alem}
\label{lem:cafpols4}
\hfill
\begin{enumerate}
\item Let $\sfr$ be a labeled cartesian square as follows
\[
\begin{CD}
X' @>F>> X \\
@VI'VV @VIVV \\
X' @>F>> X
\end{CD}
\] 
where $I,I'$ are identity maps with label $\sfP$. Then the natural 
isomorphism 
\[
 F^{\bxt} = F^{\bxt}I^{\bxt} \xto{\;\bbeta_{\sfr}\;}
I'{}^{\bxt}F^{\bxt} = F^{\bxt} 
\]
is the identity.
\item The fundamental isomorphism $\boldsymbol{\phi}_{F_1, F_2}$
of \eqref{eq:cafpols3} \vspace{-2pt} is compatible with base change on~$X$ 
by a labeled map, i.e., for any diagram of labeled cartesian squares
as follows, where the top row composes to $1_{X'}$, 
$\sfr_1$ is the square on the left and $\sfr_2$ the one on the right, 
\[
\begin{CD}
X' @>{F_1'}>> Y' @>{F_2'}>> X' \\
@V{F_3}VV @V{F_3'}VV @V{F_3}VV \\
X @>>{F_1}> Y @>>{F_2}> X
\end{CD}
\]
the following diagram of isomorphisms commutes. 
\[
\begin{CD}
F_1^{'\bxt}F_2^{'\bxt}F_3^{\bxt} @>{F_1^{'\bxt}(\bbeta_{\sfr_2})}>>  
 F_1^{'\bxt}F_3^{'\bxt}F_2^{\bxt} @>{\bbeta_{\sfr_1}(F_2^{\bxt})}>>  
 F_3^{\bxt}F_1^{\bxt}F_2^{\bxt} \\
@VV{\bphi_{F_1',F_2'}(F_3^{\bxt})}V @. @V{F_3^{\bxt}(\bphi_{F_1,F_2})}VV \\
\oneD{{X'}}F_3^{\bxt} @= F_3^{\bxt} @=  
 F_3^{\bxt}\oneD{{X}}
\end{CD}
\]
\item The base-change isomorphism 
of~\eqref{eq:cafpols2} is horizontally  and vertically transitive.
For instance, consider a
horizontal extension to the diagram~$\sfr$ in~\eqref{eq:cafpols1} 
as follows where $F_3$ has the same label $\lambda_2$ as $F_2$
and the appended square, is also cartesian.
\[
\begin{CD}
U_1 @>{F_3'}>> U @>{F_2'}>> X \\
@V{F_1''}VV @V{F_1'}VV @VV{F_1}V \\
V_1 @>>{F_3}> V @>>{F_2}> Y
\end{CD}
\]
Let $F_4, F_4'$ be the composite maps $|F_2||F_3|, |F_2'||F_3'|$
respectively with label~$\lambda_2$.
Let $\sfr_1$ be the appended square on the left and $\cfr$ the 
composite square. Then the following diagram commutes, 
where 
\[
\begin{CD}
F_3^{'\bxt}F_2^{'\bxt}F_1^{\bxt} @>>{F_3^{'\bxt}(\bbeta_{\sfr})}> 
 F_3^{'\bxt}F_1^{'\bxt}F_2^{\bxt} @>>{\bbeta_{\sfr_1}(F_2^{\bxt})}> 
 F_1^{''\bxt}F_3^{\bxt}F_2^{\bxt} \\
@V{\boldC_{F_3',F_2'}(F_1^{\bxt})}VV @. 
 @VV{F_1^{''\bxt}(\boldC_{F_3,F_2})}V \\
F_4^{'\bxt}F_1^{\bxt} @. 
 \makebox[0pt]{$\xto{\hspace{5em}\bbeta_{\cfr}\hspace{5em}}$}
 @. F_1^{''\bxt}F_4^{\bxt}
\end{CD}
\]
Similarly, transitivity holds for vertical extensions of $\sfr$.
\end{enumerate}
\end{alem}

\begin{proof}
(i). Let $F = (f, \lambda)$. If $\lambda = \sfP$ then we 
conclude using pseudofunctoriality of $(-)^{\times}$. Assume
$\lambda = \sfO$. 
Let $b$ denote the automorphism of~$F^{\bxt}$ induced by~$\bbeta_{\sfr}$.
From the vertical transitivity property of~$\beta_{-}$
(\S\ref{subsec:input}, [C](ii))
corresponding to the following diagram,
\[
\begin{CD}
X' @>f>> X \\ @| @| \\ X' @>f>> X \\ @| @| \\ X' @>f>> X 
\end{CD}
\]
we conclude that $bb = b$. Since $b$ is an isomorphism, therefore it
is the identity. 

(ii).
Recall that $\lambda_1 = \sfP$ by hypothesis in \eqref{eq:cafpols3}.
If $\lambda_2 = \sfO$, then the asserted compatibility follows from 
parts (a) or (b) of \S\ref{subsec:input}[D](i), depending on 
whether $\lambda_3$ equals $\sfP$ or $\sfO$. 
Now assume $\lambda_2 = \sfP$. If $\lambda_3 = \sfO$, then
we conclude by the vertical transitivity property 
of \S\ref{subsec:input}[C](i) and by part (i) above. 
If $\lambda_3 = \sfP$, then all the labels are $\sfP$ 
and we conclude by pseudofunctoriality. 

(iii). We argue as in (ii) by looking at various cases depending on the
value of the labels involved.
\end{proof}

\subsection{Fibered products of sequences}
\label{subsec:fpos}

We consider diagrams resulting from taking fibered products of
two sequences having the same target and associate a generalized 
base-change isomorphism to such diagrams.


\ssss
\label{sssec:cafpols5}
Let $\sigma_1,\sigma_2$ be two labeled sequences having the same target $T$.
Then one associates a cartesian diagram to these two sequences as 
follows. Suppose $\sigma_1, \sigma_2$ are given, respectively, by 
\begin{align}
&Y_1 \xto{\;(f_1,\alpha_1)\;} Y_2 \xto{\;(f_2,\alpha_2)\;} Y_3 
\xto{\quad} \cdots \xto{\quad}
Y_{n} \xto{\;(f_{n},\alpha_{n})\;} Y_{n+1} = T, \notag \\
&Z_1 \xto{\;(g_1,\gamma_1)\;} Z_2 \xto{\;(g_2,\gamma_2)\;} Z_3 
\xto{\quad} \cdots \xto{\quad}
Z_{m} \xto{\;(g_{m},\gamma_{m})\;} Z_{m+1} = T. \notag 
\end{align}
Then a \emph{fiber-product diagram on $\sigma_1, \sigma_2$} (or 
\emph{cartesian diagram on $\sigma_1, \sigma_2$})
is an $n \times m$ array $\sS = [\sfr_{i,j}]$ of labeled cartesian squares 
as follows where $\sfr_{i,j}$ is the square that is the $i$-th one from 
the top and the $j$-th one from the left,
\[
\begin{CD}
X_{1,1} @>{H_{1,1}}>> X_{1,2} @>{H_{1,2}}>> \cdots 
 @>{H_{1,m-1}}>> X_{1,m} @>{H_{1,m}}>> X_{1,m+1} \\ 
@VV{V_{1,1}}V @VV{V_{1,2}}V @. @VV{V_{1,m}}V @VV{V_{1,m+1}}V \\
X_{2,1} @>>{H_{2,1}}> X_{2,2} @>>{H_{2,2}}> \cdots 
 @>>{H_{2,m-1}}> X_{2,m} @>>{H_{2,m}}> X_{2,m+1} \\ 
@VV{V_{2,1}}V @VV{V_{2,2}}V @. @VVV @VV{V_{2,m+1}}V \\
X_{3,1} @>>> X_{3,2} @>>> \cdots @>>> X_{3,m} @>>> X_{3,m+1} \\ 
@VVV @VVV @. @VVV @VVV \\
\cdots @. \cdots @.  @. \cdots @. \cdots\\
@VVV @VVV @. @VVV @VVV \\
X_{n,1} @>>> X_{n,2} @>>> \cdots @>>> X_{n,m} @>>> X_{n,m+1} \\ 
@VV{V_{n,1}}V @VV{V_{n,2}}V @. @VV{V_{n,m}}V @VV{V_{n,m+1}}V \\
X_{n+1,1} @>{H_{n+1,1}}>> X_{n+1,2} @>{H_{n+1,2}}>> \cdots 
 @>{{H_{n+1,m-1}}}>> X_{n+1,m} @>{H_{n+1,m}}>> X_{n+1,m+1}  
\end{CD}
\]
and such that the following two conditions are satisfied.
\begin{itemize}
\item The rightmost column of maps is $\sigma_1$, 
\newline i.e., $X_{i,m+1} = Y_i$ and $V_{i,m+1} = (f_i,\alpha_i)$.
\item The bottommost row of maps is $\sigma_2$, 
\newline i.e., $X_{n+1,j} = Z_j$ and $H_{n+1,j} = (g_j,\gamma_j)$.
\end{itemize}

It follows from the definitions that $X_{i,j} = Y_i \times_T Z_j$
and that for all $i,j$, $V_{i,j}$ has label $\alpha_i$ and $H_{i,j}$
has label $\gamma_j$. We denote the underlying maps of 
$V_{i,j}$ and $H_{i,j}$ by ${\rv}_{i,j}$ and ${\rh}_{i,j}$ 
respectively. 
We denote the $i$-th row of horizontal maps by $H_i$ (i.e., $H_i$ 
is the sequence consisting of $H_{i,j}$)
and the $j$-th column of vertical maps by~$V_j$. 

The existence of a fiber-product diagram associated to $\sigma_1,\sigma_2$
is obvious. As in the case of products of single maps, the 
fiber-product diagram $\sS$ is not unique as a function of 
$\sigma_1$ and $\sigma_2$ though it is unique 
up to a unique isomorphism (at each vertex). In particular,
the terms such as $V_{i,j},H_{i,j}$ are really functions of $\sS$.


Sometimes we use draw $\sS$ in a more compact form as follows:
\[
\xymatrix{
X_{1,1} \ar@{=>}[d]_{\sigma_3} \ar@{=>}[r]^{\sigma_4} 
 & Y_1 \ar@{=>}[d]^{\sigma_1} \\
Z_1 \ar@{=>}[r]^{\sigma_2} & T
}
\]
where $\sigma_3 \set V_1(\sS)$ and $\sigma_4 \set H_1(\sS)$.
This diagram suppresses a lot of information on $\sS$
but is useful when the extra information is not directly needed.
By abuse of language, we also call $\sS$ a cartesian square
(of double arrows) on $(\sigma_1,\sigma_2)$. The usual cartesian
square on single arrows is then called a unit square.

\ssss
\label{sssec:basech}
Our next goal is to associate to the above diagram $\sS$, a 
generalized base-change isomorphism 
\begin{equation}
\label{eq:basech1}
\boldsymbol{\beta}_{\sS} \colon 
H_1^{\bxt}V_{m+1}^{\bxt} \iso V_1^{\bxt}H_{n+1}^{\bxt}. 
\end{equation}
The idea behind defining $\boldsymbol{\beta}_{\sS}$
is quite simple. Let us think of $\sS$ as a directed
planar graph in the obvious sense.
Then the sequences $H_1 \star V_{m+1}$ 
and $V_1 \star H_{n+1}$ are
the two outermost directed paths from $X_{1,1}$ to $X_{n+1,m+1}$.
Now the base-change isomorphism $\bbeta_{\sfr_{i,j}}$ 
associated to~$\sfr_{i,j}$
may be thought of as a ``flip'' from $H_{i,j} \star V_{i,j+1}$
to $V_{i,j} \star H_{i+1,j}$ (as directed paths between~$X_{i,j}$ 
and~$X_{i+1,j+1}$). Defining $\boldsymbol{\beta}_{\sS}$ now
simply amounts to choosing a sequence of flips that transforms 
$H_1 \star V_{m+1}$ to $V_1 \star H_{n+1}$. Any such sequence 
of flips/squares necessarily starts from~$\sfr_{1,m}$ and ends 
with~$\sfr_{n,1}$ and uses each of the~$mn$ squares exactly once.
But the order in which the intermediate flips can be chosen is 
not unique. 

Thus, what needs to be addressed is why the resulting definition 
of $\boldsymbol{\beta}_{\sS}$ is independent of the choice 
of the sequence of flips. As it turns out this holds for 
simple functorial reasons, i.e., no special 
property of the base-change isomorphisms for unit squares, 
other than functoriality, comes into 
play. Indeed, the question essentially reduces to an 
elementary graph-theoretic exercise.
Since we use this kind of argument on a separate 
occasion too, we provide a sketch of the details involved.

\ssss
\label{sssec:graph}
This is an interlude on some graph-theoretic notions
needed for resolving the question of well-definedness of $\bbeta_{\sS}$. 
We use the same notation as above. 
For convenience, we work more generally with directed paths 
between two arbitrary vertices in~$\sS$. Containment of paths 
is defined in the obvious manner. Now we discuss some terminology.

$\bullet$ Given two (directed) paths $p,p'$, we say $p' \ge p$ if
both share the same starting vertex and the same ending vertex and if 
for any vertex $X_{i,j}$ in $p$ there exists a vertex $X_{i',j'}$
in $p'$ such that $i' \le i$ and $j' \ge j$. Clearly $\ge$ 
is reflexive and transitive. Moreover, if $p' \ge p$
then $p'$ and $p$ do not cross each other anywhere. In fact,
for any vertex common to $p$ and $p'$, say $X_{i,j}$, the vertical
edge $V_{i-1,j}$ ending at $X_{i,j}$ is in $p$ only if it is also
in $p'$ and the horizontal edge $H_{i,j}$ starting from $X_{i,j}$
is in $p$ only if it is also in $p'$. One easily concludes
that $\ge$ is in fact antisymmetric, i.e., if $p \ge p'$ and $p' \ge p$
then $p=p'$. 

Thus $\ge$ is a partial ordering on paths. 
We say $p'>p$ if $p' \ge p$ and $p' \ne p$. 

$\bullet$ We say that a square $\sfr_{i,j}$ flips a path~$p_1$ 
to another one~$p_2$ or that $\sfr_{i,j}$ flops~$p_2$ to~$p_1$ if
the north-eastern half of the boundary of~$\sfr_{ij}$, 
viz.~$H_{i,j} \star V_{i,j+1}$, 
is contained in~$p_1$, the other half
$V_{i,j} \star H_{i+1,j}$ is in~$p_2$ and if~$p_2$ is obtained 
from~$p_1$ by interchanging these halves.

$\bullet$ A square $\sfr = \sfr_{i,j}$ is 
said to be \emph{open} for a path~$p$ if
$H_{i,j} \star V_{i,j+1}$ is contained in~$p$. This is equivalent to 
saying that there exists another path $p'$ such that $\sfr$ 
flips~$p$ to~$p'$. Note that, in this situation, if~$\rfr$ is 
another square distinct from $\sfr$ that is also open for~$p$,
then~$\rfr$ is also open for~$p'$.

$\bullet$ A path $p_1$ is said to be \emph{adjacent} to another path $p_2$ 
if $p_1 > p_2$ and for any path $p$ such that $p_1 \ge p \ge p_2$
we have $p= p_1$ or $p=p_2$. It is easily verified that $p_1$
is adjacent to~$p_2$ if and only if 
there exists a square $\sfr$, necessarily unique,
such that~$\sfr$ flips~$p_1$ to~$p_2$.

$\bullet$ A \emph{maximal} chain of paths is a 
sequence $p_1 > p_2 > \cdots > p_N$, 
($N \ge 1$) such that~$p_i$ is adjacent to~$p_{i+1}$ for all~$i$.

$\bullet$ Two maximal chains 
\[
c \set \quad p_1 > p_2 > \cdots > p_N \qquad \text{and} \qquad
c' \set \quad p_1' > p_2' > \cdots > p_N'
\] 
are said to be \emph{flip-flop neighbors} if one of the following
conditions holds:
\newline (a) $c = c'$; or 
\newline (b) $N>2$ and there exists an integer $i$ such that 
$1<i<N$ and such that $p_j = p_j'$ for $j \ne i$.
\newline In the situation of (b), with $c \ne c'$, if $\rfr$ (resp.~$\sfr$)
is the square that flips $p_{i-1}$ to $p_i$ (resp.~$p_i$ to $p_{i+1}$),
then $\sfr$ (resp.~$\rfr$) flips $p_{i-1}'$ to $p_i'$ 
(resp.~$p_i'$ to $p_{i+1}'$). Moreover, $\rfr$ and $\sfr$ do not 
share a common edge.

\medskip
Here then is the main graph-theoretic result that we use.

\begin{alem}
\label{lem:graph1}
For any two maximal chains 
\[
c\set \quad p_1 > p_2 > \cdots > p_N \qquad \text{and} \qquad
c' \set \quad p_1' > p_2' > \cdots > p_N'
\] 
satisfying $p_1 = p_1'$, $p_n = p_N'$ there exists a 
sequence of maximal chains $c_1, c_2, \ldots, c_M$ such that 
$c_1 = c$, $c_M = c'$ and such that for all $i$,
$c_i$ and $c_{i+1}$ are flip-flop neighbors. 
\end{alem}
\begin{proof}
We use induction on $N$, the length of the chains $c,c'$. 
If $N \le 3$ then the lemma holds trivially. Assume $N>3$. We 
claim that there exists another maximal chain 
$c''\set \;\; p_1'' > p_2'' > \cdots > p_N''$ satisfying
$p_1'' = p_1, p_N'' = p_N$ and there are integers
$i,j$ such that $1<i,j<N$ and such that $p_i'' = p_i$ and $p_j'' = p_j'$.

Assuming the above claim we proceed as follows. It suffices to find
a sequence of flip-flop neighbors between $c,c''$ and between $c'',c'$.
In the case of $c$ and $c''$ we may break up each of these chains into 
two pieces as follows
\begin{align}
c_{\le i} &= \quad p_1 > \cdots > p_i,& 
c_{\ge i} &= \quad p_i > \cdots > p_N, \notag \\
c_{\le i}'' &= \quad p_1'' > \cdots > p_i'',&
c_{\ge i} &= \quad p_i'' > \cdots > p_N''. \notag
\end{align}
By induction hypothesis there exists a sequence of flip-flop neighbors 
between~$c_{\le i}$ and~$c''_{\le i}$ and also one between~$c_{\ge i}$ 
and~$c''_{\ge i}$. These can be put together to give one 
for~$c,c''$. Similarly one argues for $c'',c'$.

It remains to prove the above claim. If $p_2' = p_2$ then the claim 
is proven by choosing $c''= c'$, and $i = j = 2$. Assume $p_2' \ne p_2$.
Let $\sfr_{(k)}$ (resp.~$\sfr'$) be the square that flips $p_k$ to $p_{k+1}$ 
(resp.~$p_1'$ to $p_2'$). Then $\sfr'$ is open for $p_1$ but not for 
$p_N = p_N'$. Let $i$ be the largest integer such that $\sfr'$ is 
open for $p_i$. Then $\sfr' = \sfr_{(i)}$. 
We now define~$c''$ inductively as follows.
By definition, $p_k'' = p_k$ for $k= 1,N$. Let~$p_2''$ be the 
path obtained by flipping $p_1$ via $\sfr'$. For $1 < k < i+1$,
let $p_{k+1}''$ be the path obtained by flipping $p_{k}''$ via $\sfr_{(k-1)}$
and for $k \ge i+1$, let $p_{k+1}''$ be obtained by flipping~$p_{k}''$ 
via~$\sfr_{(k)}$. Then for this choice of $c''$ and $i$ and with $j = 2$,
the claim is verified.
\end{proof}

\ssss
\label{sssec:invariance}
We return to the issue of the well-definedness of $\bbeta_{\sS}$.
For any path~$p$ from~$X_{1,1}$ to~$X_{n+1,m+1}$, let 
$p^{\bxt}$ be the obvious functor $\cD_{X_{n+1,m+1}} \to \cD_{X_{1,1}}$ 
induced via $p$. Then for any maximal chain 
$c$ given by say, $p_1 > p_2 > \cdots > p_{n+m}$, where 
$p_1 = H_1 \star V_{m+1}$ and $p_{n+m} = V_1 \star H_{n+1}$, if $\sfr_{(k)}$
is the square that flips~$p_k$ to~$p_{k+1}$ then one
obtains an isomorphism~$\bbeta_{\sS}(c)$ given by 
\[
p_1^{\bxt} \xto{\text{via }\bbeta_{\sfr_{(1)}}} p_2^{\bxt} 
\xto{\text{via }\bbeta_{\sfr_{(2)}}}
\quad \cdots \quad \xto{\text{via }\bbeta_{\sfr_{(n+m-1)}}} p_{n+m}^{\bxt}.   
\]
Thus the problem may be
rephrased as saying that for any two maximal chains
$c,c'$ that start from $H_1 \star V_{m+1}$ and end 
with $V_1 \star H_{n+1}$, we have 
$\bbeta_{\sS}(c) = \bbeta_{\sS}(c')$.

By Lemma \ref{lem:graph1}, it suffices to prove equality
in the case where $c$ and $c'$ are flip-flop neighbors and $c \ne c'$. 
Let $p_k$ (resp.~$p_k'$) be the $k$-th element of the 
chain $c$ (resp.~$c'$). Let $i$ be the unique integer such that
$p_i \ne p_i'$. If $\rfr$ is the square that flips $p_{i-1}$ to $p_i$
and~$\sfr$ the one that flips~$p_{i}$ to~$p_{i+1}$, then the
following diagram commutes for functorial reasons. 
\[
\begin{CD} 
p_{i-1}^{\bxt} @>{\text{via }\bbeta_{\rfr}}>> p_{i}^{\bxt} \\
@V{\text{via }\bbeta_{\sfr}}VV @VV{\text{via }\bbeta_{\sfr}}V \\
p_{i}'{}^{\bxt} @>>{\text{via }\bbeta_{\rfr}}> p_{i+1}^{\bxt}
\end{CD}
\]
Now $\bbeta_{\sS}(c) = \bbeta_{\sS}(c')$ follows immediately.

The following lemma generalizes \ref{lem:cafpols4}(i)
and is used in \S\ref{sec:proofsIII}.


\begin{alem}
\label{lem:cafpols20}
Let $\sigma \colon Z \Lra X$ be a labeled sequence. 
Corresponding to any cartesian diagram~$\sS$ as follows
where the vertical maps in the intermediate places are all
$\sfP$-labeled identity maps, 
\[
\xymatrix{
Z \ar@{=>}[rr]^{\sigma} \ar@{=}[d]_{I_1 \> = \> (1_Z, \sfP)}  
 & & X \ar@{=}[d]^{I_2 \> = \> (1_X, \sfP)} \\
Z \ar@{=>}[rr]^{\sigma} & & X
}
\] 
the isomorphism $\sigma^{\bxt} = \sigma^{\bxt}I_2^{\bxt}
\xto{\bbeta_{\sS}} I_1^{\bxt}\sigma^{\bxt}  = \sigma^{\bxt}$
is the identity.
\end{alem}
\begin{proof}
We prove by induction on the length of $\sigma$. The case of length 1
has been dealt with in~\ref{lem:cafpols4}(i). Now suppose 
$\sigma = \sigma_1 \star \sigma_2$. Consider the following associated 
cartesian diagram with obvious notation. Let $\sS_1$ denote the 
square on the left and~$\sS_2$ the one on the right.
\[
\xymatrix{
Z \ar@{=>}[r]^{\sigma_1} \ar@{=}[d]_{I_1} & Y 
 \ar@{=>}[r]^{\sigma_2} \ar@{=}[d]_{J} & X \ar@{=}[d]^{I_2} \\
Z \ar@{=>}[r]^{\sigma_1} & Y \ar@{=>}[r]^{\sigma_2} & X
}
\] 
The assertion of the lemma amounts to checking that the outer border
of the following diagram commutes.
\[
\begin{CD}
\sigma_1^{\bxt}\sigma_2^{\bxt} @= \sigma_1^{\bxt}\sigma_2^{\bxt} 
 @= \sigma_1^{\bxt}\sigma_2^{\bxt} \\ 
@| @| @| \\
\sigma_1^{\bxt}\sigma_2^{\bxt}I_2^{\bxt}
 {\smash{\raisebox{5ex}{\makebox[0pt]{\qquad\qquad\dg{2}}}}} 
 @>>{\sigma_1^{\bxt}(\bbeta_{\sS_2})}> \sigma_1^{\bxt}J^{\bxt}\sigma_2^{\bxt} 
 {\smash{\raisebox{5ex}{\makebox[0pt]{\qquad\qquad\dg{1}}}}}
 @>>{\bbeta_{\sS_1}(\sigma_2^{\bxt})}> 
 I_1^{\bxt}\sigma_1^{\bxt}\sigma_2^{\bxt} 
\end{CD}
\] 
Since $\sigma_i$ has length smaller than $\sigma$, 
by induction hypothesis the conclusion holds
for $\sigma_i$ and hence \dg{$i$} commutes. 
\end{proof}

\subsection{The generalized fundamental isomorphism}
\label{subsec:genfund}
We upgrade the fundamental isomorphism of~\eqref{eq:cafpols3}
to the level of sequences of arbitrary length. Modulo the issue
of canonicity, the generalization is achieved in \eqref{eq:psi3}.
The isomorphism defined here is one of the
important special cases as well as a building block of our central 
construction of $\Psi_{-,-}$ in \S\ref{subsec:psi}.

\ssss
\label{sssec:psi-1}
Consider a labeled sequence $\sigma$ given by 
\begin{equation}
\label{eq:psi1}
X = Y_1 \xto{\;(f_1,\lambda_1)\;} Y_2 \xto{\;(f_2,\lambda_2)\;} Y_3 
\xto{\quad} \cdots \xto{\quad}
Y_{n} \xto{\;(f_{n},\lambda_{n})\;} Y_{n+1} = X
\end{equation}
such that $|\sigma| \set 1_X$. Our aim is to define an isomorphism
$\sigma^{\bxt} \iso  \oneD{X}$ that generalizes the fundamental 
isomorphism of \eqref{eq:cafpols3}. This is accomplished by means
of a staircase diagram that we now define.

A \emph{staircase diagram based on} $\sigma$ 
is a collection of $n(n-1)/2$ cartesian squares stacked into 
$n$ columns as follows, (where the $i$-th column from the left 
has $n-i$ squares in it)
\begin{equation}
\label{eq:psi1a}
\begin{CD}
X \\
@VV{V_{1,1}}V \\
Y_{2,1} @>{H_{2,1}}>> X \\
@VV{V_{2,1}}V @VV{V_{2,2}}V \\
\cdots @. \cdots @. \cdots  \\
@VVV @VVV @. \\
Y_{n-1,1} @>>> Y_{n-1,2} @>>> \cdots @>>> X \\
@VV{V_{n-1,1}}V @VV{V_{n-1,2}}V @. @VV{V_{n-1,n-1}}V \\
Y_{n,1} @>{H_{n,1}}>> Y_{n,2} @>{H_{n,2}}>> \cdots 
 @>{H_{n,n-2}}>> Y_{n,n-1} @>{H_{n,n-1}}>> X \\
@VV{V_{n,1}}V @VV{V_{n,2}}V @. @VV{V_{n,n-1}}V @VV{V_{n,n}}V  \\
X @>>> Y_2 @>>> \cdots  @>>> Y_{n-1} @>>> Y_{n} @>>> X
\end{CD}
\end{equation}
satisfying the following conditions.
\begin{itemize}
\item The sequence of maps at the base is $\sigma$. Thus we may set 
$H_{n+1,i} \set (f_i, \lambda_i)$.
\item The vertical maps, $V_{i,j}$ all have label $\sfP$. 
\item Let $\rv_{i,j} = |V_{i,j}|$ and $\rh_{i,j} = |H_{i,j}|$. 
Then for every $i$, the following hold.
\[
\rh_{i+1,i}\;\rv_{i,i} = 1_X \qquad \qquad 
\rv_{n,i}\;\rv_{n-1,i} \;\cdots \;\rv_{i,i} = f_{i-1}f_{i-2} \cdots f_1
\]
\end{itemize}
A quick approach to constructing a staircase is obtained by noting that 
$Y_{i,j}$ must equal the product $Y_j \times_{Y_i} X$ where 
$Y_j \to Y_i$ (resp.~$X \to Y_i$) is the composite map $f_i \cdots f_j$ 
(resp.~$f_i \cdots f_1$);
then ${\rh}_{i,j}$ and ${\rv}_{i,j}$ are the obvious natural maps
induced by~$f_j$ and~$f_i$ respectively. It also follows that 
for a fixed $\sigma$, any two staircases based on it are isomorphic
via unique isomorphisms at each vertex. 

Nevertheless, the lack of uniqueness of a staircase, 
means that all the terms used in definition of
a staircase $\cS$ such as $V_{-,-}$, $H_{-,-}$, $Y_{-,-}$, etc.,
should be regarded as functions of~$\cS$, e.g., $V_{-,-}(\cS)$.
In specific contexts, when the choice of~$\cS$ is known, reference
to~$\cS$ in the notation may be dropped. 

For $\sigma$ of \eqref{eq:psi1}, let us assume that a staircase $\cS$
such as the one drawn above has been chosen. 
We now show how $\cS$ gives rise to an isomorphism 
${\sigma}^{\bxt} \iso \oneD{X}$.

Each ``step'' in $\cS$ gives rise to a functor isomorphic to
$\oneD{X}$, i.e., for $1 \le i \le n$, we have
${\rh}_{i+1,i} \; {\rv}_{i,i} = 1_{X}$ and hence 
\eqref{eq:cafpols3} gives an isomorphism 
\[
V_{i,i}^{\bxt}H_{i+1,i}^{\bxt} \iso \oneD{X}.
\] 
Set $\texttt{Steps} = (V_{1,1} \star H_{2,1}) 
\star (V_{2,2} \star H_{3,2}) \star \cdots \star (V_{n,n} \star H_{n+1,n})$
and let 
\begin{equation}
\label{eq:psi2}
\texttt{Steps}^{\bxt} \iso \oneD{X}
\end{equation}
be the isomorphism obtained by successively using \eqref{eq:cafpols3}.

Let $V_1$ be the sequence occurring as the leftmost column of maps 
in~$\cS$. By definition, $V_1$ composes to~$1_{X}$
and the label of each map in it is~$\sfP$. Therefore 
pseudofunctoriality of~$(-)^{\times}$ gives a natural isomorphism
\begin{equation}
\label{eq:psi2a}
V_1^{\bxt} \iso (1_X)^{\times} = \oneD{X}.
\end{equation}

The \emph{generalized} fundamental isomorphism associated to a staircase
$\cS$ based on~$\sigma$ is defined to be the isomorphism
\begin{equation}
\label{eq:psi3}
\Phi_{\sigma}(\cS) \colon {\sigma}^{\bxt} \iso \oneD{X}
\end{equation}
obtained as the following composition
\[
{\sigma}^{\bxt} \xto{\text{via }\eqref{eq:psi2a}^{-1}} 
V_1^{\bxt}{\sigma}^{\bxt} 
\xto{\quad} \texttt{Steps}^{\bxt} \xto{\eqref{eq:psi2}} \oneD{X}
\]
where the map in the middle is defined the same way as $\bbeta_{\sS}$
is defined in \eqref{eq:basech1}, viz., through the base-change isomorphisms
associated to each of the unit squares in $\cS$. 
(As in the case of $\bbeta_{\sS}$, the order in which the squares are chosen 
can vary but the end result is the same in view of the arguments in 
\ref{sssec:invariance}.)

A-priori, for a fixed $\sigma$,
the generalized fundamental isomorphism depends on the 
choice of the staircase $\cS$ based on~$\sigma$, in that, 
if $\cS'$ is another staircase based on~$\sigma$, then
it is not clear whether 
$\Phi_{\sigma}(\cS) = \Phi_{\sigma}(\cS')$. We would like to say that
the two are indeed equal so that \eqref {eq:psi3} 
defines a canonical isomorphism, independent of the choice 
of $\cS$. Fortunately, this is true, though it is not altogether
trivial to verify. We take up this issue in the next subsection. 

Let us verify that $\Phi_{\sigma}(\cS)$ does indeed generalize 
the earlier definitions of a fundamental isomorphism.
\begin{alem}
\label{lem:psi8a}
Let $\sigma$ be the sequence $X \xto{\;F_1\;} Y \xto{\;F_2\;} X$
of~\textup{\ref{sssec:cafpols3z}}, i.e., $F_1$ has label $\>\sfP$ and 
$|F_1 \star F_2| = 1_X$. Let $\cS$ be a staircase 
based on $\sigma$. Then $\Phi_{\sigma}(\cS) = \bphi_{F_1,F_2}$.
\end{alem}
\begin{proof}
Suppose $\cS$ is given as follows. Let $\sfr$ be
the cartesian square in it.
\[
\begin{CD}
X \\
@V{\Delta}VV \\
X^2  @>{Q}>> X \\
@V{P}VV @VV{F_1'= F_1}V \\
X @>{F_1}>> Y @>{F_2}>> X
\end{CD}
\] 
By construction, $\Delta$, $P$, $Q$ and $F_1$, all  have label $\sfP$. 
Hence the following composite isomorphism $\theta$,
\[
F_1^{\bxt} \xto{\text{via $\boldC^{-1}_{\Delta,P}$}} 
\Delta^{\!\bxt}P^{\bxt}F_1^{\bxt} 
\xto{\text{via $\bbeta_{\sfr}$}} 
\Delta^{\!\bxt}Q^{\bxt}F_1^{'\bxt} 
\xto{\text{via $\bphi_{\Delta,Q}$}} F_1^{'\bxt},
\]
is the identity for pseudofunctorial reasons. 
By definition, $\Phi_{\sigma}(\cS)$
is the following composition and thus the lemma follows.
\[
F_1^{\bxt}F_2^{\bxt} \xto{\text{via $\theta$ above}} F_1^{'\bxt}F_2^{\bxt}
\xto{\bphi_{F_1',F_2}} \oneD{X}
\]
\end{proof}

\subsection{Canonicity of the fundamental isomorphism}
\label{subsec:canon}
Our main concern here is of showing that the 
generalized fundamental isomorphism of~\eqref{eq:psi3} is 
independent of the choice of the staircase diagram chosen.
We accomplish this in Proposition~\ref{prop:canon8} below.

The general idea is of comparing the base-change isomorphism
$\bbeta_{-}$ resulting from isomorphic cartesian squares.
We begin by setting up a comparison isomorphism $\boldC_{-,-}$
for the composition of two maps where one of them is an isomorphism
with label $\sfP$.

In what follows, in the context of labeled maps, a vertically 
drawn ``equal sign'' is downward pointing by default and 
a horizontally drawn ``equal sign'' points rightward.

\begin{adefi}
\label{def:canon1}
Let $i \colon Y \to Y'$ and $j \colon Y' \to Y$ be inverse isomorphisms. 
Set $I \set (i,\sfP), J \set (j, \sfP)$. 
\begin{enumerate}
\item For any labeled map $X \xto{F\>=\>(f,\lambda)} Y$, 
with $F'$ as the map $X \xto{if} Y'$ with  
label $\lambda$, 
(since $i$ is in both, $\sfP$ and $\sfO$, hence~$F'$ 
can have the same label as~$F$)
we define $\boldC_{F,I} \colon F^{\bxt}I^{\bxt} \iso F'{}^{\bxt}$
to be the base-change isomorphism~$\bbeta_{\sfr}$ associated to
the cartesian square~$\sfr$ drawn as follows.
\[
\xymatrix{
X \ar@{=}[d]_{(1_X, \sfP)} \ar[r]^F & Y \ar[d]^{I} \\
X \ar[r]^{F'} & Y'
}
\] 
\item For any labeled map $Y \xto{G\>=\>(g,\lambda)} Z$, with 
$G' \set (gj,\lambda) \colon Y' \to Z$, we define 
$\boldC_{J,G} \colon J^{\bxt}G^{\bxt} \iso G'{}^{\bxt}$
to be the base-change isomorphism $\bbeta_{\sfr}$ associated to
the cartesian square~$\sfr$ drawn as follows.
\[
\xymatrix{
Y' \ar[d]_{G'} \ar[r]^J & Y \ar[d]^{G} \\
Z \ar@{=}[r]^{(1_Z, \sfP)} & Z
}
\]
\end{enumerate}
\end{adefi}

\begin{arems}
\label{rem:canon2}
\hfill

1. In \ref{def:canon1}(i), by interchanging the role of $i$ with 
$j$ and $F$ with $F'$ we obtain in an analogous fashion an isomorphism
$\boldC_{F',J} \colon F'{}^{\bxt}J^{\bxt} \iso F^{\bxt}$. 
Likewise, in (ii) we obtain an isomorphism
$\boldC_{I,G'} \colon I^{\bxt}G'{}^{\bxt} \iso G^{\bxt}$.

2. In \ref{def:canon1}(i), if $F$ has label $\sfP$, then 
every map in $\sfr$ has label $\sfP$ and hence by pseudofunctoriality 
of $(-)^{\times}$ we see that 
the above definition of $\boldC_{F,I}$ agrees with the 
one in~\ref{sssec:cafpols0}. 
A similar remark holds for $\boldC_{G,J}$ of \ref{def:canon1}(ii).
\end{arems}

\begin{alem}
\label{lem:canon3}
Let notation be as in \textup{\ref{def:canon1}}. Let the following  
isomorphisms
\[
\varphi \colon I^{\bxt}J^{\bxt} \iso \oneD{Y}
\qquad \text{and} \qquad 
\psi \colon J^{\bxt}I^{\bxt} \iso \oneD{Y'}
\] 
be the ones obtained from the pseudofunctoriality of $(-)^{\times}$.
Then the following four diagrams of isomorphisms commute.
\[
\xymatrix{
& F^{\bxt}I^{\bxt}J^{\bxt} \ar[d]_{\textup{via}}^{\!\boldC_{F,I}}   
 \ar[ld]_{\textup{via } \varphi}  \\
F^{\bxt} & F'{}^{\bxt}J^{\bxt} \ar[l]^{\boldC_{F',J}}  
} \!
\xymatrix{
& F'{}^{\bxt}J^{\bxt}I^{\bxt} \ar[d]_{\textup{via}}^{\!\boldC_{F'\!,J}}   
 \ar[ld]_{\textup{via } \psi}  \\
F'{}^{\bxt} & F^{\bxt}I^{\bxt} \ar[l]^{\boldC_{F,I}}  
} \!\!
\xymatrix{
& I^{\bxt}J^{\bxt}G^{\bxt} \ar[d]_{\textup{via}}^{\!\boldC_{J,G}}   
 \ar[ld]_{\textup{via } \varphi}  \\
G^{\bxt} & I^{\bxt}G'{}^{\bxt} \ar[l]^{\boldC_{I,G'}}  
} \!
\xymatrix{
& J^{\bxt}I^{\bxt}G'{}^{\bxt} \ar[d]_{\textup{via}}^{\!\boldC_{I,G'}}   
 \ar[ld]_{\textup{via } \psi}  \\
G'{}^{\bxt} & J^{\bxt}G^{\bxt} \ar[l]^{\boldC_{J,G}}  
}
\]
\end{alem}
\begin{proof}
All the four diagrams are handled in a similar way. For the first two
diagrams from the left, we use vertical transitivity of base-change 
(Lemma \ref{lem:cafpols4}(iii))
corresponding to the following two diagrams respectively.
\[
\xymatrix{
X \ar@{=}[d]_{(1_X, \sfP)} \ar[r]^F & Y \ar[d]^{I} \\
X \ar@{=}[d]_{(1_X, \sfP)} \ar[r]^{F'} & Y' \ar[d]^{J} \\
X \ar[r]^F & Y  \\
}
\hspace{8em}
\xymatrix{
X \ar@{=}[d]_{(1_X, \sfP)} \ar[r]^{F'} & Y' \ar[d]^{J} \\
X \ar@{=}[d]_{(1_X, \sfP)} \ar[r]^{F} & Y \ar[d]^{I} \\
X \ar[r]^{F'} & Y'  \\
}
\]
A similar argument works for the other two cases.
\end{proof}

For the next lemma, the basic premise is that there is a cartesian square
$\sfr$ drawn as follows where $F$ and $F_1$ 
have label $\sfP$.
\[
\begin{CD}
W @>{G_1}>> X \\
@V{F_1}VV @VVFV \\
Z @>G>> Y
\end{CD}
\]
We shall now consider cases where exactly one of its vertices 
is modified via a $\sfP$-labeled
isomorphism. In each case the new edges shall be 
called $F',G',F_1'$ or~$G_1'$ as applicable and each of these new maps shall 
have the same label as that of its old counterpart.
Also, in each case, the modified square shall be called~$\sfr'$. 

\begin{alem}
\label{lem:canon4}
With notation and convention as above, the following hold.
\begin{enumerate}
\item Let $I \colon W \to W'$ and $J \colon W' \to W$ be  
$\sfP$-labeled inverse isomorphisms. 
Then the following diagram commutes.
\[
\begin{CD}
G_1^{\bxt}F^{\bxt} @>{\bbeta_{\sfr}}>> F_1^{\bxt}G^{\bxt} \\
@V{\textup{via } \boldC_{I,G'_1}^{-1}}VV 
 @VV{\textup{via } \boldC_{I,F_1'}^{-1}}V \\
I^{\bxt}G'_1{}^{\bxt}F^{\bxt} @>{\textup{via}}>{\bbeta_{\sfr'}}> 
 I^{\bxt}F_1'{}^{\bxt}G^{\bxt}
\end{CD}
\]
\item Let $I \colon Z \to Z'$ and $J \colon Z' \to Z$ be  
$\sfP$-labeled inverse isomorphisms. 
Then the following diagram commutes.
\[
\xymatrix{
& & F_1^{\bxt}G^{\bxt} \ar[lld]_{\bbeta_{\sfr}^{-1}} 
 \ar[rd]^{\textup{via }\boldC_{I,G'}^{-1}} \\
G_1^{\bxt}F^{\bxt} \ar[rrd]_{\bbeta_{\sfr'}} & F_1'{}^{\bxt}J^{\bxt}G^{\bxt} 
 \ar[ru]_{\textup{via }\boldC_{F_1',J}} \ar[rd]^{\textup{via }\boldC_{J,G}} 
 & & F_1^{\bxt}I^{\bxt}G'{}^{\bxt} \\
& & F_1'{}^{\bxt}G'{}^{\bxt}  \ar[ru]_{\textup{via }\boldC_{F_1,I}^{-1}}
}
\]
\item Let $I \colon X \to X'$ and $J \colon X' \to X$ be  
$\sfP$-labeled inverse isomorphisms. 
Then the following diagram commutes.
\[
\xymatrix{
 & G_1^{\bxt}F^{\bxt} \ar[rd]_{\textup{via }\boldC_{I,F'}^{-1}} 
 \ar[rrd]^{\bbeta_{\sfr}} \\
G_1'{}^{\bxt}J^{\bxt}F^{\bxt} \ar[ru]^{\textup{via }\boldC_{G_1',J}} 
 \ar[rd]_{\textup{via }\boldC_{J,F}} 
 & & G_1^{\bxt}I^{\bxt}F'{}^{\bxt} & F_1^{\bxt}G^{\bxt} 
 \ar[lld]^{\bbeta_{\sfr'}^{-1}} \\
 & G_1'{}^{\bxt}F'{}^{\bxt} \ar[ru]^{\textup{via }\boldC_{G_1,I}^{-1}}
}
\]
\item Let $I \colon Y \to Y'$ and $J \colon Y' \to Y$ be  
$\sfP$-labeled inverse isomorphisms. 
Then the following diagram commutes.
\[
\begin{CD}
G_1^{\bxt}F^{\bxt} @>{\bbeta_{\sfr}}>> F_1^{\bxt}G^{\bxt} \\
@V{\textup{via }\boldC_{F',J}^{-1}}VV @VV{\textup{via }\boldC_{G',J}^{-1}}V \\
G_1^{\bxt}F'{}^{\bxt}J^{\bxt} @>{\textup{via}}>{\bbeta_{\sfr'}}> 
 F_1^{\bxt}G'{}^{\bxt}J^{\bxt}
\end{CD}
\]
\end{enumerate}
\end{alem}
\begin{proof}
We only give proofs of (i) and (ii). The remaining two cases are handled 
similarly.

In (i), commutativity follows by transitivity of base-change associated
to the following diagram where $\sfr'$ is the square at the bottom
and $\sfr$ is the composite square. (For $\boldC_{I,F_1'}$ we also refer to 
part 2 of \ref{rem:canon2}.)
\[
\xymatrix{
W \ar[d]_{I} \ar[r]^{G_1} & X \ar@{=}[d]^{(1_X,\sfP)} \\
W' \ar[d]_{F_1'} \ar[r]^{G_1'} & X \ar[d]^F \\
Z \ar[r]^G & Y
}
\]

In (ii), the diagram has two parts. The non-convex part on the left 
commutes because of transitivity of base-change associated to the
following diagram where $\sfr'$ is the square at the top
and $\sfr$ is the composite square. (For~$\boldC_{F_1',J}$ we also refer to 
part 2 of \ref{rem:canon2}.)
\[
\xymatrix{
W \ar[d]_{F_1'} \ar[r]^{G_1} & X \ar[d]^F \\
Z' \ar[d]_{J} \ar[r]^{G'} & Y \ar@{=}[d]^{(1_Y,\sfP)} \\
Z \ar[r]^G & Y
}
\]
The rhombus on the right in (ii) is the outer border of the following diagram
\[
\xymatrix{
& & & F_1^{\bxt}G^{\bxt} 
 \ar[rrrd]^{\textup{via }\boldC_{I,G'}^{-1}} \\
F_1'{}^{\bxt}J^{\bxt}G^{\bxt} \ar[rrru]^{\textup{via }\boldC_{F_1',J}} 
 \ar[rrrd]_{\textup{via }\boldC_{J,G}} 
 \ar[rrr]^{\textup{via }\boldC_{F_1,I}^{-1}} 
 & & & F_1^{\bxt}I^{\bxt}J^{\bxt}G^{\bxt} 
 \ar[u]^{\textup{via}}_{I^{\bxt}J^{\bxt} \cong \mathbf{1}}
 \ar[rrr]^{\textup{via }\boldC_{J,G}} & & & F_1^{\bxt}I^{\bxt}G'{}^{\bxt} \\
& & & F_1'{}^{\bxt}G'{}^{\bxt}  \ar[rrru]_{\textup{via }\boldC_{F_1,I}^{-1}}
}
\]
where the upper two triangles commute by \ref{lem:canon3} (see its first and 
third diagrams from the left) while 
the lower part commutes for functorial reasons. 
\end{proof}

We now consider some situations that arise in the proof 
of Proposition \ref{prop:canon8}. 
In each case, one looks at a local piece of a staircase
diagram. The goal is to describe the cumulative effect resulting 
from replacing a vertex by an isomorphic one. Here is the
setup of the first situation.


Consider the following two cartesian diagrams, differing only via 
the middle vertex and satisfying conditions described below.
\[
\begin{CD}
X_1 @>>> X_2 @>>> X_3 \\
@VVV @VVV @VVV \\
X_4 @>>> X_5 @>>> X_6 \\
@VVV @VVV @VVV \\
X_7 @>>> X_8 @>>> X_9
\end{CD}
\qquad \qquad \qquad
\begin{CD}
X_1 @>>> X_2 @>>> X_3 \\
@VVV @VVV @VVV \\
X_4 @>>> X_5' @>>> X_6 \\
@VVV @VVV @VVV \\
X_7 @>>> X_8 @>>> X_9
\end{CD}
\]
The maps in the above diagrams are named as follows.
For the diagram on the left, we use $F_{i,j}$ for
the drawn map $X_i \to X_j$. For the diagram on the right,
if $i \ne 5$ and $j \ne 5$ then the map $X_i \to X_j$ is the
same as~$F_{i,j}$ and for the remaining cases we use~$F_{i,j}'$. 
Finally, we name the diagram on the left as~$\sS$ and the one on the 
right as~$\sS'$.

The above two diagrams $\sS$ and $\sS'$ are assumed to further satisfy 
the following properties.
\begin{itemize}
\item The vertical maps are all $\sfP$-labeled.
\item Whenever $F_{i,j}'$ exists, it has the same label that $F_{i,j}$ has.
\item There are $\sfP$-labeled inverse isomorphisms
$I \colon X_5 \to X_5'$ and $J \colon X_5' \to X_5$ such that 
whenever $F_{i,j}'$ exists, $F_{i,j}$ and $F_{i,j}'$ are related 
via $I, J$ in the obvious manner. 
\end{itemize}

Using the base-change isomorphisms for each of the unit squares in
$\sS$ and likewise for $\sS'$, one obtains 
generalized base-change isomorphisms 
$\bbeta_{\sS}$ and $\bbeta_{\sS'}$ respectively, each of the 
following type
\[
F_{1,2}^{\bxt}F_{2,3}^{\bxt}F_{3,6}^{\bxt}F_{6,9}^{\bxt} \iso
F_{1,4}^{\bxt}F_{4,7}^{\bxt}F_{7,8}^{\bxt}F_{8,9}^{\bxt}.
\]

\begin{alem}
\label{lem:canon5}
In the above situation, $\bbeta_{\sS} = \bbeta_{\sS'}$.
\end{alem}

\begin{proof}
(cf. Remark \ref{rem:canon5a} below.) 
It suffices to prove that the following diagram commutes. Here the column
of four maps on the left end defines~$\bbeta_{\sS}$ through the  
base-change isomorphisms associated to the four cartesian squares in~$\sS$.
Likewise, the column at the right end defines~$\bbeta_{\sS'}$.
\[
\xymatrix{
F_{1,2}^{\bxt}F_{2,3}^{\bxt}F_{3,6}^{\bxt}F_{6,9}^{\bxt} 
 \ar[d]^{\hspace{10em} \blacksquare_1} \ar@{=}[rr]
 & & F_{1,2}^{\bxt}F_{2,3}^{\bxt}F_{3,6}^{\bxt}F_{6,9}^{\bxt} \ar[d] \\
F_{1,2}^{\bxt}F_{2,5}^{\bxt}F_{5,6}^{\bxt}F_{6,9}^{\bxt} 
 \ar[d]^{\hspace{4em} \blacksquare_2} \ar[r]^{6\quad}
 & F_{1,2}^{\bxt}F_{2,5}^{\bxt}I^{\bxt}F'_{5,6}{}^{\!\!\!\!\bxt}F_{6,9}^{\bxt} 
 \ar[d] \ar[r]^{2} & F_{1,2}^{\bxt}F'_{2,5}{}^{\!\!\!\!\bxt}
 F'_{5,6}{}^{\!\!\!\!\bxt}F_{6,9}^{\bxt} \ar[d] \\
F_{1,4}^{\bxt}F_{4,5}^{\bxt}F_{5,6}^{\bxt}F_{6,9}^{\bxt} 
 \ar[d]^{\hspace{4em} \blacksquare_3} \ar[r]^{6\quad}
 & F_{1,4}^{\bxt}F_{4,5}^{\bxt}I^{\bxt}F'_{5,6}{}^{\!\!\!\!\bxt}F_{6,9}^{\bxt} 
 \ar[d] \ar[r]^{4}
 & F_{1,4}^{\bxt}F'_{4,5}{}^{\!\!\!\!\bxt}
 F'_{5,6}{}^{\!\!\!\!\bxt}F_{6,9}^{\bxt} \ar[d] \\
F_{1,4}^{\bxt}F_{4,5}^{\bxt}F_{5,8}^{\bxt}F_{8,9}^{\bxt} \ar[d] \ar[r]^{8\quad}
 & F_{1,4}^{\bxt}F_{4,5}^{\bxt}I^{\bxt}F'_{5,8}{}^{\!\!\!\!\bxt}F_{8,9}^{\bxt} 
 \ar[r]^{4}
 & F_{1,4}^{\bxt}F'_{4,5}{}^{\!\!\!\!\bxt}
 F'_{5,8}{}^{\!\!\!\!\bxt}F_{8,9}^{\bxt} \ar[d] \\
F_{1,4}^{\bxt}F_{4,7}^{\bxt}F_{7,8}^{\bxt}F_{8,9}^{\bxt} \ar@{=}[rr]
 & & F_{1,4}^{\bxt}F_{4,7}^{\bxt}F_{7,8}^{\bxt}F_{8,9}^{\bxt} 
}
\]
The horizontal maps marked as 6 (resp.~2, 4, 8) are the obvious ones
induced by $\boldC_{I,F'_{5,6}}^{-1}$
(resp. $\boldC_{F_{2,5},I}$, $\boldC_{F_{4,5},I}$, $\boldC_{I,F'_{5,8}}^{-1}$).
The vertical maps are all induced by $\bbeta_{-}$ with an obvious choice for
the subscript in each case (corresponding to one of the unit squares 
in $\sS$ or $\sS'$). 

Commutativity of all the rectangles are proved using \ref{lem:canon4}
or functoriality considerations. We concentrate only on 
$\sbsq_1,\sbsq_2$ and $\sbsq_3$ by way of illustration.

In $\sbsq_1$, we may cancel the common factors 
$F_{1,2}^{\bxt}$ on the left and $F_{6,9}^{\bxt}$ on the 
right from each vertex, and so we reduce to checking commutativity of
the following diagram.
\[
\begin{CD}  
F_{2,3}^{\bxt}F_{3,6}^{\bxt}  @>>> 
 F'_{2,5}{}^{\!\!\!\!\bxt}F'_{5,6}{}^{\!\!\!\!\bxt} \\
@VVV @VVV \\
F_{2,5}^{\bxt}F_{5,6}^{\bxt} @>>> 
 F_{2,5}^{\bxt}I^{\bxt}F'_{5,6}{}^{\!\!\!\!\bxt}
\end{CD}
\]
This diagram commutes by \ref{lem:canon4}(ii) using
$G_1 = F_{2,3}, F_1 = F_{2,5}, F = F_{3,6}, G = F_{5,6}$
and $F_1' = F_{2,5}', G' = F_{5,6}'$.

In $\sbsq_3$, we reduce to checking commutativity of the 
following diagram
\[
\begin{CD}  
F_{5,6}^{\bxt}F_{6,9}^{\bxt} @>>> 
 I^{\bxt}F'_{5,6}{}^{\!\!\!\!\bxt}F_{6,9}^{\bxt} \\
@VVV @VVV \\
F_{5,8}^{\bxt}F_{8,9}^{\bxt} @>>>
 I^{\bxt}F'_{5,8}{}^{\!\!\!\!\bxt}F_{8,9}^{\bxt} 
\end{CD}
\]
and we conclude by \ref{lem:canon4}(i) using 
$G_1 = F_{5,6}, F_1 = F_{5,8}, F = F_{6,9}, G = F_{8,9}$
and $F_1' = F_{5,8}', G_1' = F_{5,6}'$.

Finally, $\sbsq_2$ commutes for functorial reasons.
\end{proof}

\begin{arem}
\label{rem:canon5a}
We have given a longer proof above, because a good portion of it is used for 
tackling \ref{lem:canon6} below. 
A quicker way to prove \ref{lem:canon5} would be to note that since 
all the vertical maps in~$\sS$ and~$\sS'$ have the same label 
(which is~$\sfP$), one can use vertical transitivity 
of base-change for any pair of squares stacked vertically. Now all 
the vertical composites of~$\sS$ and~$\sS'$, 
including the ones in the middle, viz.,  
$|F_{2,5}\star F_{5,8}|$ and~$|F_{2,5}'\star F_{5,8}'|$, are equal. 
Thus \ref{lem:canon5} follows. 
\end{arem}

The next situation involves the following two cartesian 
diagrams $\sS$ and $\sS'$ differing only at the middle vertex.
\[
\begin{CD}
X_1 @>>> X_2 \\
@VVV @VVV  \\
X_4 @>>> X_5 @>>> X_6 \\
@VVV @VVV @VVV \\
X_7 @>>> X_8 @>>> X_9
\end{CD}
\qquad \qquad \qquad
\begin{CD}
X_1 @>>> X_2 \\
@VVV @VVV \\
X_4 @>>> X_5' @>>> X_6 \\
@VVV @VVV @VVV \\
X_7 @>>> X_8 @>>> X_9
\end{CD}
\]
We adopt the same convention for naming the maps as before.
The properties satisfied by $\sS$ and $\sS'$ are the same as before 
except for the following additional condition: We assume that 
$X_2 = X_6$ (each being henceforth given the common name $X$),
and that $|F_{2,5}\star F_{5,6}| = 1_X = |F_{2,5}'\star F_{5,6}'|$.

Under these assumptions, we obtain two natural candidates for an
isomorphism of the following form, one via the cartesian squares 
of $\sS$ and the fundamental isomorphism $\bphi_{F_{2,5}, F_{5,6}}$, 
and the other via $\sS'$ \ldots 
\[
F_{1,2}^{\bxt}F_{6,9}^{\bxt}  = F_{1,2}^{\bxt} \oneD{X} F_{6,9}^{\bxt} \iso 
F_{1,4}^{\bxt}F_{4,7}^{\bxt}F_{7,8}^{\bxt}F_{8,9}^{\bxt}
\]

\begin{alem}
\label{lem:canon6}
The above-mentioned two natural candidates are the same.
\end{alem}

\begin{proof}
It suffices to check commutativity of the following replacement 
of the diagram in the proof of \ref{lem:canon5}.
The topmost rectangle there, $\sbsq_1$, has been replaced by 
$\widetilde{\sbsq}$. The remaining portion remains the same and 
only its outer border is shown below.
\[
\xymatrix{
F_{1,2}^{\bxt}F_{6,9}^{\bxt} 
 \ar[d]^{\hspace{10em} \widetilde{\sbsq}} \ar@{=}[rr]
 & & F_{1,2}^{\bxt}F_{6,9}^{\bxt} \ar[d] \\
F_{1,2}^{\bxt}F_{2,5}^{\bxt}F_{5,6}^{\bxt}F_{6,9}^{\bxt} 
 \ar[d]^{\hspace{2em}
 \text{use lower part of diagram in proof of \ref{lem:canon5}}} 
 \ar[r]_{6\quad}
 & F_{1,2}^{\bxt}F_{2,5}^{\bxt}I^{\bxt}F'_{5,6}{}^{\!\!\!\!\bxt}F_{6,9}^{\bxt} 
 \ar[r]_{2} & F_{1,2}^{\bxt}F'_{2,5}{}^{\!\!\!\!\bxt}
 F'_{5,6}{}^{\!\!\!\!\bxt}F_{6,9}^{\bxt} \ar[d] \\
F_{1,4}^{\bxt}F_{4,7}^{\bxt}F_{7,8}^{\bxt}F_{8,9}^{\bxt} \ar@{=}[rr]
 & & F_{1,4}^{\bxt}F_{4,7}^{\bxt}F_{7,8}^{\bxt}F_{8,9}^{\bxt} 
}
\]
Commutativity of $\widetilde{\sbsq}$, after the canceling of the 
common factor $F_{1,2}^{\bxt}F_{6,9}^{\bxt}$ from each vertex,
reduces to checking that of the following.
\[
\xymatrix{
\mathbf{1}_{\cD_{X}} \ar@{=}[rr] & & \mathbf{1}_{\cD_{X}} \\
F_{2,5}^{\bxt}F_{5,6}^{\bxt} \ar[u] \ar[r] & 
 F_{2,5}^{\bxt}I^{\bxt}F'_{5,6}{}^{\!\!\!\!\bxt} \ar[r] &
 F'_{2,5}{}^{\!\!\!\!\bxt} F'_{5,6}{}^{\!\!\!\!\bxt} \ar[u]
}
\]
This diagram commutes because of compatibility of the fundamental 
isomorphism with base-change (Lemma~\ref{lem:cafpols4}(ii)) 
in the context of the following diagram. 
(Recall that $X = X_2 = X_6$.)

\[
\begin{CD}
X_2 @>>> X_5 @>>> X_6 \\
@| @VV{I}V @| \\
X_2 @>>> X_5' @>>> X_6 
\end{CD}
\]
\end{proof}

For the next two situations we look at truncated versions of 
diagrams encountered in the previous two results.
Consider the following two conjugate pairs of cartesian diagrams.
\[
\begin{CD}
X_2 \\ @VVV \\ X_5 @>>> X_6 \\ @VVV @VVV \\ X_8 @>>> X_9
\end{CD}
\qquad
\begin{CD}
X_2 \\ @VVV \\ X_5' @>>> X_6 \\ @VVV @VVV \\ X_8 @>>> X_9
\end{CD}
\qquad\qquad
\begin{CD}
X_2 @>>> X_3 \\ @VVV @VVV \\ X_5 @>>> X_6 \\ @VVV @VVV \\ X_8 @>>> X_9
\end{CD}
\qquad
\begin{CD}
X_2 @>>> X_3 \\ @VVV @VVV \\ X_5' @>>> X_6 \\ @VVV @VVV \\ X_8 @>>> X_9
\end{CD}
\]
Let us call the diagrams $\sS_1,\sS_1', \sS_2,\sS_2'$ respectively 
starting from the left. 
We use $F_{i,j}, F'_{i,j}$ for maps in $\sS_1,\sS_1'$ and 
$G_{i,j}, G'_{i,j}$ for the ones in $\sS_2,\sS_2'$. In $\sS_1,\sS_1'$,
the assumptions are the same as that from \ref{lem:canon6}, 
(in particular, $X_2 = X_6 = X$, etc. \ldots) while in $\sS_2,\sS_2'$
we use the ones from \ref{lem:canon5}. 

One point of departure from the earlier cases is that here we shall 
also look at the composite vertical maps. Thus 
$F_{2,8}$, $G_{3,9}$ and $G_{2,8} = G_{2,8}'$ are the obvious composite
maps with label $\sfP$.

Proceeding as before, the above diagrams lead to 
two conjugate natural candidates each for isomorphisms of the following form.
\[
F_{6,9}^{\bxt} = \oneD{X}F_{6,9}^{\bxt}\iso F_{2,8}^{\bxt}F_{8,9}^{\bxt} 
\qquad \qquad
G_{2,3}^{\bxt}G_{3,9}^{\bxt} \iso G_{2,8}^{\bxt}G_{8,9}^{\bxt}
\]

\begin{alem}
\label{lem:canon7}
In each of the two cases above, the conjugate isomorphisms are the same.
\end{alem}
\begin{proof}
See proofs of \ref{lem:canon5} and \ref{lem:canon6}. Also see 
Remark \ref{rem:canon5a}.
\end{proof}

We are finally in a position to prove the main result of this subsection.

\begin{aprop}
\label{prop:canon8}
The fundamental isomorphism of \eqref{eq:psi3} is canonical, i.e., does not 
depend on the choice of the staircase. 
\end{aprop}
\begin{proof}
Suppose $\cS$ and $\cS'$ are two staircase diagrams 
based on $\sigma$. We proceed by successively
modifying $\cS$, one vertex at a time, till it is transformed into~$\cS'$.
The starting point for these modifications
is the southeastern end of the staircase. Suppose
$\cS$ is represented by the diagram in~\eqref{eq:psi1a} and the conjugate 
staircase~$\cS'$ is represented by a corresponding $(-)'$ version. 
Then we consider the vertex~$Y_{n,n-1}$ of~$\cS$ and replace
it by~$Y_{n,n-1}'$, while also modifying the arrows coming in and 
out of it via the unique isomorphism  
$Y_{n,n-1} \cong Y_{n,n-1}'$ arising out of the fiber-product property. 
If $n>2$, then by \ref{lem:canon6} we see that both $\cS$ and 
the newly modified staircase
give rise to the same isomorphism and if $n=2$ we use \ref{lem:canon7}
for the same conclusion. Thus, without loss of generality, we may assume 
that~$Y_{n,n-1}$ and the arrows going out of it are all shared
by~$\cS$ and~$\cS'$. (This automatically forces 
$V_{n-1,n-1} = V_{n-1,n-1}'$ too.) 
Next, we replace $Y_{n,n-2}$. Arguing the same way
as above we continue westward till the bottom row is assumed to be shared  
by~$\cS$ and~$\cS'$. Then we start with~$Y_{n-1,n-2}$. Continuing
in this fashion, the proposition is proved.  
\end{proof}

\begin{adefi}
\label{def:canon9}
For any sequence $\sigma$ such that $|\sigma|$ is an identity map,
we define $\Phi_{\sigma}$ to be $\Phi_{\sigma}(\cS)$ of \eqref{eq:psi3}
for any choice of a staircase $\cS$ based on $\sigma$.
\end{adefi}

The following lemma is used in \S\ref{sec:proofsIII}.

\begin{alem}
\label{lem:canon10}
Let $\sigma \colon X \Lra X$ be a sequence such that $|\sigma| = 1_X$ and
set $J \set (1_X, \sfP)$. Then via the identification 
$J^{\bxt}\sigma^{\bxt} = \sigma^{\bxt}$, we have 
$\Phi_{J\star\sigma} = \Phi_{\sigma}$. 
\end{alem}
\begin{proof}
Let $\cS$ be a staircase based on $\sigma$.
Then a staircase based on $J \star \sigma$ can be 
obtained as shown below where $V= V_1(\cS)$ is the leftmost column of $\cS$
and $V'= V, J'' = J' = J$.
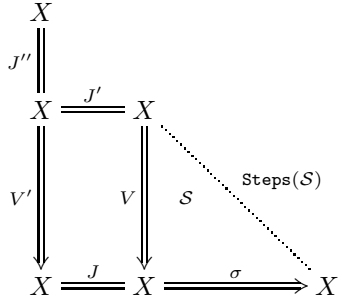
\[
\xymatrix{
X \ar@{=}[d]_{J''} \\ 
X \ar@{=}[r]^{J'} \ar@{=>}[dd]_{V'} & X    
 \ar@{=>}[dd]^{\hspace{1em}\cS}_{V} 
 \ar@{.}[rrdd]^{\texttt{Steps}(\cS)} \\ 
\\
X \ar@{=}[r]^{J} & X \ar@{=>}[rr]^{\sigma} & &  X
}
\]
By pseudofunctoriality of $(-)^{\times}$ 
we conclude that the isomorphism $\theta$ given by composing the following 
ones
\[
\oneD{X} = J^{\bxt} \iso J''{}^{\bxt}V'{}^{\bxt}J^{\bxt} \iso 
J''{}^{\bxt}J'{}^{\bxt}V^{\bxt} \iso V^{\bxt}
\]
is the same as the obvious one resulting from composing all the 
(necessarily $\sfP$-labeled) maps in $V$. The desired result now follows 
by comparing the definition of~$\Phi_{\sigma}$ with that of 
$\Phi_{J \star \sigma}$.
\end{proof}

\section{Proofs II (the cocycle condition)}
\label{sec:proofsII}
This section contains some of the core results used in the proof of 
the abstract output results of the paper. We begin with the definition
of $\Psi_{-,-}$ that gives an isomorphism 
$\sigma_1^{\bxt} \iso \sigma_2^{\bxt}$ whenever $|\sigma_1| = |\sigma_2|$.
Then we show that $\Psi_{-,-}$ satisfies the cocycle condition.
The proof of this is somewhat long and occupies most of this section.
We reduce the proof into a few lemmas, which are then verified 
over different subsections.

\subsection{The definition of $\Psi_{-,-}$}
\label{subsec:psi}
The isomorphism $\Psi_{-,-}$ that we define here is 
easy to obtain using results defined from the previous section.
Its canonicity however will not be immediately
obvious owing to dependence on certain fibered-product diagrams 
needed in the construction. It is only after proving the cocycle condition, 
the topic of the next few subsections, that we get to 
conclude canonicity.

\ssss
\label{sssec:psi5z}
Let $\sigma_1, \sigma_2$ be two labeled sequences such that 
$|\sigma_1| = |\sigma_2|$. Then we want to define an isomorphism 
\[
\Psi_{\sigma_1, \sigma_2} \colon \sigma_2^{\bxt} \iso \sigma_1^{\bxt}.
\]
This is achieved as follows. 
Let $X$ be the source of $\sigma_i$ and~$Y$ the target.
Consider the following diagram where the square $\sS$ of double arrows
is cartesian (\S\ref{sssec:cafpols5}) and~$\delta$ is the diagonal map, which, 
by \S\ref{subsec:input}[A](ii), is in~$\sfP$.  
\[
\xymatrix{
X \ar[rr]^{\Delta\>=\>(\delta, \sfP)\qquad} 
 & & X \times_Y X \ar@{=>}[d]_{\sigma_1'}^{\hspace{3.6em} \sS} 
 \ar@{=>}[rr]^{\quad\sigma_2'} & & X \ar@{=>}[d]^{\sigma_1}  \\
& & X \ar@{=>}[rr]_{\sigma_2} & &  Y 
}
\]
We define $\Psi_{\sigma_1, \sigma_2}^{\sS}$ 
to be the composition of the following sequence of isomorphisms.
\begin{align}
\sigma_2^{\bxt} 
 \xrightarrow[\ref{def:canon9}]{\text{via } 
 (\Phi_{\Delta \star \sigma_1'})^{-1}} 
 (\Delta \star \sigma_1')^{\bxt}\sigma_2^{\bxt}
 &= \Delta^{\!\bxt} (\sigma_1' \star \sigma_2)^{\bxt} \label{eq:psi5} \\
&\cong \Delta^{\!\bxt} (\sigma_2' \star \sigma_1)^{\bxt} \qquad 
 \text{(via $\bbeta_{\sS}^{-1}$)} \notag \\
&= (\Delta \star \sigma_2')^{\bxt} \sigma_1^{\bxt}
 \xrightarrow[\ref{def:canon9}]{\text{via }\Phi_{\Delta \star \sigma_2'}}  
 \sigma_1^{\bxt} \notag
\end{align}

\begin{arem}
\label{rem:psi6}
As in the case of the generalized fundamental isomorphism, 
we would like to say that $\Psi_{\sigma_1, \sigma_2}^{\sS}$
does not depend on the choice of $\sS$ which is determined 
only up to isomorphism. However, the approach that we took 
in \S\ref{subsec:canon} 
cannot be adapted here because there we exploited the fact 
that vertical maps of a staircase are all $\sfP$-labeled. 
What is perhaps needed is a more elaborate set of 
compatibilities over the category $\sfI$ of isomorphisms
in $\sfC$. Instead, we resort to the following strategy. We show that 
the cocycle condition for $\Psi_{-,-}$, a property that we are 
interested in for various reasons, holds even at the level of the 
extrinsically defined $\Psi_{-,-}^{-}$.
Then, using the cocycle condition we deduce that the definition is indeed 
canonical. 
The price we pay for this indirect approach 
is of having to retain the cumbersome notation
such as $\Psi_{\sigma_1, \sigma_2}^{\sS}$ instead of 
$\Psi_{\sigma_1, \sigma_2}$, and in particular, of keeping track of the
choice of~$\sS$ whenever necessary, till the end of the proof of the 
cocycle condition. 
\end{arem}

\ssss
\label{sssec:psi7}
The following particular cases of \eqref{eq:psi5} are noteworthy.

(i) Let $\sigma$ be a sequence such that $f \set |\sigma| \in \sfP$. 
Set $F \set (f, \sfP)$. Then one obtains an isomorphism
(for any available choice of $\sS$)
$\Psi_{F, \sigma}^{\sS} \colon \sigma^{\bxt} \iso f^{\times}$.

(ii) Similarly, if we assume instead, that $f = |\sigma| \in \sfO$, then, 
with $F \set (f, \sfO)$, we obtain an isomorphism
$\Psi_{F, \sigma}^{\sS} \colon \sigma^{\bxt} \iso f^{\ssbox}$.  

It can be shown that $\Psi_{-,-}$ recovers the fundamental isomorphism,
in that if we let $f$ in (i) above to be an identity map, then 
for a suitable choice of~$\sS$, (and hence for any choice) 
the isomorphism in~(i) gives the
fundamental isomorphism~$\Phi_{\sigma}$. However, we postpone the 
verification of this fact to the next section.

\subsection{The cocycle condition}
\label{subsec:cocycle}

Perhaps the most important result lying at the technical heart of 
this paper is the following.

\begin{athm}
\label{thm:cocycle1}
The isomorphism $\Psi_{-,-}^{-}$ of \eqref{eq:psi5} satisfies the 
cocycle condition, i.e., if $\sigma_1, \sigma_2, \sigma_3$ are labeled 
sequences such that 
$|\sigma_1| = |\sigma_2| = |\sigma_3|$, then 
the following condition holds
\[
\Psi_{\sigma_3,\sigma_2}^{\sT}\Psi_{\sigma_2,\sigma_1}^{\sS} 
= \Psi_{\sigma_3,\sigma_1}^{\sU}
\]
for any available choice of fibered product diagrams $\sS, \sT, \sU$. 
\end{athm}

\begin{arem}
\label{rem:cocycle1a}
We use the cocycle condition as stated above in two stages.
In the first one, the presence of superscripts of $\Psi$ is important.
This is really catered towards proving that 
$\Psi_{\sigma_1,\sigma_2}^{\sS}$ of \eqref{eq:psi5}
is independent of the choice of~$\sS$, which would then give us a canonical
definition of $\Psi_{\sigma_1,\sigma_2}$ (\ref{def:conseq2}). 
Once this is accomplished,
so that we have the luxury of dropping the superscripts for $\Psi$,
the true cocycle condition emerges. The trimmed version leads to 
many interesting consequences including the existence of a pseudofunctor
over $\sfQ = \ov{\{\sfO,\sfP\}}$ 
that generalizes $(-)^{\ssbox}$ and $(-)^{\times}$.
\end{arem}

Before going into the proof of \ref{thm:cocycle1}, we make some 
general remarks on fibered products.

\ssss
\label{sssec:cocycle1b}
For the remainder of this paper, we resort to 
the following somewhat ambiguous convention. 
For any cartesian square $\sfr \set (F_1,F_2,F_1',F_2')$
we use $\bbeta_{F_1,F_2}$ in place of $\bbeta_{\sfr}$.
To eliminate any ambiguity, this notation is employed 
only in situations where all the components of $\sfr$ 
are already chosen unambiguously. 
Note that under such circumstances
we may write $\bbeta_{F_1,F_2} = \bbeta_{F_2,F_1}^{-1}$.

A similar convention also applies in the more general situation of products
of sequences. Thus given two sequences $\sigma_1,\sigma_2$ having the same 
target and given a choice of a fibered square $\sS$ on $(\sigma_1, \sigma_2)$,
we use $\bbeta_{\sigma_1, \sigma_2}$ 
(resp.~$\bbeta_{\sigma_2, \sigma_1}$) to denote $\bbeta_{\sS}$
(resp.~$\bbeta_{\sS}^{-1}$).

\ssss
\label{sssec:cocycle1c}
Given three maps $X_i \xto{F_i} Y$ for $i = 1,2,3$, a
cartesian cube (or fibered cube) on $F_1$, $F_2$ and $F_3$ has the 
obvious meaning (each of the six faces is a cartesian square)
and is drawn as follows
where the symbol $\bullet$ stands for the appropriate fibered products. 
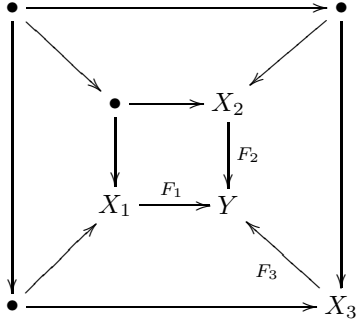
\[
\xymatrix{
\bullet \ar[ddd] \ar[dr] \ar[rrr] & & & \bullet \ar[ddd] \ar[dl] \\
& \bullet \ar[d] \ar[r] & X_2 \ar[d]^{F_{2}} & \\
& X_1 \ar[r]^{F_1} & Y  & \\
\bullet \ar[ur] \ar[rrr] & & & X_3  \ar[ul]^{F_3}
}
\]
Unlike the case of cartesian squares, here we avoid
imposing an ordering on the edges since our interest
is mainly on base-change isomorphisms associated to squares, 
for which the convention in \ref{sssec:cocycle1b} above suffices.

The existence of such a cartesian cube on $F_1, F_2$ and $F_3$ 
is a formal exercise on fibered products. In fact, the 
following also holds. For any choice of cartesian squares
$\sfr_{12}$ (resp.~$\sfr_{13}$, resp.~$\sfr_{23}$) on $F_1,F_2$
(resp.~$F_1,F_3$, resp.~$F_2,F_3$) there is a cartesian cube
on $F_1$, $F_2$ and $F_3$ which has $\sfr_{12}$, $\sfr_{13}$ and
$\sfr_{23}$ for three of its faces. 
 
The notion of a cartesian cube also generalizes to sequences. 
Thus, given three sequences $\sigma_i \colon X_i \to Y$ for $i = 1,2,3$, 
a cartesian cube on $\sigma_1$, $\sigma_2$ and $\sigma_3$
has the obvious meaning. Moreover, for any choice of 
squares $\sS_{1,2}, \sS_{1,3}, \sS_{2,3}$ obtained as fibered product of
$\sigma_1$ and $\sigma_2$, etc., there is a cartesian cube having
these squares for three of its faces. As in the case of cartesian
squares over sequences, we usually resort to drawing a cartesian cube 
only through its outermost edges, each represented by a double arrow.

Coming back to the topic of proving Theorem \ref{thm:cocycle1}, we now state
the three lemmas which are used in it. 

We call the first of the following lemmas as the cube lemma.

\begin{alem}
\label{lem:cocycle2}
Let $\sigma_1, \sigma_2, \sigma_3$ be labeled sequences 
having the same target. Let~$X_i$ be the source of~$\sigma_i$
and~$Y$ the target. Consider a cartesian cube on $\sigma_1,\sigma_2$
and $\sigma_3$ as follows where
$Z$ is the product of $X_1, X_2, X_3$ over~$Y$.  
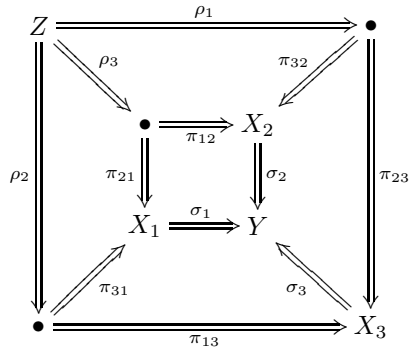
\[
\xymatrix{
Z \ar@{=>}[ddd]_{\rho_2}\ar@{=>}[dr]^{\rho_3} 
 \ar@{=>}[rrr]^{\rho_1} & & & \bullet \ar@{=>}[ddd]^{\pi_{23}} 
 \ar@{=>}[dl]_{\pi_{32}} \\
& \bullet \ar@{=>}[d]_{\pi_{21}} \ar@{=>}[r]_{\pi_{12}}  
 & X_2 \ar@{=>}[d]^{\sigma_{2}} & \\
& X_1  \ar@{=>}[r]^{\sigma_1} & Y  & \\
\bullet \ar@{=>}[ur]_{\pi_{31}} \ar@{=>}[rrr]_{\pi_{13}} & & & 
 X_3  \ar@{=>}[ul]^{\sigma_3}
}
\]
Then the following hexagon of isomorphisms commutes. 
Here each object corresponds to one of 
the six different paths from $Z$ to~$Y$ in the above diagram.
For the morphisms the omitted subscripts of the $\bbeta$'s  
are the obvious ones.
\[
\xymatrix{
& \rho_3^{\bxt}\pi_{21}^{\bxt}\sigma_1^{\bxt}  
 \ar[rr]^{\rho_3^{\bxt}(\bbeta)} & & 
 \rho_3^{\bxt}\pi_{12}^{\bxt}\sigma_2^{\bxt} 
 \ar[dr]^{\bbeta(\sigma_2^{\bxt})}  \\
\rho_2^{\bxt}\pi_{31}^{\bxt}\sigma_1^{\bxt} 
 \ar[ur]^{\bbeta(\sigma_1^{\bxt})} \ar[dr]_{\rho_2^{\bxt}(\bbeta)} & & & & 
 \rho_1^{\bxt}\pi_{32}^{\bxt}\sigma_2^{\bxt} \\
& \rho_2^{\bxt}\pi_{13}^{\bxt}\sigma_3^{\bxt} 
 \ar[rr]_{\bbeta(\sigma_3^{\bxt})} & & 
 \rho_1^{\bxt}\pi_{23}^{\bxt}\sigma_3^{\bxt} \ar[ur]_{\rho_1^{\bxt}(\bbeta)}
}
\]
\end{alem}

\begin{alem}
\label{lem:cocycle3}
Let $X \stackrel{\sigma_1\;}{\Lra} Y \stackrel{\sigma_2\;}{\Lra} X$ 
be sequences such that $|\sigma_2\star\sigma_1| = 1_X$,
so that, by \textup{[A](i) and [A](ii) of \S\ref{subsec:input},} 
$|\sigma_1| \in \sfP$. 
Set $G \set (|\sigma_1|, \sfP)$. 
Then the following diagram of isomorphisms commutes for any 
available choice of $\sS$. 
\[
\xymatrix{
\sigma_1^{\bxt}\sigma_2^{\bxt} 
 \ar[rr]^{\Psi_{G,\sigma_1}^{\sS}(\sigma_2^{\bxt})} 
 \ar[rd]_{\Phi_{\sigma_1\star\sigma_2}} & & 
 G^{\bxt}\sigma_2^{\bxt} \ar[ld]^{\;\;\Phi_{G\star\sigma_2}} \\
& \xyoneD{X}
}
\]
\end{alem}

\begin{alem}
\label{lem:cocycle4}
Consider the following diagram 
\[
\xymatrix{
X \ar@{=>}[r]^{\sigma_1} & A \ar@{=>}[d]_{\sigma_3} \ar@{=>}[r]^{\sigma_2} 
 & B \ar@{=>}[d]^{\sigma_4}  \\
& C \ar@{=>}[r]_{\sigma_5} & D \ar@{=>}[r]_{\sigma_6} & X 
}
\]
where the square is cartesian and 
$|\sigma_1 \star \sigma_3 \star \sigma_5 \star \sigma_6| = 1_X$.
Then the following diagram of isomorphisms commutes 
where the~$\Phi$'s have the obvious subscripts.
\[
\xymatrix{
\sigma_1^{\bxt} \sigma_3^{\bxt} \sigma_5^{\bxt} \sigma_6^{\bxt}
 \ar[rr]^{\text{via } \;\bbeta_{\sigma_5, \sigma_4}} \ar[rd]_{\Phi} 
 & & \sigma_1^{\bxt} \sigma_2^{\bxt} \sigma_4^{\bxt} \sigma_6^{\bxt} 
 \ar[ld]^{\Phi} \\
& \mathbf{1}_{\cD_X} 
}
\]
\end{alem}


\ssss
\label{sssec:cocycle5}
Assuming \ref{lem:cocycle2}, \ref{lem:cocycle3} and \ref{lem:cocycle4},
one proves \ref{thm:cocycle1} as follows. Let $X$ be the 
source of~$\sigma_i$ and~$Y$ the target. Consider the following 
diagram where the cube of double arrows is cartesian 
having $\sS,\sT,\sU$ as three of its faces in the obvious places.
Also, $\Delta$ and $\Delta_{ij}$ are the appropriate diagonal maps
with label $\sfP$.
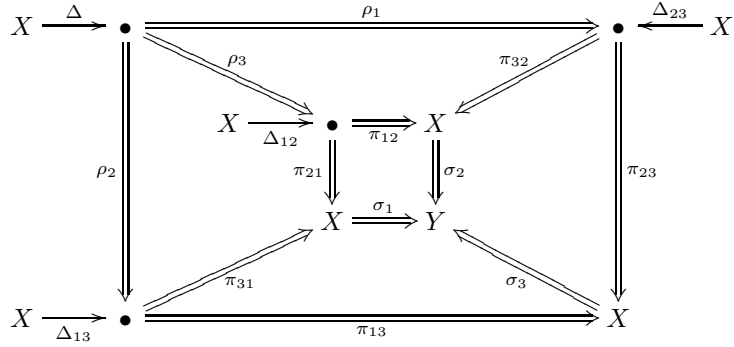
\[
\xymatrix{
X \ar[r]^{\Delta} & {\phantom{l}\!\!\bullet\!\!\phantom{l}} 
 \ar@{=>}[ddd]_{\rho_2}\ar@{=>}[drr]^{\rho_3} \ar@{=>}[rrrrr]^{\rho_1} 
 & & & & &  {\phantom{l}\!\!\bullet\!\!\phantom{l}} \ar@{=>}[ddd]^{\pi_{23}} 
 \ar@{=>}[dll]_{\pi_{32}} & \ar[l]_{\Delta_{23}} X \\
& & X \ar[r]_{\Delta_{12}} & {\phantom{l}\!\!\bullet\!\!\phantom{l}}
 \ar@{=>}[d]_{\pi_{21}} \ar@{=>}[r]_{\pi_{12}} 
 & X \ar@{=>}[d]^{\sigma_{2}} \\
& & & X  \ar@{=>}[r]^{\sigma_1} & Y & & \\
X \ar[r]_{\Delta_{13}} & {\phantom{l}\!\!\bullet\!\!\phantom{l}}
 \ar@{=>}[urr]_{\pi_{31}} \ar@{=>}[rrrrr]_{\pi_{13}} & & & & &
 X  \ar@{=>}[ull]^{\sigma_3}
}
\]
Now consider the following induced diagram of isomorphisms.
To avoid clutter, we only describe the maps in \ddag,
the rest being defined analogously. In \ddag, for the 
downward arrow on the left, $\Phi$ has subscript
$\Delta\star\rho_2\star\pi_{31}$, while for the one on
the right, the subscript is $\Delta\star\rho_2\star\pi_{13}$.
The remaining definitions are self-evident. 
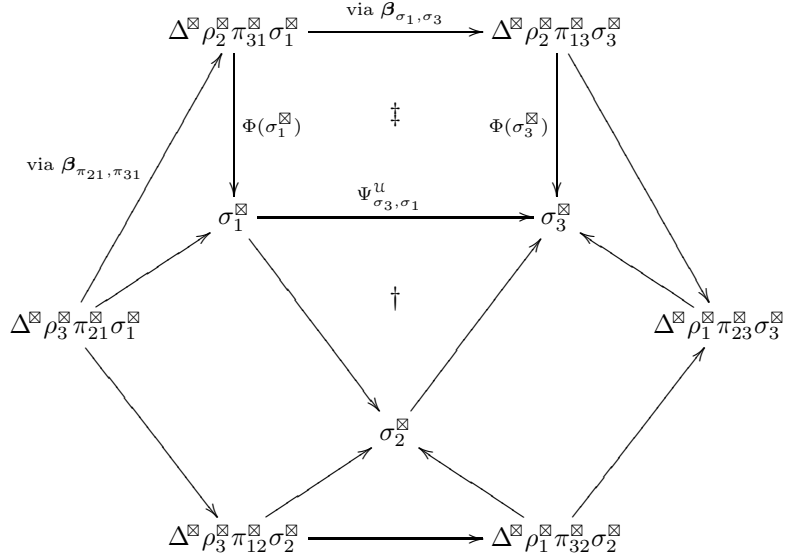
\[
\xymatrix{
& \Delta^{\!\bxt}\rho_2^{\bxt}\pi_{31}^{\bxt}\sigma_1^{\bxt}
 \ar[dd]^{\Phi(\sigma_1^{\bxt})} 
 \ar[rr]^{\text{via }\bbeta_{\sigma_1,\sigma_3}} & & 
 \ar[dd]_{\Phi(\sigma_3^{\bxt})}  \ar[rddd] 
 \Delta^{\!\bxt}\rho_2^{\bxt}\pi_{13}^{\bxt}\sigma_3^{\bxt} \\
\\
& \sigma_1^{\bxt} \ar[rdd] 
 \ar[rr]^{\makebox[0pt]{\scriptsize{$\Psi_{\sigma_3, \sigma_1}^{\sU}$}} 
 \smash{\raisebox{7ex}{\ddag}}}_{\smash{
 \raisebox{-7ex}{\dag}}} & & \sigma_3^{\bxt} \\
\makebox[0pt]{\qquad 
 $\Delta^{\!\bxt}\rho_3^{\bxt}\pi_{21}^{\bxt}\sigma_1^{\bxt}$} 
 \phantom{A^B_D} \ar[ruuu]^{\text{via } \bbeta_{\pi_{21}, \pi_{31}}} 
 \ar[ru] \ar[rdd] & & & & \ar[ul] \phantom{A^B_D} \makebox[0pt]{
 $\Delta^{\!\bxt}\rho_1^{\bxt}\pi_{23}^{\bxt}\sigma_3^{\bxt}$\qquad} \\
& & \sigma_2^{\bxt} \ar[ruu] \\
& \Delta^{\!\bxt}\rho_3^{\bxt}\pi_{12}^{\bxt}\sigma_2^{\bxt} 
 \ar[ru] \ar[rr] & & 
 \Delta^{\!\bxt}\rho_1^{\bxt}\pi_{32}^{\bxt}\sigma_2^{\bxt}
 \ar[lu] \ar[ruu]  
}
\]
The outer border, which is a hexagon, commutes, as is seen by 
first canceling~$\Delta^{\!\bxt}$ 
from each vertex and then using the cube lemma (\ref{lem:cocycle2}). 
The three outer triangles commute by~\ref{lem:cocycle4}. 
For the rectangles we use~\ref{lem:cocycle3}. For example, 
let us consider~\ddag. We expand it as follows.
The bottom row spells out the definition 
of~$\Psi_{\sigma_3, \sigma_1}^{\sU}$. The~$\Phi$'s have the obvious
subscripts and we use  
$\Psi' \set \Psi_{\Delta_{13}, \Delta\star\rho_2}^{\sV}$ 
for some fixed available choice of~$\sV$.
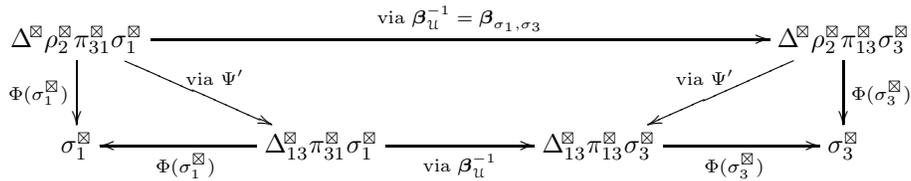
\[
\xymatrix{
\Delta^{\!\bxt}\rho_2^{\bxt}\pi_{31}^{\bxt}\sigma_1^{\bxt}
 \ar[d]_{\Phi(\sigma_1^{\bxt})} 
 \ar[rd]^{\text{via }\Psi'} 
 \ar[rrr]^{\text{via } \bbeta_{\sU}^{-1} \; = \; \bbeta_{\sigma_1,\sigma_3}} 
 & \text{\phantom{AAAAAAAAAA}} & \text{\phantom{AAAAAAAAAA}} & 
 \Delta^{\!\bxt}\rho_2^{\bxt}\pi_{13}^{\bxt}\sigma_3^{\bxt} 
 \ar[ld]_{\text{via }\Psi'} \ar[d]^{\Phi(\sigma_3^{\bxt})} \\
\sigma_1^{\bxt} & \Delta_{13}^{\!\bxt}\pi_{31}^{\bxt}\sigma_1^{\bxt}
 \ar[l]^{\Phi(\sigma_1^{\bxt})\quad} 
 \ar[r]_{\text{via }\bbeta_{\sU}^{-1}} & 
 \Delta_{13}^{\!\bxt}\pi_{13}^{\bxt}\sigma_3^{\bxt} 
 \ar[r]_{\quad\Phi(\sigma_3^{\bxt})} & \sigma_3^{\bxt}
}
\]
The trapezium in the middle commutes for functorial reasons. 
The triangle on the left, upon canceling of the common 
term~$\sigma_1^{\bxt}$ from each vertex, results in 
the following one whose commutativity
now follows from~\ref{lem:cocycle3}. 
\[
\xymatrix{
\Delta^{\!\bxt}\rho_2^{\bxt}\pi_{31}^{\bxt} 
 \ar[d]_{\Phi_{(\Delta\star\rho_2)\star\pi_{31}}} 
 \ar[rd]^{\quad\text{via }\Psi_{\Delta_{13}, \Delta\star\rho_2}^{\sV}} & 
 \text{\phantom{AAAAAAAAAA}} \\
\mathbf{1}_{\cD_{X}} & \Delta_{13}^{\!\bxt}\pi_{31}^{\bxt}  
 \ar[l]^{\Phi_{\Delta_{13}\star\pi_{31}}} 
}
\] 
A similar argument works for the triangle on the right.
Thus \ddag\; commutes. A similar proof works for the remaining 
rectangles.

It now follows that \dag\; also commutes. This proves the cocycle
condition modulo
Lemmas \ref{lem:cocycle2}, \ref{lem:cocycle3} and 
\ref{lem:cocycle4}.

\subsection{Proof of the cube lemma}
\label{subsec:cube}
This subsection is devoted to the proof of Lemma~\ref{lem:cocycle2}.

We proceed by induction on the lengths of the $\sigma_i$.
To be precise, we use the sum of the lengths for induction. 
The proof is divided into two parts. The first concerns the basis 
of induction and the second, the inductive step. 

\ssss
\label{sssec:unitcube}
The basis of induction is 
the case when $\sigma_1$, $\sigma_2$ and 
$\sigma_3$, all are assumed to have length one. 
In this case we use the single-arrow
notation for drawing the~$\sigma_i$'s. Note that two of the~$\sigma_i$'s
must have the same label. Let us assume, without loss of generality,
that~$\sigma_2$ and~$\sigma_3$
have the same label. Let~$X_4$ be the fibered product of~$X_2$ and~$X_3$
over~$Y$ and let~$\sigma_4$ be the induced map $X_4 \to Y$, having the 
same label as that of~$\sigma_2$. (Note that $X_4$ and the projection maps
$\pi_{32},\pi_{23}$ are assumed to have been specified at the outset 
as per the hypothesis of \ref{lem:cocycle2}.) 
Consider the following diagram where $\rho_4$ is the natural composite map
$Z \to X_1$ having the same label as that of~$\rho_2$. Note that 
$(\sigma_4,\sigma_1,\rho_4,\rho_1)$ is also a cartesian square.

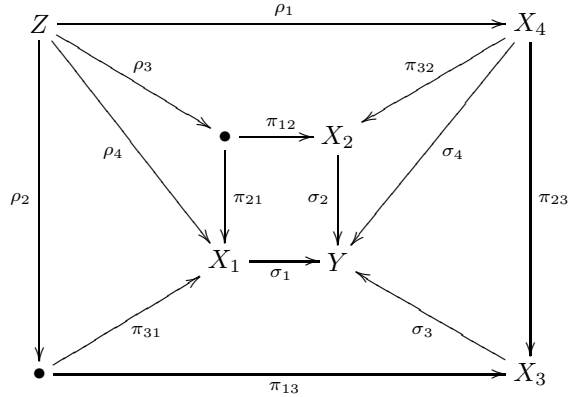
\[
\xymatrix{
Z \ar[ddd]_{\rho_2} \ar[drr]^{\rho_3} \ar[ddrr]_{\;\;\rho_4}
 \ar[rrrrr]^{\rho_1} & & & & & X_4 \ar[ddd]^{\pi_{23}} 
 \ar[dll]_{\pi_{32}} \ar[ddll]^{\sigma_4} \\
\raisebox{4ex}{ } & & \bullet \ar[d]^{\pi_{21}} \ar[r]^{\pi_{12}}  
 & X_2 \ar[d]_{\sigma_{2}}  \\
\raisebox{4ex}{ } & & X_1  \ar[r]_{\sigma_1} & Y  \\
\bullet \ar[urr]_{\pi_{31}} \ar[rrrrr]_{\pi_{13}} & & & & & 
 X_3  \ar[ull]^{\sigma_3}
}
\]
Consider the following induced diagram of isomorphisms 
where each object corresponds to one of the eight different paths
from $Z$ to $Y$ and the omitted subscripts of the $\bbeta$'s 
and the $\boldC$'s (\ref{sssec:cafpols0}) are
the obvious ones. 
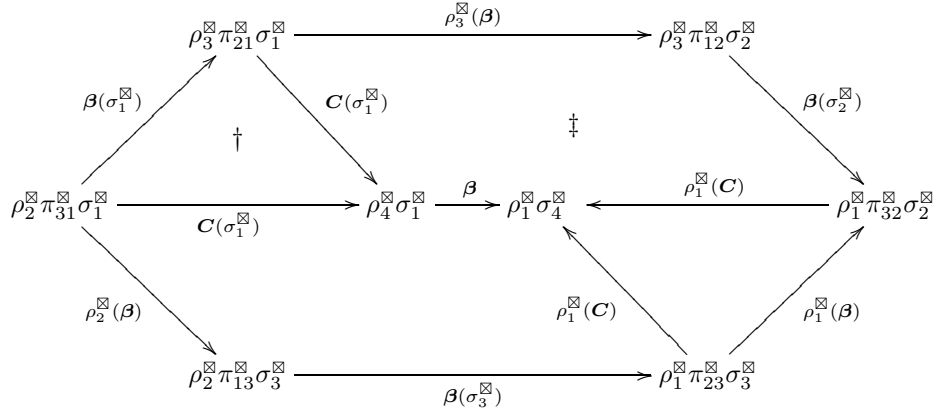
\[
\xymatrix{
& \rho_3^{\bxt}\pi_{21}^{\bxt}\sigma_1^{\bxt} 
 \ar[dr]^{\boldC(\sigma_1^{\bxt})} 
 \ar[rrr]^{\rho_3^{\bxt}(\bbeta)} & & & 
 \rho_3^{\bxt}\pi_{12}^{\bxt}\sigma_2^{\bxt} 
 \ar[dr]^{\bbeta(\sigma_2^{\bxt})}  \\
\rho_2^{\bxt}\pi_{31}^{\bxt}\sigma_1^{\bxt} \ar[rr]_{\boldC(\sigma_1^{\bxt})}  
 \ar[ur]^{\bbeta(\sigma_1^{\bxt})} \ar[dr]_{\rho_2^{\bxt}(\bbeta)} 
 & \raisebox{12ex}{\dag} &
 \rho_4^{\bxt}\sigma_1^{\bxt} \ar[r]^{\bbeta} & \rho_1^{\bxt}\sigma_4^{\bxt} 
  \smash{\raisebox{7ex}{\ddag}} & & \ar[ll]_{\rho_1^{\bxt}(\boldC)}
 \rho_1^{\bxt}\pi_{32}^{\bxt}\sigma_2^{\bxt} \\
& \rho_2^{\bxt}\pi_{13}^{\bxt}\sigma_3^{\bxt} 
 \ar[rrr]_{\bbeta(\sigma_3^{\bxt})} & & & \ar[ul]^{\rho_1^{\bxt}(\boldC)}
 \rho_1^{\bxt}\pi_{23}^{\bxt}\sigma_3^{\bxt} \ar[ur]_{\rho_1^{\bxt}(\bbeta)}
}
\]
Since $\sigma_2$, $\sigma_3$ and $\sigma_4$ have the same labels, 
by construction, this label is also shared by 
$\pi_{23},\pi_{32}, \pi_{31}, \rho_2, \rho_3, \rho_4$. Therefore,
\dag, upon the canceling of the common term~$\sigma_1^{\bxt}$
from each vertex, commutes for pseudofunctorial reasons. For \ddag,
we deduce commutativity using the transitivity property of the $\bbeta$'s
(\ref{lem:cafpols4}(iii)). The remaining subdiagrams above commute 
analogously. From the outer border of the preceding diagram, we therefore
conclude that the cube lemma holds in the case when all the~$\sigma_i$'s
have length one.

\ssss
\label{sssec:doublecube}
Now consider the case where one of the $\sigma_i$'s has 
length greater than one so that we may decompose it as a 
concatenation of two sequences; without loss of generality, let us 
assume that $\sigma_3$ factors as 
$X_3 \stackrel{\sigma_4\;}{\Lra} W \stackrel{\sigma_5\;}{\Lra} Y$.
The corresponding cartesian diagram, a double cube as shown below,
is the decomposition of the original cube into two cubes.
(Once again, recall that all the faces of this double cube are 
assumed to have been already specified by hypothesis.)

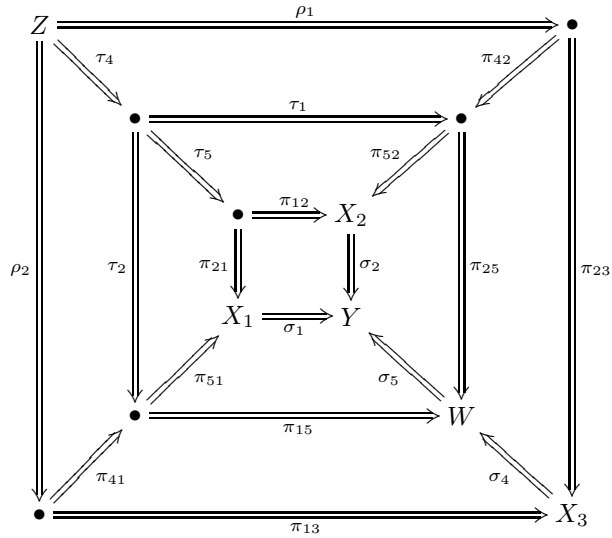
\[
\xymatrix{
Z \ar@{=>}[ddddd]_{\rho_2} \ar@{=>}[rd]^{\tau_4} \ar@{=>}[rrrrr]^{\rho_1} 
 & & & & & \bullet \ar@{=>}[ld]_{\pi_{42}} \ar@{=>}[ddddd]^{\pi_{23}} \\
& \bullet \ar@{=>}[ddd]_{\tau_2} \ar@{=>}[rd]^{\tau_5} \ar@{=>}[rrr]^{\tau_1} 
 & & & \bullet \ar@{=>}[ld]_{\pi_{52}} \ar@{=>}[ddd]^{\pi_{25}} \\
& & \bullet \ar@{=>}[d]_{\pi_{21}} \ar@{=>}[r]^{\pi_{12}} 
 & X_2 \ar@{=>}[d]^{\sigma_2} \\
& & X_1 \ar@{=>}[r]_{\sigma_1} & Y \\
& \bullet \ar@{=>}[ru]_{\pi_{51}} \ar@{=>}[rrr]_{\pi_{15}} 
 & & & W \ar@{=>}[lu]^{\sigma_5} \\
\bullet \ar@{=>}[ru]_{\pi_{41}} \ar@{=>}[rrrrr]_{\pi_{13}}  
 & & & & & X_3 \ar@{=>}[lu]^{\sigma_4}
}
\]
Note that the following hold.
\[
\rho_3 = \tau_4\star\tau_5, \qquad \pi_{31} = \pi_{41}\star\pi_{51}, \qquad
\pi_{32} = \pi_{42}\star\pi_{52}, \qquad \sigma_3 = \sigma_4\star\sigma_5
\]
Consider the following diagram of isomorphisms where each
unframed object corresponds to one of the twelve different paths from~$Z$ 
to~$Y$ in the preceding diagram. For convenience, we have dropped
the common superscript ${}^{\bxt}$ from each term. 
The morphisms are determined using the following thumb rule: 
to determine any drawn morphism between two objects, first cancel
the common terms from the two objects involved. 
The remaining terms then uniquely determine 
a base-change isomorphism~$\bbeta_{-,-}$. The morphism under 
consideration is then the obvious one induced by~$\bbeta_{-,-}$.

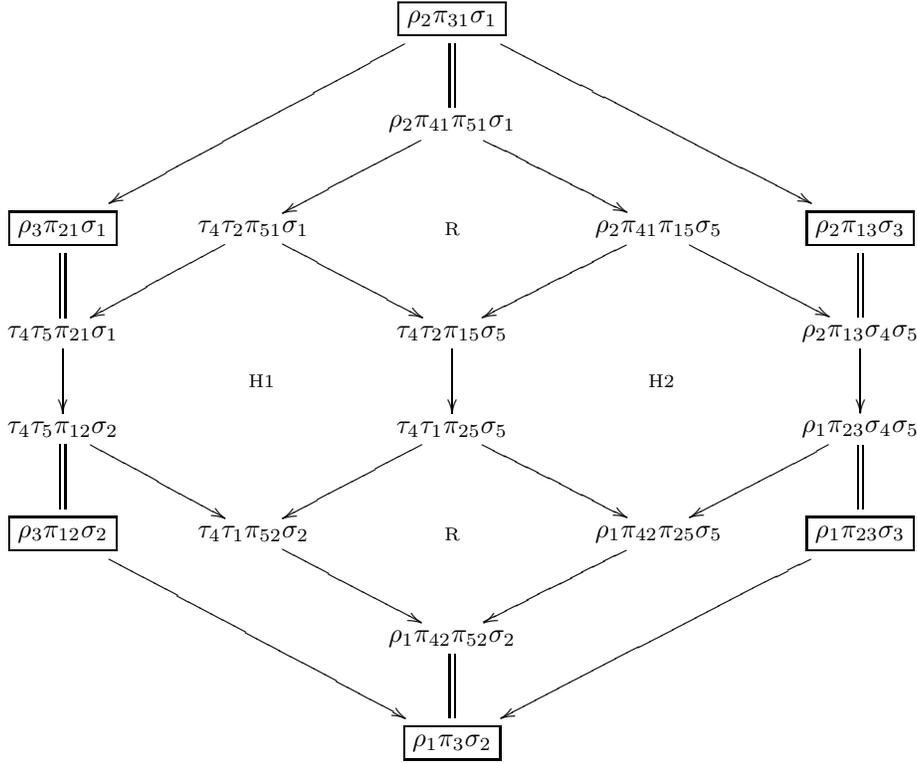
\[
\xymatrix{
& & \framebox{$\rho_2\pi_{31}\sigma_1$} \ar[ddll] \ar@{=}[d] \ar[ddrr] \\
& & \rho_2\pi_{41}\pi_{51}\sigma_1 \ar[dl] \ar[dr]  \\
\framebox{$\rho_3\pi_{21}\sigma_1$} \ar@{=}[d] & \tau_4\tau_2\pi_{51}\sigma_1
 \ar[dl] \ar[dr] & \text{\scriptsize R} & 
 \rho_2\pi_{41}\pi_{15}\sigma_5 \ar[dl] \ar[dr] & 
 \framebox{$\rho_2\pi_{13}\sigma_3$} \ar@{=}[d] \\
\tau_4\tau_5\pi_{21}\sigma_1 \ar[d]^{\hspace{7em}\!\!\text{H1}} 
 & & \tau_4\tau_2\pi_{15}\sigma_5 \ar[d]
 & & \rho_2\pi_{13}\sigma_4\sigma_5 \ar[d]_{\text{H2}\!\!\hspace{7em}} \\
\tau_4\tau_5\pi_{12}\sigma_2 \ar@{=}[d] \ar[dr] & & 
 \tau_4\tau_1\pi_{25}\sigma_5 \ar[dl] \ar[dr] & & 
 \rho_1\pi_{23}\sigma_4\sigma_5\ar[dl] \ar@{=}[d] \\
\framebox{$\rho_3\pi_{12}\sigma_2$} \ar[ddrr] & 
 \tau_4\tau_1\pi_{52}\sigma_2 \ar[dr] 
 & \text{\scriptsize R} & \rho_1\pi_{42}\pi_{25}\sigma_5 \ar[dl] & 
 \framebox{$\rho_1\pi_{23}\sigma_3$} \ar[ddll] \\
& & \rho_1\pi_{42}\pi_{52}\sigma_2 \ar@{=}[d] \\
& & \framebox{$\rho_1\pi_{3}\sigma_2$} 
}
\]
The two rhombuses, denoted by R, commute for functorial 
reasons. The two hexagons, H1 and H2, may be assumed to 
commute by the induction hypothesis; for H1 we use the cube 
constructed from $\pi_{15}$, $\pi_{25}$ and $\sigma_4$,
while for H2 we use the cube constructed from 
$\sigma_1$, $\sigma_2$ and $\sigma_5$. The remaining four outer
subdiagrams commute by definition of the
base-change isomorphism (cf.~\ref{sssec:invariance}); 
the equalities in these subdiagrams
are obtained by comparing the double cube with the 
original one. 
 
The outer border of the preceding diagram comprising the 
framed vertices and the induced maps
now proves commutativity of the hexagon 
corresponding to the original cube. 

By induction~\ref{lem:cocycle2} follows. 

\subsection{Blocks of staircase}
\label{subsec:strcs}
We give a proof of Lemma \ref{lem:cocycle3} below.

\ssss
\label{sssec:strcs1}
Let us recall \ref{lem:cocycle3} for convenience.
\emph{
Let $X \stackrel{\sigma_1\;}{\Lra} Y \stackrel{\sigma_2\;}{\Lra} X$ 
be sequences such that $|\sigma_1\star\sigma_2| = 1_X$.
Then with $G \set (|\sigma_1|, \sfP)$,
the following diagram of isomorphisms commutes for any available choice 
of $\sS$.
\[
\xymatrix{
\sigma_1^{\bxt}\sigma_2^{\bxt} 
 \ar[rr]^{\Psi_{G,\sigma_1}^{\sS}(\sigma_2^{\bxt})} 
 \ar[rd]_{\Phi_{\sigma_1\star\sigma_2}} & & 
 G^{\bxt}\sigma_2^{\bxt} \ar[ld]^{\Phi_{G\star\sigma_2}} \\
& \xyoneD{X}
}
\]
}

Consider the following cartesian diagrams Sc1-Sc3 explained below.
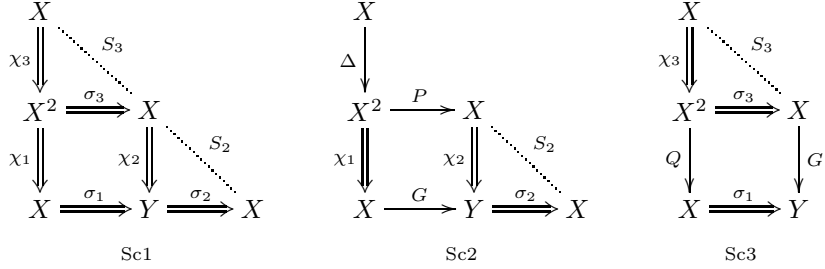
\[
\substack{
\xymatrix{
X \ar@{=>}[d]_{\chi_3} \ar@{.}[rd]^{S_3} \\
X^2 \ar@{=>}[d]_{\chi_1} \ar@{=>}[r]^{\sigma_3} & X 
 \ar@{=>}[d]_{\chi_2} \ar@{.}[rd]^{S_2} \\
X \ar@{=>}[r]^{\sigma_1} & Y \ar@{=>}[r]^{\sigma_2} & X
} 
\vspace{.1in} \\  \text{Sc1}} \vspace{.1cm}
\qquad
\substack{
\xymatrix{
X \ar[d]_{\Delta} \\
X^2 \ar@{=>}[d]_{\chi_1} \ar[r]^{P} & X \ar@{=>}[d]_{\chi_2} 
 \ar@{.}[rd]^{S_2} \\
X \ar[r]^{G} & Y \ar@{=>}[r]^{\sigma_2} & X
}
\vspace{.1in} \\  \text{Sc2}} \vspace{.1cm}
\qquad
\substack{
\xymatrix{
X \ar@{=>}[d]_{\chi_3} \ar@{.}[rd]^{S_3} \\
X^2 \ar[d]_{Q} \ar@{=>}[r]^{\sigma_3} & X \ar[d]^{G}  \\
X \ar@{=>}[r]^{\sigma_1} & Y
}
\vspace{.1in} \\  \text{Sc3}} \vspace{.1cm}
\]
We begin by choosing Sc1 to be a staircase based
on $\sigma_1 \star \sigma_2$.
Thus $S_2$ (resp.~$S_3$) is the sequence of steps induced by $\sigma_2$
(resp.~$\sigma_3$) and $S_3 \star S_2$, is the sequence of 
steps induced by $\sigma_1 \star \sigma_2$.
The other two diagrams,~Sc2 and~Sc3, are obtained by collapsing certain
portions of~Sc1: In~Sc2 we set $\Delta \set (|\chi_3|, \sfP)$ 
and $P \set (|\sigma_3|, \sfP)$ while in~Sc3 we set 
$Q \set (|\chi_1|, \sfP)$; these labels make sense. Indeed,
by assumption $|\sigma_1| \in \sfP$, so that by base change,
$|\sigma_3| \in \sfP$ and by definition of a staircase, $|\chi_1| \in \sfP$. 
Also note that by definition of a staircase,
$|\chi_2| = |\sigma_1|$ so that $G = (|\chi_2|, \sfP)$.

It follows that the squares occurring in~Sc2 and~Sc3 are cartesian and the 
triangles in them represent staircases based on $\sigma_2$ and $\sigma_3$
respectively.

In order to prove \ref{lem:cocycle3}, let us fix a 
cartesian diagram $\sS$ for the purpose of 
defining $\Psi_{G,\sigma_1}^{\sS}$. Given $\sS$, could we have 
chosen Sc1 above in
such a manner that the resulting square in Sc3 is precisely $\sS$ ?
The answer is yes; see \ref{lem:strcs1a} below. We shall therefore 
assume that Sc1 has been chosen so that 
Sc1 and $\sS$ are compatible in the said manner.

For $i = 2,3$, let $\alpha_i \colon \chi_i^{\bxt}\sigma_i^{\bxt} 
\iso S_i^{\bxt} \iso \oneD{X}$ be
the obvious natural isomorphisms (see last two maps in the 
composition defining~\eqref{eq:psi3}). Consider the following diagram
of isomorphisms where reduced notation has been employed for
naming the maps (\S\ref{subsec:conv}\eqref{conv2}).
\[
\begin{CD}
@. \framebox{$\sigma_1^{\bxt}\sigma_2^{\bxt}$} 
 @= \sigma_1^{\bxt}\sigma_2^{\bxt} \\
@. @V{\bphi_{\chi_3\star\chi_1}}V{\hspace{4em} \ddag_1}V 
 @VV{\bphi_{\chi_3\star Q}}V \\
@. \chi_3^{\bxt}\chi_1^{\bxt}\sigma_1^{\bxt}\sigma_2^{\bxt} 
 @<{\Psi_{\chi_1, Q}}<< \chi_3^{\bxt}Q^{\bxt}\sigma_1^{\bxt}\sigma_2^{\bxt} \\
@. @V{\bbeta_{\sigma_1,\chi_2}}V{\hspace{4em} \ddag_2}V 
 @VV{\bbeta_{\sigma_1,G}}V \\
@. \chi_3^{\bxt}\sigma_3^{\bxt}\chi_2^{\bxt}\sigma_2^{\bxt} 
 @<{\Psi_{\chi_2, G}}<< \chi_3^{\bxt}\sigma_3^{\bxt}G^{\bxt}\sigma_2^{\bxt} \\
@. @V{\alpha_3}VV @VV{\alpha_3}V \\
S_3^{\bxt}  @<{\alpha_2}<< 
 S_3^{\bxt}\chi_2^{\bxt}\sigma_2^{\bxt} @<{\Psi_{\chi_2, G}}<< 
 S_3^{\bxt}G^{\bxt}\sigma_2^{\bxt} \\
@V{\bphi_{S_3}}VV  @V{\bphi_{S_3}}VV @VV{\bphi_{S_3}}V \\
\framebox{$\boldsymbol{1}_{\cD_X}$} 
 @<{\alpha_2}<< \chi_2^{\bxt}\sigma_2^{\bxt} @<{\Psi_{\chi_2, G}}<< 
 \framebox{$G^{\bxt}\sigma_2^{\bxt}$} \\
@A{\bphi_{\Delta \star P}}AA 
 @A{\bphi_{\Delta \star P}}A{\hspace{4em} \ddag_3}A 
 @VV{\bphi_{\Delta \star \chi_1}^{-1}}V \\
\Delta^{\bxt} P^{\bxt} 
 @<{\alpha_2}<< \Delta^{\bxt} P^{\bxt}\chi_2^{\bxt}\sigma_2^{\bxt} 
 @<<{\bbeta_{G,\chi_2}}< \Delta^{\bxt}\chi_1^{\bxt}G^{\bxt}\sigma_2^{\bxt}
\end{CD}
\]
The unnamed rectangles commute by functoriality.
In $\ddag_1$, upon canceling off the common factor 
$\sigma_1^{\bxt}\sigma_2^{\bxt}$ from each vertex, we are left with  
objects in which every sequence is composed only of maps 
labeled~$\sfP$. Therefore $\ddag_1$ commutes by 
pseudofunctoriality. Similarly, for $\ddag_3$, we first
cancel off $\sigma_2^{\bxt}$ from each vertex and then
pseudofunctoriality of $(-)^{\times}$ applies. 
Finally, for $\ddag_2$, we first cancel 
$\chi_3^{\bxt}$ on the left and $\sigma_2^{\bxt}$ on
the right and then use transitivity of $\bbeta$. Thus the 
preceding diagram commutes. Now Lemma~\ref{lem:cocycle3} follows 
by looking at the outer border comprising the framed objects 
and the composite maps between them.

We used the following lemma in the above proof.

\begin{alem}
\label{lem:strcs1a}
Let $\sS$ be a cartesian square of labeled sequences 
as follows where~$G$ and~$Q$ have length one.
\[
\xymatrix{
X' \ar[d]_{Q} \ar@{=>}[rr]^{\sigma_3} & &  Y' \ar[d]^G \\
X \ar@{=>}[rr]^{\sigma_1} & & Y
}
\]
Let $\sigma_2 \colon Y' \Lra Y$ 
be a sequence such that $|\sigma_2| = |G|$. 
Then there exists a cartesian square $\sT$ as follows 
\[
\xymatrix{
X' \ar@{=>}[d]_{\sigma_4} \ar@{=>}[rr]^{\sigma_3} 
 & &  Y' \ar@{=>}[d]_{\sigma_2} \\
X \ar@{=>}[rr]^{\sigma_1} & & Y
}
\]
such that $|\sigma_4| = |Q|$. 
\end{alem}
\begin{proof}
If $\sigma_2$ has length 
one then the result holds trivially with \mbox{$\sT = \sS$}. Otherwise, we
may decompose $\sigma_2$ as $Y' \xto{H} Y_1 \stackrel{\rho\;}{\Lra} Y$.
We then proceed by induction on the length of $\sigma_1$.
Suppose $\sigma_1$ has length more than two. Then we may 
decompose~$\sS$ as shown on the left below. 
\[
\xymatrix{
X' \ar[d] \ar@{=>}[r] & W' \ar[d] \ar[r]^{E'} & 
 Z' \ar[d]^K \ar[r]^{F'} & Y' \ar[d]^{G} \\
X \ar@{=>}[r] & W \ar[r]^E & Z \ar[r]^F & Y
}
\qquad \qquad
\xymatrix{
\widetilde{X} \ar[d] \ar@{=>}[r] & \widetilde{W} 
 \ar[d] \ar[r]^{\widetilde{E}} & \widetilde{Z} 
 \ar[d]^{\widetilde{K}} \ar[r]^{\widetilde{F}} & Y' \ar[d]^H \\
X_1 \ar@{=>}[d] \ar@{=>}[r] & 
 W_1 \ar@{=>}[d] \ar[r] & 
 Z_1 \ar@{=>}[d]^{\kappa} \ar[r] & Y_1 \ar@{=>}[d]^{\rho} \\
X \ar@{=>}[r] & W \ar[r]^E & Z \ar[r]^F & Y
}
\]
Let us choose any cartesian diagram $\sT_0$ as shown on the right above.
Then $Z'$ and~$\widetilde{Z}$ are universal
for the same fibered diagram, and so there are canonical 
inverse isomorphisms
$i \colon \widetilde{Z} \iso Z'$ and $j \colon Z' \to \widetilde{Z}$
such that $|\widetilde{F}|j = |F'|$ and
$|\kappa||\widetilde{K}|j = |K|$. Let $\sT_1$ be the diagram obtained
by making the following replacements in $\sT_0$
\[
\widetilde{Z} \leadsto Z, \qquad \widetilde{F} \leadsto F',
\qquad \widetilde{K} \leadsto K, \qquad \widetilde{E} \leadsto 
(i|\widetilde{E}|, \lambda),
\]
where $\lambda$ is the label of $E$. 
Then $\sT_1$ is also a cartesian diagram and it agrees with $\sS$ 
in the rightmost column of squares. We proceed in a similar manner 
replacing one by one the vertices and edges in the top row of 
$\sT_i$ with that of $\sS$. The case when $\sigma_1$ has length one 
or two is handled using similar (and simpler) arguments.
\end{proof}

\subsection{Proof of Lemma \ref{lem:cocycle4}}
\label{subsec:PL}
This is the third and last of the three lemmas used in the proof
of Theorem~\ref{thm:cocycle1}. Its proof is somewhat longer. 
The compatibility of base-change isomorphisms with the
fundamental isomorphism plays a role here.

\ssss
\label{sssec:PL1}
Let us begin by recalling what has to be proven.
The cartesian diagram under consideration is  
\begin{equation}
\label{eq:PL2}
\xymatrix{
X \ar@{=>}[r]^{\sigma_1} & A \ar@{=>}[d]_{\sigma_3} \ar@{=>}[r]^{\sigma_2} 
 & B \ar@{=>}[d]^{\sigma_4}  \\
& C \ar@{=>}[r]_{\sigma_5} & D \ar@{=>}[r]_{\sigma_6} & X 
}
\end{equation}
where the induced composite map $X \to X$ is the identity.
Then the assertion is that
the following diagram of isomorphisms commutes where the~$\Phi$'s 
have the obvious subscripts.
\begin{equation}
\label{eq:PL2a}
\xymatrix{
\sigma_1^{\bxt} \sigma_3^{\bxt} \sigma_5^{\bxt} \sigma_6^{\bxt}
 \ar[rr]^{\text{via } \;\bbeta_{\sigma_5, \sigma_4}} \ar[rd]_{\Phi} 
 & & \sigma_1^{\bxt} \sigma_2^{\bxt} \sigma_4^{\bxt} \sigma_6^{\bxt} 
 \ar[ld]^{\Phi} \\
& \mathbf{1}_{\cD_X}
}
\end{equation}
A quick remark on conventions: in the above context, we shall also
allow $\sigma_6$ to be the empty sequence, 
and in such a case we declare 
$D = X$ and $\sigma_6^{\bxt} = \oneD{X}$. 

Our immediate goal is to show that the following restrictions
may be imposed on the hypothesis of Lemma \ref{lem:cocycle4} 
(or \eqref{eq:PL2} above) without loss of generality.
\begin{enumerate}
\item The sequences $\sigma_2$, $\sigma_3$, $\sigma_4$ and~$\sigma_5$, 
all have length one.
\item The tail sequence $\sigma_6$ is empty.
\item The head sequence $\sigma_1$ is a single map with label $\sfP$.
\end{enumerate}

\ssss
\label{sssec:PL3}
\emph{Reduction to} (i): We prove this by induction on the sum of
the lengths of~$\sigma_4$ and~$\sigma_5$. 

Suppose $\sigma_5 = \mu_1\star\mu_2$ and let 
\[
\xymatrix{
X \ar@{=>}[r]^{\sigma_1} & A \ar@{=>}[d]_{\sigma_3} \ar@{=>}[r]^{\nu_1} 
 & B_1 \ar@{=>}[d]_{\theta} \ar@{=>}[r]^{\nu_2} & B \ar@{=>}[d]^{\sigma_4} \\
& C \ar@{=>}[r]_{\mu_1} & D_1 \ar@{=>}[r]_{\mu_2} & D 
 \ar@{=>}[r]_{\sigma_6} & X 
}
\]
be the corresponding induced cartesian diagram. Consider the following
diagram of isomorphisms
\[
\xymatrix{
\sigma_1^{\bxt}\sigma_3^{\bxt}\sigma_5^{\bxt}\sigma_6^{\bxt} 
 \ar[rr]^{\text{via } \;\bbeta_{\sigma_5, \sigma_4}} \ar@{=}[dd] 
 \ar[dr]_{\Phi}^{\hspace{8.5em} \ddag} 
 & \text{\phantom{AAAAAAAAAAAAAAA}} & 
 \sigma_1^{\bxt}\sigma_2^{\bxt}\sigma_4^{\bxt}\sigma_6^{\bxt} 
 \ar[dl]^{\Phi} \ar@{=}[dd] \\
& \mathbf{1}_{\cD_X} \\
\sigma_1^{\bxt}\sigma_3^{\bxt}\mu_1^{\bxt}\mu_2^{\bxt}\sigma_6^{\bxt} 
 \ar[r]^{\text{via }}_{\bbeta_{\mu_1,\theta}} \ar[ur]^{\Phi} &
 \sigma_1^{\bxt}\nu_1^{\bxt}\theta^{\bxt}\mu_2^{\bxt}\sigma_6^{\bxt} 
 \ar[r]^{\text{via }}_{\bbeta_{\mu_2,\sigma_4}} 
 \ar[u]_{\hspace{1.5em}\dag_2}^{\dag_1\hspace{1.5em}} 
 & \sigma_1^{\bxt}\nu_1^{\bxt}\nu_2^{\bxt}\sigma_4^{\bxt}\sigma_6^{\bxt} 
 \ar[ul]_{\Phi}
}
\]
where the $\Phi$'s have the obvious subscripts.
The outer border of the above diagram commutes by definition 
of~$\bbeta_{\sigma_5, \sigma_4}$. Therefore, if by induction hypothesis,
we assume that $\dag_1, \dag_2$ commute, then so does \ddag. 

A similar argument works if we begin with a decomposition of~$\sigma_4$ into 
two sequences. Thus by induction, we obtain reduction to (i).

\ssss
\label{sssec:PL4}
\emph{Reduction to} (ii):
We may assume that the (i) has been applied.
We proceed in two stages.

Set $g \set |\sigma_1\star\sigma_3\star\sigma_5|$ 
and $G \set (g, \sfP)$. Consider the following diagram
\begin{equation}
\label{eq:PL4a}
\xymatrix{
\sigma_1^{\bxt}\sigma_3^{\bxt}\sigma_5^{\bxt} \ar[rd]_{\Psi_3} 
 \ar[rr]^{\text{via }\;\bbeta_{\sigma_5, \sigma_4}}  
 & & \sigma_1^{\bxt}\sigma_2^{\bxt}\sigma_4^{\bxt} \ar[ld]^{\Psi_2} \\
& G^{\bxt}
}
\end{equation}
where $\Psi_3 \set \Psi_{G, \;\sigma_1\star\sigma_3\star\sigma_5}^{\sS}$
and $\Psi_2 \set \Psi_{G, \;\sigma_1\star\sigma_2\star\sigma_4}^{\sT}$
for some suitable choice of~$\sS$ and~$\sT$. 
We claim that to prove that \eqref{eq:PL2a} commutes it suffices to 
find a choice for $\sS$ and $\sT$ such that 
\eqref{eq:PL4a} commutes. The claim follows by looking at the 
following diagram
\[
\xymatrix{
\qquad\sigma_1^{\bxt}\sigma_3^{\bxt}\sigma_5^{\bxt}\sigma_6^{\bxt} 
 \ar[rdd]_{\Phi} \ar[rd]^{\Psi_3(\sigma_6^{\bxt})} 
 \ar[rr]^{\text{via }\;\bbeta_{\sigma_5, \sigma_4}} 
 & \text{\phantom{AAAAAAAAAAA}} 
 & \sigma_1^{\bxt}\sigma_2^{\bxt}\sigma_4^{\bxt}\sigma_6^{\bxt}\qquad 
 \ar[ld]_{\Psi_2(\sigma_6^{\bxt})} \ar[ldd]^{\Phi} \\
& G^{\bxt}\sigma_6^{\bxt} \ar[d]^{\!\Phi} \\
& F^{\bxt}
}
\]
where the outer border is \eqref{eq:PL2a}, 
the lower two sub-triangles on either side commute by \ref{lem:cocycle3},
and commutativity of the subtriangle on the top reduces to
that of~\eqref{eq:PL4a}. 

Now consider a 
cartesian diagram as follows, which is obtained by taking
a product of \eqref{eq:PL2}
(with a truncated tail) with $G$.
\[
\xymatrix{
X \ar[r]^{\Delta} & X^2 \ar@{=>}[r]^{\rho_1} 
 \ar[ddd]^{\mu_1} & A' \ar[ddd]^{\mu_2} \ar[dr]^{\rho_2} 
 \ar[rrr]^{\rho_3} & & & C' \ar[ddd]^{\mu_5} \ar[dl]_{\rho_5} \\ 
& & & B' \ar[d]^{\mu_3} \ar[r]^{\rho_4}  
 & X \ar[d]^{G} & \\ 
& & & B  \ar[r]_{\sigma_4} & D  & \\ 
& X \ar@{=>}[r]_{\sigma_1} & A \ar[ur]_{\sigma_2} \ar[rrr]_{\sigma_3}
 & & & C  \ar[ul]^{\sigma_5}
}
\]
Let $\Psi_3$ and $\Psi_2$ (and in particular, $\sS$ and $\sT$) 
of~\eqref{eq:PL4a} be 
defined via the cartesian squares in the preceding diagram. 
Now consider the following diagram of isomorphisms.
\[
\xymatrix{
\sigma_1^{\bxt}\sigma_3^{\bxt}\sigma_5^{\bxt} 
 \ar[ddd]_{\Phi^{-1}_{\Delta\star\mu_1}}^{\hspace{3em} \lozenge_3} 
 \ar[rd]^{\Psi_3} \ar[rrr]^{\bbeta_{\sigma_5, \sigma_4}}  
 & \hspace{7em} & \hspace{7em} & 
 \sigma_1^{\bxt}\sigma_2^{\bxt}\sigma_4^{\bxt} \ar[ld]_{\Psi_2} 
 \ar[ddd]^{\Phi^{-1}_{\Delta\star\mu_1}}_{\lozenge_4 \hspace{3em}} \\ 
&  G^{\bxt} 
\ar@{=}[r]^{\smash{\raisebox{3ex}{$\scriptstyle \lozenge_1$}}}_{\smash{
 \raisebox{-3ex}{$\scriptstyle \lozenge_2$}}} 
 & G^{\bxt} \\ 
& \Delta^{\!\bxt}\rho_1^{\bxt}\rho_3^{\bxt}\rho_5^{\bxt}G^{\bxt} 
 \ar[u]^{\Phi_1} \ar[r]^{\bbeta_{\rho_5, \rho_4}} &  
 \Delta^{\!\bxt}\rho_1^{\bxt}\rho_2^{\bxt}\rho_4^{\bxt}G^{\bxt} 
 \ar[u]_{\Phi_2} \\ 
\Delta^{\!\bxt}\mu_1^{\bxt}\sigma_1^{\bxt}\sigma_3^{\bxt}\sigma_5^{\bxt} 
 \ar[ru]^{\bbeta_1} 
 \ar[rrr]_{\bbeta_{\sigma_5, \sigma_4}}^{\smash{
 \raisebox{3ex}{$\scriptstyle \lozenge_5$}}} & & & 
 \Delta^{\!\bxt}\mu_1^{\bxt}\sigma_1^{\bxt}\sigma_2^{\bxt}\sigma_4^{\bxt} 
 \ar[lu]_{\bbeta_2}
}
\]
The morphisms are indicated in reduced notation using the following.
\[
\Phi_1 = \Phi_{\Delta\star\rho_1\star\rho_3\star\rho_5}
\quad
\Phi_2 = \Phi_{\Delta\star\rho_1\star\rho_2\star\rho_4}
\quad
\bbeta_1 = \bbeta_{\sigma_1\star\sigma_3\star\sigma_5, \; G}
\quad
\bbeta_2 = \bbeta_{\sigma_1\star\sigma_2\star\sigma_4, \; G}
\]
The outer border of the preceding diagram commutes for functorial 
reasons. From the definition of~$\Psi_3$ and~$\Psi_2$ we see that
$\lozenge_3$ and $\lozenge_4$ commute (cf.~\ref{sssec:psi5z}). 
In $\lozenge_5$, if we first cancel off $\Delta^{\!\bxt}$ on the right, 
then the transpose of the resulting diagram may be expanded 
as follows.
\[
\begin{CD}
\mu_1^{\bxt}\sigma_1^{\bxt}\sigma_3^{\bxt}\sigma_5^{\bxt} 
 @>{\bbeta_{\sigma_1, \mu_2}}>> 
 \rho_1^{\bxt}\mu_2^{\bxt}\sigma_3^{\bxt}\sigma_5^{\bxt} 
 @>{\bbeta_{\sigma_3, \mu_5}}>>
 \rho_1^{\bxt}\rho_3^{\bxt}\mu_5^{\bxt}\sigma_5^{\bxt} 
 @>{\bbeta_{\sigma_5, G}}>> \rho_1^{\bxt}\rho_3^{\bxt}\rho_5^{\bxt}G^{\bxt} \\
@V{\bbeta_{\sigma_5, \sigma_4}}VV @VV{\bbeta_{\sigma_5, \sigma_4}}V 
 @. @V{\bbeta_{\rho_5, \rho_4}}VV \\
\mu_1^{\bxt}\sigma_1^{\bxt}\sigma_2^{\bxt}\sigma_4^{\bxt} 
 @>{\bbeta_{\sigma_1, \mu_2}}>> 
 \rho_1^{\bxt}\mu_2^{\bxt}\sigma_2^{\bxt}\sigma_4^{\bxt} 
 @>{\bbeta_{\sigma_2, \mu_3}}>>
 \rho_1^{\bxt}\rho_2^{\bxt}\mu_3^{\bxt}\sigma_4^{\bxt} 
 @>{\bbeta_{\sigma_4, G}}>> \rho_1^{\bxt}\rho_2^{\bxt}\rho_4^{\bxt}G^{\bxt} \\
\end{CD}
\]
Once again, reduced notation has been applied here. The rectangle on the
left commutes for functorial reasons while the one on the right
commutes by the cube lemma. Thus $\lozenge_5$ commutes.

It follows that $\lozenge_1$ commutes if and only if 
$\lozenge_2$ commutes. Now consider the following two diagrams. 
\begin{equation}
\label{eq:PL5a}
\xymatrix{
X \ar@{=>}[r]^{\Delta\star\rho_1} & A' \ar[d]_{\rho_3} 
 \ar[r]^{\rho_2} & B' \ar[d]^{\rho_4}  \\
\qquad & C' \ar[r]_{\rho_5} & D'\makebox[0pt]{$\quad\quad\> = X$} 
}
\qquad \quad
\xymatrix{
\Delta^{\!\bxt}\rho_1^{\bxt}\rho_3^{\bxt}\rho_5^{\bxt} \ar[rd]_{\Phi} 
 \ar[rr]^{\text{via } \bbeta_{\rho_5, \rho_4}} & & 
 \Delta^{\!\bxt}\rho_1^{\bxt}\rho_2^{\bxt}\rho_4^{\bxt} \ar[ld]^{\Phi} \\
& \boldsymbol{1}_{\cD_X}
}
\end{equation}
The one on the left is an instance of the cartesian diagram in \eqref{eq:PL2} 
with tail sequence empty. The one on the 
right corresponds to~\eqref{eq:PL2a}. Now $\lozenge_2$, 
upon the canceling of~$G^{\bxt}$ from each vertex, results in the 
diagram shown on the right in~\eqref{eq:PL5a}. Therefore, proving that
$\lozenge_1$ commutes, reduces to proving that the diagram on the 
right in~\eqref{eq:PL5a} commutes. 

In summary, if the diagram on the 
right in~\eqref{eq:PL5a} commutes then so does~$\lozenge_1$ which
is nothing but \eqref{eq:PL4a}, and hence so does \eqref{eq:PL2a}. 
We have therefore achieved reduction to (ii).

\ssss
\label{sssec:PL6}
\emph{Reduction to} (iii):
We may now assume that the restrictions of~(i) and~(ii) apply.

Let us set $H \set (|\sigma_1|, \sfP)$.
Now consider the following diagram of isomorphisms 
where reduced notation has been used
for the morphisms, the $\Phi$'s have the obvious subscripts
and $\sS$ is a suitable fixed diagram used in defining 
$\Psi_{H, \sigma_1}^{\sS}$.
\[
\xymatrix{
\sigma_1^{\bxt}\sigma_3^{\bxt}\sigma_5^{\bxt} 
 \ar[dd]_{\Psi_{H, \sigma_1}^{\sS}} 
 \ar[rd]^{\Phi} \ar[rr]^{\bbeta_{\sigma_5, \sigma_4}}
 & \hspace{4em} & 
 \sigma_1^{\bxt}\sigma_2^{\bxt}\sigma_4^{\bxt} \ar[ld]_{\Phi} 
 \ar[dd]^{\Psi_{H, \sigma_1}^{\sS}}  \\
& \boldsymbol{1}_{\cD_X} \\
H^{\bxt}\sigma_3^{\bxt}\sigma_5^{\bxt} \ar[ru]_{\Phi} 
 \ar[rr]_{\bbeta_{\sigma_5, \sigma_4}} & & 
 H^{\bxt}\sigma_2^{\bxt}\sigma_4^{\bxt} \ar[lu]^{\Phi}
}
\]
The outer border commutes for functorial reasons. The two triangles on 
the left and right sides commute by~\ref{lem:cocycle3}. Therefore for the 
top triangle to commute, it suffices that the bottom one commutes. 
This is precisely reduction to (iii).

\ssss
\label{sssec:PL7}
Once all the reductions (i)-(iii)
have been carried out, we may
replace the diagrams of \eqref{eq:PL2} and~\eqref{eq:PL2a} 
by the following ones where 
each $\sigma_i$ has length one, $\sigma_1$ has label~$\sfP$ and 
$|\sigma_1\star\sigma_3\star\sigma_5| = 1_X$.
\begin{equation} 
\label{eq:PL8}
\xymatrix{
X \ar[r]^{\sigma_1} & A \ar[d]_{\sigma_3} 
 \ar[r]^{\sigma_2} & B \ar[d]^{\sigma_4}  \\
 & C \ar[r]_{\sigma_5} & X
}
\qquad\qquad
\xymatrix{
\sigma_1^{\bxt}\sigma_3^{\bxt}\sigma_5^{\bxt} 
 \ar[rd]^{\bphi} \ar[rr]^{\bbeta_{\sigma_5, \sigma_4}}
 & & \sigma_1^{\bxt}\sigma_2^{\bxt}\sigma_4^{\bxt} \ar[ld]_{\bphi} \\
& \boldsymbol{1}_{\cD_X}
}
\end{equation}
Set $t_4 \set |\sigma_1\star\sigma_2|\colon X \to B$ and
$t_5 \set |\sigma_1\star\sigma_3|\colon X \to C$, so that 
$t_4,t_5$ are sections of $|\sigma_4|,|\sigma_5|$ respectively.
Therefore $t_4, t_5 \in \sfP$. For $i = 4,5$ set $\theta_i \set (t_i, \sfP)$.
Now taking a fibered product of the sequence
$X \xto{\theta_4} B \xto{\sigma_4} X$ with 
$X \xto{\theta_5} C \xto{\sigma_5} X$ results in a diagram as
follows for suitably determined~$\theta_2$ and~$\theta_3$.
\begin{equation} 
\label{eq:PL8a}
\begin{CD}
X @>{\theta_5}>> C @>{\sigma_5}>> X \\
@VV{\theta_4}V @VV{\theta_3}V @VV{\theta_4}V  \\
B @>{\theta_2}>> A @>{\sigma_2}>> B \\
@VV{\sigma_4}V @VV{\sigma_3}V @VV{\sigma_4}V  \\
X @>{\theta_5}>> C @>{\sigma_5}>> X 
\end{CD}
\end{equation} 
Set $t \set |\theta_4\star\theta_2| = |\theta_5\star\theta_3| \colon X \to A$.
We claim that $t = |\sigma_1|$. Indeed, note that for $i = 2,3$, 
we have $|\sigma_i|t = |\sigma_i||\sigma_1|$ because
\[
|\sigma_2|t = |\sigma_2||\theta_2||\theta_4| = |\theta_4| \set 
|\sigma_2||\sigma_1|,
\qquad
|\sigma_3|t = |\sigma_3||\theta_3||\theta_5| = |\theta_5| \set 
|\sigma_3||\sigma_1|. 
\]
Therefore, by the universal property of fibered products for
the cartesian square in \eqref{eq:PL8},
we conclude that $t = |\sigma_1|$. In particular, since 
$\sigma_1$ and $\theta_i$ for $i = 2,3,4,5$ have label $\sfP$,
by pseudofunctoriality there result natural isomorphisms 
\begin{equation}
\label{eq:PL9} 
\sigma_1^{\bxt} \iso \theta_4^{\bxt}\theta_2^{\bxt} \iso 
\theta_5^{\bxt}\theta_3^{\bxt}.
\end{equation}

Now consider the following diagram of isomorphisms. 
The superscript ${}^{\bxt}$ has been omitted from each term and 
reduced notation has been used for the morphisms.
The morphisms denoted by~$\pi_i$ are obtained using~\eqref{eq:PL9}.
The unlabeled ones are obtained using the 
suitably determined fundamental isomorphism. 
\begin{equation}
\label{eq:PL10}
\xymatrix{
\sigma_1\sigma_3\sigma_5 \ar[d]_{\pi_1}
 \ar[rrrr]^{\bbeta_{\sigma_5, \sigma_4}}  
 & & & & \sigma_1\sigma_2\sigma_4  \ar[d]^{\pi_2} \\ 
\theta_4\theta_2\sigma_3\sigma_5 \ar[dddd]_{\pi_3}^{\hspace{3em} \ddag_2} 
 \ar[rd]^{\bbeta_{\sigma_3,\theta_5}} 
 \ar[rrrr]^{\bbeta_{\sigma_5, \sigma_4}}_{\smash{
 \raisebox{-3ex}{$\scriptstyle \ddag_1$}}} 
 & & & & \theta_4\theta_2\sigma_2\sigma_4 
 \ar[dddd]^{\pi_4}_{\ddag_3 \hspace{3em}} \ar[ld] \\ 
& \theta_4\sigma_4\theta_5\sigma_5 \ar[dd] \ar[rr]
 & & \theta_4\sigma_4 \ar[dd] \ar[ld] \\ 
& & \boldsymbol{1}_{\cD_X} \ar[ld] \\ 
& \theta_5\sigma_5 \ar[ld] \ar[rr] & & 
 \theta_5\sigma_5\theta_4\sigma_4 \ar[rd]^{\bbeta_{\theta_4,\sigma_2}} \\ 
\theta_5\theta_3\sigma_3\sigma_5 
 \ar[rrrr]_{\bbeta_{\sigma_5, \sigma_4}}^{\smash{
 \raisebox{3ex}{$\scriptstyle \ddag_4$}}} 
 & & & & \theta_5\theta_3\sigma_2\sigma_4 
}
\end{equation}
Let us verify that this diagram commutes. The unmarked subdiagrams
commute for functorial reasons and so only $\ddag_{\>i}$ need to 
be considered. For these we 
use the compatibility of the fundamental isomorphism
with base change~(\ref{lem:cafpols4}(ii)). 
For instance, in~$\ddag_1$, we first cancel the common 
term~$\theta_4^{\bxt}$ and then use~\ref{lem:cafpols4}(ii)
for the following subdiagram of~\eqref{eq:PL8a}.
\[
\begin{CD}
B @>{\theta_2}>> A @>{\sigma_2}>> B \\
@VV{\sigma_4}V @VV{\sigma_3}V @VV{\sigma_4}V  \\
X @>{\theta_5}>> C @>{\sigma_5}>> X 
\end{CD}
\]
A similar argument works for $\ddag_2$, $\ddag_3$ and $\ddag_4$. 

Since \eqref{eq:PL10} commutes, to finish the proof of commutativity 
of the diagram on 
right in \eqref{eq:PL8}, it suffices to show that 
the following composites obtained from \eqref{eq:PL10} give 
$\Phi_{\sigma_1\star\sigma_3\star\sigma_5}$ and 
$\Phi_{\sigma_1\star\sigma_2\star\sigma_4}$ respectively.
\begin{gather}
\sigma_1^{\bxt}\sigma_3^{\bxt}\sigma_5^{\bxt} \xto{\;\pi_1\;}
\theta_4^{\bxt}\theta_2^{\bxt}\sigma_3^{\bxt}\sigma_5^{\bxt} 
 \xrightarrow[{\bbeta_{\sigma_3,\theta_5}}]{\text{via}}
\theta_4^{\bxt}\sigma_4^{\bxt}\theta_5^{\bxt}\sigma_5^{\bxt} \xto{\;\Phi\;}
\theta_5^{\bxt}\sigma_5^{\bxt} \xto{\;\Phi\;}\oneD{X} \\ 
\sigma_1^{\bxt}\sigma_2^{\bxt}\sigma_4^{\bxt} \xto{\pi_4\pi_2}
\theta_5^{\bxt}\theta_3^{\bxt}\sigma_2^{\bxt}\sigma_4^{\bxt} 
 \xrightarrow[{\bbeta_{\sigma_2,\theta_4}}]{\text{via}} 
\theta_5^{\bxt}\sigma_5^{\bxt}\theta_4^{\bxt}\sigma_4^{\bxt} \xto{\;\Phi\;}
\theta_4^{\bxt}\sigma_4^{\bxt} \xto{\;\Phi\;} \oneD{X} \label{eq:PL11}
\end{gather}
We only consider the second case, viz., that of showing
$\Phi_{\sigma_1\star\sigma_2\star\sigma_4} = \eqref{eq:PL11}$. 
The first one is proved in a similar manner.

Consider the following cartesian diagram where $\kappa_i$ for $i = 2,3,4,5$
is constructed via the cartesian condition and~$\kappa_1$ is the diagonal 
map with label~$\sfP$.
\[
\begin{CD}
X \\
@VV{\kappa_1}V \\
\bullet @>{\kappa_2}>> X \\
@VV{\kappa_3}V  @VV{\theta_5}V \\
\bullet @>{\kappa_4}>> C @>{\sigma_5}>> X \\
@VV{\kappa_5}V  @VV{\theta_3}V  @VV{\theta_4}V \\
X @>{\sigma_1}>> A @>{\sigma_2}>> B @>{\sigma_4}>> X 
\end{CD}
\]
By definition, $|\theta_4| \set |\sigma_1 \star \sigma_2|$ and 
as verified above, $|\theta_5\star\theta_3| = |\sigma_1|$ and 
$|\theta_5\star\sigma_5| = 1_X$. Therefore, the preceding
diagram is a staircase based 
on~$\sigma_1\star\sigma_2\star\sigma_4$.  

Now consider the following diagram of isomorphisms.
\[
\begin{CD}
\sigma_1^{\bxt}\sigma_2^{\bxt}\sigma_4^{\bxt} @>>> 
 \theta_5^{\bxt}\theta_3^{\bxt}\sigma_2^{\bxt}\sigma_4^{\bxt} @>>>
 \theta_5^{\bxt}\sigma_5^{\bxt}\theta_4^{\bxt}\sigma_4^{\bxt} @>>>
 \oneD{X} \\
@VVV @VVV @VVV @| \\
 \kappa_1^{\bxt}\kappa_3^{\bxt}
 \kappa_5^{\bxt}\sigma_1^{\bxt}\sigma_2^{\bxt}\sigma_4^{\bxt} 
 @>>> \kappa_1^{\bxt}\kappa_2^{\bxt}
 \theta_5^{\bxt}\theta_3^{\bxt}\sigma_2^{\bxt}\sigma_4^{\bxt} 
 @>>> \kappa_1^{\bxt}\kappa_2^{\bxt}
 \theta_5^{\bxt}\sigma_5^{\bxt}\theta_4^{\bxt}\sigma_4^{\bxt} 
 @>>> \oneD{X}
\end{CD}
\]
Here the top row gives \eqref{eq:PL11} and the rest of the 
outer border spells out the 
definition of~$\Phi_{\sigma_1\star\sigma_2\star\sigma_4}$.
Each vertical map is induced in the obvious way by a 
suitably determined fundamental isomorphism.
The rectangle on the left, upon the cancellation of the
common term $\sigma_2^{\bxt}\sigma_4^{\bxt}$, consists only
of $\sfP$-labeled maps and hence commutes by pseudofunctoriality.
The other two rectangles commute for functorial reasons. Thus 
the preceding diagram commutes. In particular, 
$\Phi_{\sigma_1\star\sigma_2\star\sigma_4} = \eqref{eq:PL11}$. 

Thus we have shown that the diagram on the right in~\eqref{eq:PL8}
commutes. This completes the proof of Lemma~\ref{lem:cocycle4}.

\subsection{Canonicity of $\Psi_{-,-}$}
\label{subsec:conseq}
Now that the proof of Theorem \ref{thm:cocycle1} has been completed, 
we discuss some of the immediate consequences.

\begin{aprop}
\label{prop:conseq1}
The isomorphism $\Psi^{\sS}_{\sigma_1,\sigma_2}$ of \eqref{eq:psi5} 
is independent of the choice of $\sS$.
\end{aprop}
\begin{proof}
First we claim that for any sequence $\sigma$ and any available choice
of diagram $\sT$, the isomorphism $\Psi^{\sT}_{\sigma, \sigma}$ 
is the identity. Indeed, the cocycle condition says that
\[
\Psi^{\sT}_{\sigma, \sigma}\Psi^{\sT}_{\sigma, \sigma} = 
\Psi^{\sT}_{\sigma, \sigma}
\]
and since $\Psi^{\sT}_{\sigma, \sigma}$ is an isomorphism, it is necessarily
the identity.

Now suppose $\sS$ and $\sS'$ are two diagrams through which 
$\Psi^{-}_{\sigma_1,\sigma_2}$ is defined. Let~$\sT$ be any diagram 
through which $\Psi^{-}_{\sigma_2,\sigma_2}$ is defined. The cocycle
condition implies that the following holds
\[
\Psi^{\sS}_{\sigma_1,\sigma_2}\Psi^{\sT}_{\sigma_2,\sigma_2} = 
\Psi^{\sS'}_{\sigma_1,\sigma_2}
\] 
and since $\Psi^{\sT}_{\sigma_2,\sigma_2}$ is the identity, therefore
$\Psi^{\sS}_{\sigma_1,\sigma_2} = \Psi^{\sS'}_{\sigma_1,\sigma_2}$.
\end{proof}

\begin{adefi}
\label{def:conseq2}
For any two sequences $\sigma_1,\sigma_2$ such that 
$|\sigma_1| = |\sigma_2|$ we define $\Psi_{\sigma_1,\sigma_2}$ to be 
$\Psi_{\sigma_1,\sigma_2}^{\sS}$ of \eqref{eq:psi5} for any available $\sS$.
\end{adefi}

\begin{aprop}
\label{prop:conseq3}
The isomorphism $\Psi_{-,-}$ is reflexive and symmetric, i.e., the following
hold.
\begin{enumerate}
\item $\Psi_{\sigma,\sigma}$ is the identity for any $\sigma$.
\item $\Psi_{\sigma_1,\sigma_2} = \Psi_{\sigma_2,\sigma_1}^{-1}$
for any $\sigma_1,\sigma_2$ such that $|\sigma_1| = |\sigma_2|$.
\end{enumerate}
\end{aprop}
\begin{proof}
This is an immediate consequence of the cocycle rule.
See proof of~\ref{prop:conseq1}.
\end{proof}

\section{Proofs III (old isomorphisms and linearity)}
\label{sec:proofsIII}
There are two objectives in this section.
The first is to show that the isomorphism~$\Psi_{-,-}$
defined in the previous section recovers the isomorphisms
of the input data in~\S\ref{subsec:input} and more generally their
labeled counterparts. This is achieved in~\S\ref{subsec:recold}
and~\S\ref{subsec:recoldbc}. The other objective, 
discussed in~\S\ref{subsec:meiosis},
is of proving a linearity rule satisfied by~$\Psi_{-,-}$. 


\subsection{Recovering $\bC_{-,-}$ and $\Phi_{-}$ }
\label{subsec:recold}
Here record that the isomorphisms
$\bC_{-,-}$ and $\Phi_{-}$ are  
expressible in terms of $\Psi_{-,-}$.


Before we begin with the case of $\bC_{-,-}$, we need a lemma.
\begin{alem}
\label{lem:recC0}
Let $\sigma_1$ be the sequence $X \xto{G_1} Y \xto{G_2} Z \xto{G_3} X$ 
and $\sigma_2$ the sequence $X \xto{G_1} Y \xto{G_4} X$ such that
$G_2, G_3,G_4$ have label $\sfO$, $|G_3||G_2| = |G_4|$ and 
$|\sigma_i| = 1_X$. Then the following diagram commutes
\[
\xymatrix{
G_1^{\bxt}G_2^{\bxt}G_3^{\bxt} \ar[rd]_{\Phi_{\sigma_1}} 
 \ar[rr]^{\textup{via } \bC_{G_2,G_3}} 
 & &  G_1^{\bxt}G_4^{\bxt} \ar[ld]^{\Phi_{\sigma_2}} \\
& \mathbf{1}_{\cD_X}
}
\] 
\end{alem}
\begin{proof}
Consider the following staircases where $\sfr$ is obtained by vertically
composing~$\pfr$ and~$\qfr$.
\[
\begin{CD}
X \\
@V{\alpha_1}VV \\
X'' @>{\alpha_2}>> X \\
@V{\alpha_3}V{\hspace{2em} \pfr}V @VV{\alpha_4}V \\
X' @>{\alpha_5}>> Y' @>{\alpha_6}>> X \\
@V{\alpha_7}V{\hspace{2em} \qfr}V @VV{\alpha_8}V @VV{\alpha_9}V \\
X @>>{G_1}> Y @>>{G_2}> Z @>>{G_3}> X
\end{CD}
\qquad\qquad
\begin{CD}
X \\
@V{\alpha_1}VV \\
X'' @>{\alpha_2}>> X \\
@V{\alpha_{37}}V{\hspace{2em} \sfr}V @VV{\alpha_{48}}V \\
X @>>{G_1}> Y @>>{G_4}> X
\end{CD}
\]
Set $a_i \set |\alpha_i|$, $a_{ij} \set |\alpha_{ij}|$ and $g_i \set |G_i|$. 
By following through the definitions 
of~$\Phi_{\sigma_i}$ and upon using transitivity of base-change 
for $\pfr + \qfr = \sfr$, we see that the lemma reduces to checking that
the following diagram of isomorphisms commutes. 
\[
\begin{CD}
a_4^{\times}a_8^{\times}g_2^{\ssbox}g_3^{\ssbox} 
 @>>> a_4^{\times}a_6^{\ssbox}a_9^{\times}g_3^{\ssbox} \\
@VVV  @VVV \\
a_{48}^{\times}g_4^{\ssbox} @>>>  \oneD{X}
\end{CD}
\]
This diagram commutes by \S\ref{subsec:input},[D](ii).
\end{proof}

\begin{aprop}
\label{prop:recC}
\emph{(Recovering $\bC_{-,-}$).}
Let $X \xto{F_1} Y \xto{F_2} Z$ 
be maps having the same label~$\lambda$. 
Set $F_3 \set (|F_2||F_1|, \lambda)$. Then 
$\bC_{F_1,F_2} = \Psi_{F_3, F_1 \star F_2}$.
\end{aprop}
\begin{proof}
The definition of $\Psi_{F_3, F_1 \star F_2}$ leads
us to a cartesian diagram as follows.
\[
\begin{CD}
X @>{\Delta}>> X\times_ZX @>{\widetilde{F}_1}>> Y \times_Z X 
 @>{\widetilde{F}_2}>> X \\
@. @VV{F_3'}V @VVV @VV{F_3}V \\
@. X @>{F_1}>> Y @>{F_2}>> Z
\end{CD}
\] 
Now we argue via two cases. If $\lambda = \sfP$, then  
all the maps in the preceding diagram are $\sfP$-labeled.
Moreover any staircase on
$\Delta \star \widetilde{F}_1 \star \widetilde{F}_2$ also 
consists only of $\sfP$-labeled maps. Therefore the proposition 
follows from pseudofunctoriality of $(-)^{\times}$.

Suppose $\lambda = \sfO$. Set 
$\widetilde{F}_3 \set (|\widetilde{F}_2||\widetilde{F}_1|, \sfO)$. 
Consider the following diagram 
of isomorphisms where the top row gives $\Psi_{F_3,F_1 \star F_2}$
and the bottom row gives $\Psi_{F_3,F_3}$.
\[
\begin{CD}
F_1^{\bxt}F_2^{\bxt} @>>> \Delta^{\!\bxt}F_3'{}^{\bxt}F_1^{\bxt}F_2^{\bxt} 
 @>>> \Delta^{\!\bxt}\widetilde{F}_1{}^{\bxt}
 \widetilde{F}_2{}^{\bxt}F_3^{\bxt} @>>> F_3^{\bxt} \\
@V{\bC_{F_1,F_2}}VV @V{\text{via }}V{\bC_{F_1,F_2}}V 
 @V{\text{via }}V{\bC_{\widetilde{F}_1,\widetilde{F}_2}}V @| \\
F_3^{\bxt} @>>> \Delta^{\!\bxt}F_3'{}^{\bxt}F_3^{\bxt} @>>>
 \Delta^{\!\bxt}\widetilde{F}_3{}^{\bxt}F_3^{\bxt} @>>> F_3^{\bxt}
\end{CD}
\]
The square on the left commutes for functorial reasons, the one in the 
middle commutes by transitivity of base-change, and the one on the 
right commutes by~\ref{lem:recC0}. Thus the diagram commutes.
By reflexivity of $\Psi_{-,-}$ (\ref{prop:conseq3}(i)) the bottom row is 
the identity and hence the proposition follows by looking at the outer 
border of the diagram. 
\end{proof}

Next, we look at the fundamental isomorphism.

\begin{alem}
\label{lem:psi8b}
\emph{(Recovering $\Phi_{-}$).}
Let $\sigma$ be a sequence, such that $|\sigma|$ is an identity map,
say $1_X$. Set $I \set (1_X, \sfP)$. Then 
$\Phi_{\sigma} = \Psi_{I, \sigma}$.  
\end{alem}

\begin{proof}
By definition, $\Psi_{I, \sigma}$ is defined via the 
following diagram where we may assume that intermediate vertical
maps in the cartesian square $\sS$ are all $\sfP$-labeled identity maps.
\[
\xymatrix{
X \ar[r]^{\Delta} & X  \ar@{=}[d]_{I'=I}^{\hspace{3em} \sS} 
 \ar@{=>}[rr]^{\sigma'= \sigma} & & X \ar@{=}[d]^{I}  \\
& X \ar@{=>}[rr]_{\sigma}  & &  X 
}
\]
In this case the diagonal map is the identity so that 
$\Delta^{\!\bxt} = I^{\bxt} = \oneD{X}$. 
By looking at the outer border between the framed vertices below
we see that it suffices to prove that the following diagram commutes.
\[
\xymatrix{
\Delta^{\!\bxt}I'{}^{\bxt}\sigma^{\bxt} 
 \ar[dd]_{\text{via }\Phi_{\Delta \star I'}}^{\hspace{2em} \dag} 
 \ar@{=}[rd] \ar[rr]^{\text{via }\bbeta_{\sS}^{-1}} 
 & & \Delta^{\!\bxt}\sigma'{}^{\bxt}I^{\bxt} \ar@{=}[d]^{\hspace{4em}\ddag_1} 
 \ar[rr]^{\text{via }\Phi_{\Delta \star \sigma'}} & & I^{\bxt} \ar@{=}[d] \\
& I'{}^{\bxt}\sigma^{\bxt} \ar[r]^{\text{via }\bbeta_{\sS}^{-1}}
 & \sigma'{}^{\bxt}I^{\bxt} \ar@{=}[d]_{\ddag_2 \hspace{4em}} 
 \ar[rr]_{\text{via }\Phi_{\sigma'}} & & I^{\bxt} \ar@{=}[d] \\
\framebox{$\sigma^{\bxt}$} \ar@{=}[ru] \ar@{=}[rr] & & \sigma'{}^{\bxt} 
 \ar[rr]^{\Phi_{\sigma'}} & & \framebox{$\mathbf{1}_{\cD_X}$}
}
\]
Since $\Delta$ and $I'$ are $\sfP$-labeled, $\scriptstyle \dag$ commutes by
pseudofunctoriality of $(-)^{\times}$.
For~$\scriptstyle \ddag_2$, we use \ref{lem:cafpols20} while 
$\scriptstyle \ddag_1$ commutes by \ref{lem:canon10}. 
The remaining subdiagrams commute for functorial reasons.
\end{proof}

\subsection{Recovering the base-change isomorphisms}
\label{subsec:recoldbc}

Now, we move on to recovering the base-change isomorphism.
Again, we need some preliminaries. These are also used 
in~\S\ref{subsec:meiosis}.

\begin{alem} 
\label{lem:strcs1b}
Let $\sigma_i \colon X \Lra X$ for $i = 1,2$ be sequences such that 
$|\sigma_i| = 1_X$. Then the following diagram commutes.
\[
\xymatrix{
\sigma_1^{\bxt}\sigma_2^{\bxt} \ar[rr]^{\textup{via }\Phi_{\sigma_2}} 
 \ar[d]_{\textup{via }\Phi_{\sigma_1}} \ar[rrd]^{\Phi_{\sigma_1\star\sigma_2}} 
 & & \sigma_1^{\bxt} \ar[d]^{\Phi_{\sigma_1}} \\
\sigma_2^{\bxt} \ar[rr]_{\Phi_{\sigma_2}} & & \mathbf{1}_{\cD_X}
}
\]
\end{alem}
\begin{proof}
The outer diagram commutes for functorial reasons. For the lower triangle
we use \ref{lem:cocycle3} with $G$ as the identity map on $X$
and then use \ref{lem:psi8b} followed by \ref{lem:canon10}.
\end{proof}

\begin{alem}
\label{lem:recold1}
Consider a cartesian diagram as follows where $|\sigma_1|= 1_X$
and $|\sigma_3|= 1_Z$.
\[
\xymatrix{
Z \ar@{=>}[d]_{\sigma_2} \ar@{=>}[r]^{\sigma_3} & Z \ar@{=>}[d]^{\sigma_2} \\
X \ar@{=>}[r]^{\sigma_1} & X
}
\]
Then the following diagram commutes. 
\[
\xymatrix{
\sigma_2^{\bxt}\sigma_1^{\bxt} \ar[rd]_{\text{via } \Phi_{\sigma_1}} 
 \ar[rr]^{\bbeta_{\sigma_1, \sigma_2}} 
 & & \sigma_3^{\bxt}\sigma_2^{\bxt} \ar[ld]^{\text{via } \Phi_{\sigma_3}} \\
& \sigma_2^{\bxt} 
}
\]
\end{alem}
\begin{proof}
As a first step, we reduce, by induction on the length of $\sigma_2$,
to the case where $\sigma_2$ has length one. Let $n>1$ be an integer
and assume that the lemma has been verified when 
$\sigma_2$ has length less than~$n$. Now suppose $\sigma_2$
has length $n$. We can decompose $\sigma_2$ as $\sigma_4\star\sigma_5 \colon 
Z \stackrel{\sigma_4\;}{\Lra} Z_1 \stackrel{\sigma_5\;}{\Lra} X$.
Both, $\sigma_4, \sigma_5$ have length less than $n$.
Consider the following cartesian diagram,
\[
\xymatrix{
Z \ar@{=>}[d]_{\sigma_4} \ar@{=>}[r]^{\sigma_3} 
 & Z \ar@{=>}[d]^{\sigma_4} \\
Z_1 \ar@{=>}[d]_{\sigma_5} \ar@{=>}[r]^{\sigma_6} 
 & Z_1 \ar@{=>}[d]^{\sigma_5} \\
X \ar@{=>}[r]^{\sigma_1} & X 
}
\] 
corresponding to which, there results the following diagram of 
isomorphisms.
\[
\xymatrix{
\sigma_4^{\bxt}\sigma_5^{\bxt}\sigma_1^{\bxt} 
 \ar[rd]_{\text{via } \Phi_{\sigma_1}} 
 \ar[r]^{\text{via }\bbeta_{\sigma_1, \sigma_5}} 
 & \sigma_4^{\bxt}\sigma_6^{\bxt}\sigma_5^{\bxt}  
 \ar[r]^{\text{via }\bbeta_{\sigma_6, \sigma_4}} 
 \ar[d]_{\text{via }}^{\Phi_{\sigma_6}}
 & \sigma_3^{\bxt}\sigma_4^{\bxt}\sigma_5^{\bxt} 
 \ar[ld]^{\text{via } \Phi_{\sigma_3}} \\
\hspace{9em} & \sigma_4^{\bxt}\sigma_5^{\bxt} & \hspace{9em} 
}
\]
By the induction hypothesis, the two subtriangles 
commute and hence from the outer border we get the 
desired reduction.

Henceforth, we assume that $\sigma_2$ has length one. For 
carrying out the next step we need the following definition.

For any labeled sequence $\sigma$, given by, say,
\[
X_1 \xto{(f_1, \lambda_1)} X_2 \xto{(f_2, \lambda_2)} \quad \cdots \quad
\xto{(f_n, \lambda_n)} X_{n+1}
\]
we define its $\sfO$-depth by 
\[
\text{$\sfO$-depth}(\sigma) = 
\begin{cases} 
 0, &\text{if $\lambda_i = \sfP$ for all $i$;}  \\
 n+1- \inf\>\{i\;|\; \lambda_i \in \sfO \}, &\text{otherwise.}  
\end{cases}
\]

Returning to the proof of the lemma, we 
now argue by induction on the $\sfO$-depth of~$\sigma_1$.
Let us verify that the lemma holds 
if $\sigma_1$ has $\sfO$-depth~0. Indeed, first note that in 
this case, the entire staircase associated to~$\sigma_1$ consists 
of $\sfP$-labeled maps only and hence $\Phi_{\sigma_1}$ is the
same as the isomorphism $\sigma_1^{\bxt} \iso (1_X)^{\times} = \oneD{X}$ 
resulting from the pseudofunctoriality of~$(-)^{\times}$. 
A similar description holds for~$\Phi_{\sigma_3}$.
Therefore, using transitivity of base-change successively over the 
length of $\sigma_1, \sigma_3$ we reduce to proving the lemma when
$\sigma_1, \sigma_3$ are assumed to be identity maps with label $\sfP$. 
Now we conclude using \ref{lem:cafpols4}(i).

Now let $n>0$ be an integer, assume that 
the lemma holds whenever $\sigma_1$ has $\sfO$-depth less than $n$
and assume that $\sigma_1$ has $\sfO$-depth $n$.
Let us decompose
$\sigma_1$ as $\sigma_4 \star G \colon 
X \stackrel{\sigma_4\;}{\Lra} Y \xto{G} X$; here~$G$ is 
a labeled map. Note that $\sigma_4$ has $\sfO$-depth $n-1$. 
Set $F \set (|\sigma_4|,\sfP)$. Consider the
following cartesian diagram where the previously undefined maps
are determined by the cartesian condition.
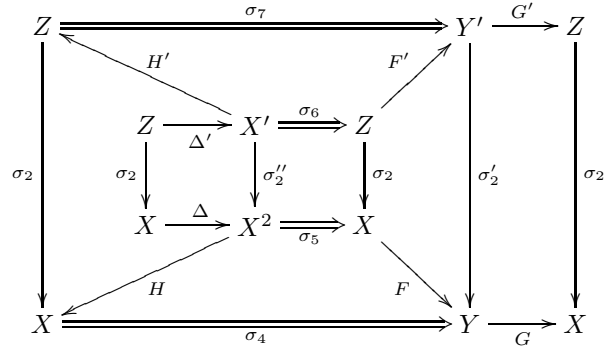
\[
\xymatrix{
Z \ar[ddd]_{\sigma_2}
 \ar@{=>}[rrrr]^{\sigma_7}  & & & &  Y' \ar[ddd]^{\sigma_2'} 
  \ar[r]^{G'} & Z \ar[ddd]^{\sigma_2} \\
 & Z \ar[d]_{\sigma_2} \ar[r]_{\Delta'} & X' \ar[d]^{\sigma_2''} 
 \ar@{=>}[r]^{\sigma_6} \ar[llu]_{H'} & Z \ar[d]^{\sigma_2} \ar[ru]^{F'} \\
 & X \ar[r]^{\Delta} & X^2 \ar@{=>}[r]_{\sigma_5} \ar[lld]^{H} 
 & X \ar[rd]_{F} \\
X \ar@{=>}[rrrr]_{\sigma_4} 
 & & & &  Y \ar[r]_{G} & X
}
\]
Correspondingly, we have the following diagram of isomorphisms
where the common superscript ${}^{\bxt}$ has been omitted from each
term. 
\begin{equation}
\label{eq:recold2}
\xymatrix{
\framebox{$\sigma_2\sigma_4G$} \ar[dd]^{\bbeta_4} 
 \ar[rd]^{\Phi_2} \ar[rr]^{\Phi_1} & & \sigma_2\Delta H\sigma_4G 
 \ar[d]_{\blacksquare_1\hspace{5em}\bbeta_2} \ar[rr]^{\bbeta_1} 
 & & \sigma_2\Delta\sigma_5FG \ar[ld]_{\bbeta_2} 
 \ar[dd]^{\Phi_5}_{\blacksquare_3\hspace{2em}} \\
& \Delta\!'H'\sigma_2\sigma_4G \ar[d]^{\bbeta_4\hspace{6.5em}\blacksquare_2} 
 \ar[r]^{\bbeta_3}
 & \Delta\!'\sigma_2''H\sigma_4G \ar[r]^{\bbeta_1} 
 & \Delta\!'\sigma_2''\sigma_5FG \ar[d]_{\bbeta_5} \\
\sigma_7\sigma_2'G \ar[d]^{\bbeta_6} \ar[r]^{\Phi_2} 
 & \Delta\!'H'\sigma_7\sigma_2'G \ar[d]^{\bbeta_6} \ar[r]^{\bbeta_7} 
 & \Delta\!'\sigma_6F'\sigma_2'G \ar[d]^{\bbeta_6} \ar[rd]^{\Phi_6} 
 \ar[r]^{\bbeta_8} & \Delta'\sigma_6\sigma_2FG \ar[r]^{\Phi_6} 
 & \sigma_2FG  \ar[d]^{\Phi_3} \\
\framebox{$\sigma_7G'\sigma_2$} \ar[r]_{\Phi_2} 
 & \Delta\!'H'\sigma_7G'\sigma_2 \ar[r]_{\bbeta_7} & 
 \Delta\!'\sigma_6F'G'\sigma_2 \ar[rd]_{\Phi_6} 
 & F'\sigma_2'G \makebox[0pt]{\hspace{6em}$\scriptstyle \blacksquare_4$} 
 \ar[ru]^{\bbeta_8} \ar[d]^{\bbeta_6} & \framebox{$\sigma_2$}  \\
& & & F'G'\sigma_2 \ar[ru]_{\Phi_4}
}
\end{equation}
The morphisms are denoted in the reduced notation using 
the following table.
\begin{align}
\Phi_1 &= \Phi_{\Delta\star H}^{-1} 
 &\Phi_5 &= \Phi_{\Delta\star \sigma_5} 
 &\bbeta_1 &= \bbeta_{\sigma_4,F}
 &\bbeta_5 &= \bbeta_{\sigma_5,\sigma_2}  \notag \\
\Phi_2 &= \Phi_{\Delta'\star H'}^{-1}
 &\Phi_6 &= \Phi_{\Delta'\star \sigma_6} 
 &\bbeta_2 &= \bbeta_{\Delta, \sigma_2''}
 &\bbeta_6 &= \bbeta_{G,\sigma_2} \notag \\
\Phi_3 &= \Phi_{F \star G} 
 & & 
 &\bbeta_3 &= \bbeta_{\sigma_2,H}
 &\bbeta_7 &= \bbeta_{\sigma_7,F'} \notag \\
\Phi_4 &= \Phi_{F' \star G'} 
 & & 
 &\bbeta_4 &= \bbeta_{\sigma_4,\sigma_2'}
 &\bbeta_8 &= \bbeta_{\sigma_2',F} \notag
\end{align}

Let us verify that \eqref{eq:recold2} commutes. 
The unnamed rectangles commute for functorial reasons.
By \ref{lem:cafpols4}(ii), we deduce commutativity of $\sbsq_4$ and $\sbsq_1$, 
where in the latter case we first cancel the common factor 
$\sigma_4^{\bxt}G^{\bxt}$ on the right from each vertex. 
In $\sbsq_2$, we cancel
$\Delta^{\!'\bxt}$ on the left and $G^{\bxt}$ on the right 
and then use the cube lemma. Finally, for $\sbsq_3$, we 
first cancel the common factor $F^{\bxt}G^{\bxt}$ on the right.
Now note that $\Delta\star\sigma_5$ has $\sfO$-depth $n-1$.
Therefore, by the induction hypothesis, commutativity follows.

To conclude the lemma consider the outer border 
of~\eqref{eq:recold2} redrawn, with restored superscripts, as follows.
\[
\begin{CD}
\framebox{$\sigma_2^{\bxt}\sigma_4^{\bxt}G^{\bxt}$} @>{\text{via}}>{\Phi_1}> 
 \sigma_2^{\bxt}\Delta^{\!\bxt} H^{\bxt}\sigma_4^{\bxt}G^{\bxt} 
 @>{\text{via}}>{\bbeta_1}> 
 \sigma_2^{\bxt}\Delta^{\!\bxt}\sigma_5^{\bxt}F^{\bxt}G^{\bxt} 
 @>{\text{via}}>{\Phi_5}> \sigma_2^{\bxt}F^{\bxt}G^{\bxt} \\
@V{\text{via}}V{\bbeta_4}V @. @. @V{\text{via}}V{\Phi_3}V \\
\sigma_7^{\bxt}\sigma_2^{'\bxt}G^{\bxt} @. @. @. 
 \framebox{$\sigma_2^{\bxt}$} \\
@V{\text{via}}V{\bbeta_6}V @. @. @A{\text{via}}A{\Phi_4}A \\
\framebox{$\sigma_7^{\bxt}G^{'\bxt}\sigma_2^{\bxt}$} @>{\text{via}}>{\Phi_2}> 
 \Delta^{\!'\bxt}H^{'\bxt}\sigma_7^{\bxt}G^{'\bxt}\sigma_2^{\bxt} 
 @>{\text{via}}>{\bbeta_7}> 
 \Delta^{\!'\bxt}\sigma_6^{\bxt}F^{'\bxt}G^{'\bxt}\sigma_2^{\bxt} 
 @>{\text{via}}>{\Phi_6}> F^{'\bxt}G^{'\bxt}\sigma_2^{\bxt} 
\end{CD}
\]
Observe that in the top row, upon the canceling of $\sigma_2^{\bxt}$ on 
the left and $G^{\bxt}$ on the right from each vertex, the resulting 
sequence composes to $\Psi_{F,\sigma_4}$. Therefore, by 
Lemma \ref{lem:cocycle3},
the composition $\sigma_2^{\bxt}\sigma_4^{\bxt}G^{\bxt} \xto{\text{top row}}
\sigma_2^{\bxt}F^{\bxt}G^{\bxt} \xto{\text{via }\Phi_3} \sigma_2^{\bxt}$
is the same as 
$\sigma_2^{\bxt}(\Phi_{\sigma_4\star G}) = \sigma_2^{\bxt}(\Phi_{\sigma_1})$.
By a similar argument we also see that the composition 
$\sigma_7^{\bxt}G^{'\bxt}\sigma_2^{\bxt} \xto{\text{bottom row}}
F^{'\bxt}G^{'\bxt}\sigma_2^{\bxt} \xto{\text{via }\Phi_4} \sigma_2^{\bxt}$ 
is the same as 
$\Phi_{\sigma_7\star G'}(\sigma_2^{\bxt}) = \Phi_{\sigma_3}(\sigma_2^{\bxt})$.

Since the preceding diagram commutes, from its outer border 
we deduce the lemma.
\end{proof}

\begin{alem}
\label{lem:recold2a}
Let $F \colon X \xto{\quad} Y$ be a $\sfP$-labeled map and 
${\sigma}\colon Y {\Lra} Y$ a sequence
such that $|\sigma| = 1_Y$.
Then the following two isomorphisms are equal
\[
F^{\bxt}\sigma^{\bxt} \xto{F^{\bxt}(\Phi_{\sigma})} F^{\bxt},
\qquad \qquad
F^{\bxt}\sigma^{\bxt} \xto{\Psi_{F, F\star\sigma}} F^{\bxt}. 
\] 
\end{alem}
\begin{proof}
We construct a cartesian diagram as follows where $\Delta$ is the diagonal map
with label $\sfP$.
\[
\xymatrix{
X \ar[r]^{\Delta} & X^2 \ar[d]_{P}^{\hspace{1.5em}\pfr} \ar[r]^{Q} 
 & X \ar[d]^{F\hspace{2.5em}\sfr} \ar@{=>}^{\sigma'}[rr] & & X \ar[d]^F \\
& X \ar[r]_{F} & Y \ar@{=>}[rr]_{\sigma} & & Y   
}
\]
It suffices to prove that the outer border of the following
diagram of isomorphisms commutes.
\[
\xymatrix{
F^{\bxt}\sigma^{\bxt} \ar[dd]_{\text{via }\Phi_{\sigma}} 
 & & & & \Delta^{\bxt}P^{\bxt}F^{\bxt}\sigma^{\bxt} 
 \ar[d]^{\text{via }\bbeta_{\pfr}^{-1}} 
 \ar[llll]_{\text{via }\Phi_{\Delta \star P}} \\
& \sigma'{}^{\bxt}F^{\bxt} \ar[lu]^{\bbeta_{\sfr}} 
 \ar[ld]^{\text{via }\Phi_{\sigma'}} 
 & & \ddag & \Delta^{\bxt}Q^{\bxt}F^{\bxt}\sigma^{\bxt} 
 \ar[d]^{\text{via }\bbeta_{\sfr}^{-1}}  
 \ar[llllu]_{\text{via }\Phi_{\Delta \star Q}} \\
F^{\bxt} & & & & \Delta^{\bxt}Q^{\bxt}\sigma'{}^{\bxt}F^{\bxt} 
 \ar[llll]^{\text{via }\Phi_{\Delta \star Q \star \sigma'}} 
 \ar[lllu]^{\text{via }\Phi_{\Delta \star Q}}  
}
\]
The triangle on the left commutes by \ref{lem:recold1}. The uppermost triangle
commutes by pseudofunctoriality since all the maps involved are
$\sfP$-labeled. Commutativity of~$\ddag$ is obvious.
Finally, the lowermost triangle commutes by \ref{lem:strcs1b}.
\end{proof}

Now we are in a position to show that $\Psi_{-,-}$ recovers 
the base change isomorphisms.

\begin{aprop}
\label{prop:recbeta}
\emph{(Recovering $\bbeta_{-,-}$)}.
Consider a cartesian diagram as follows.
\[
\xymatrix{
X \ar@{=>}[d]_{\sigma_3} \ar@{=>}[r]^{\sigma_4} 
 & Y_1 \ar@{=>}[d]^{\sigma_1}  \\
Y_2 \ar@{=>}[r]_{\sigma_2} & Z
}
\]
Then $\Psi_{\sigma_4\star\sigma_1, \sigma_3\star\sigma_2}
= \bbeta_{\sigma_2,\sigma_1}$.
\end{aprop}
\begin{proof}
Our goal may be rephrased as saying that corresponding to the 
cartesian diagram shown on the left below where 
$\rho_1 = \sigma_4\star\sigma_1$, $\rho_2 = \sigma_3\star\sigma_2$ and
$\Delta$ is the diagonal map with label~$\sfP$, the 
diagram of isomorphisms on the right commutes.
\begin{equation}
\label{eq:recold3}
\xymatrix{
X \ar[r]^{\Delta} & X^2 \ar@{=>}[d]_{\pi_1} \ar@{=>}[r]^{\pi_2} 
 & X \ar@{=>}[d]^{\rho_1}  \\
& X \ar@{=>}[r]_{\rho_2} & Z
}
\qquad \qquad
\xymatrix{
\rho_2^{\bxt} \ar[d]^{\bbeta_{\sigma_2,\sigma_1}} & & 
 \Delta^{\bxt}\pi_1^{\bxt}\rho_2^{\bxt} 
 \ar[ll]_{\text{via }}^{\Phi_{\Delta\star\pi_1}}  
 \ar[d]_{\text{via }}^{\bbeta_{\rho_2, \rho_1}}  \\
\rho_1^{\bxt} & & 
 \Delta^{\bxt}\pi_2^{\bxt}\rho_1^{\bxt}
 \ar[ll]_{\text{via }}^{\Phi_{\Delta\star\pi_2}}
}
\end{equation}

We expand and extend the diagram on the left in \eqref{eq:recold3} as follows.
Here the $\Delta_i$'s are the diagonal maps with label $\sfP$ and the 
rest are determined by the cartesian condition. Note that 
$|\Delta| = |\Delta_4\star\Delta_3'|$, $\pi_1 = \sigma_4''' \star \sigma_3'$
and $\pi_2 = \sigma_3''' \star \sigma_4'$.
\[
\xymatrix{
X \ar[d]_{\Delta_4} & & X \ar[d]_{\Delta_4} \\
Y_4 \ar@{=>}[d]_{\sigma_4''}  \ar[r]^{\Delta_3'} & X^2 
 \ar@{=>}[d]_{\sigma_4'''} \ar@{=>}[r]^{\sigma_3'''} 
 & Y_4\ar@{=>}[d]^{\sigma_4''} 
 \ar@{=>}[r]^{\sigma_4'} & X \ar@{=>}[d]^{\sigma_4} \\
X \ar[r]_{\Delta_3} & Y_3 \ar@{=>}[d]^{\sigma_3'} \ar@{=>}[r]^{\sigma_3''} 
 & X \ar@{=>}[d]^{\sigma_3} \ar@{=>}[r]^{\sigma_4} 
 & Y_1 \ar@{=>}[d]^{\sigma_1} \\
& X \ar@{=>}[r]_{\sigma_3} & Y_2 \ar@{=>}[r]_{\sigma_2} & Z 
}
\]
Correspondingly we expand the diagram on the right 
in \eqref{eq:recold3} as follows, where the superscript
${}^{\bxt}$ has been omitted for convenience.
\begin{equation}
\label{eq:recold4}
\xymatrix{
& & \framebox{$\Delta\pi_1\rho_2$} \ar@{=}[r]
 & \Delta\sigma_4'''\sigma_3'\sigma_3\sigma_2 
 \ar[d]^{\boldC_{\Delta_4,\Delta_3'}^{-1}} \\
\framebox{$\rho_2$} \ar@{=}[d] 
 & \Delta_3\sigma_3'\sigma_3\sigma_2 \ar[ld]_{\Phi_3'} 
 \ar[d]^{a}_{\raisebox{-3ex}{$\blacktriangle$} \hspace{1em}} 
 \ar[r]^{\Phi_4''\quad} &
 \Delta_4\sigma_4''\Delta_3\sigma_3'\sigma_3\sigma_2 \ar[d]^{a} 
 \ar[r]^{\bbeta_{\Delta_3,\sigma_4'''}} &
 \Delta_4\Delta_3'\sigma_4'''\sigma_3'\sigma_3\sigma_2 \ar[d]^{a} \\
\sigma_3\sigma_2 \ar[dd]^{\bbeta_{\sigma_2, \sigma_1}} \ar[rd]_{\Phi_4''}  
 & \Delta_3\sigma_3''\sigma_3\sigma_2 
 \ar[l]_{\Phi_3''} \ar[r]^{\Phi_4''\quad} &
 \Delta_4\sigma_4''\Delta_3\sigma_3''\sigma_3\sigma_2  
 \ar[ld]^{\Phi_3''\makebox[0pt]{\hspace{7em}\smash{\raisebox{2ex}{$
 \scriptstyle \blacksquare$}}}} 
 \ar[r]^{\bbeta_{\Delta_3,\sigma_4'''}} &
 \Delta_4\Delta_3'\sigma_4'''\sigma_3''\sigma_3\sigma_2 
 \ar[ld]^{\bbeta_{\sigma_3'', \sigma_4''}} \\
\hspace{6em} & \Delta_4\sigma_4''\sigma_3\sigma_2 
 \ar[d]^{\bbeta_{\sigma_2, \sigma_1}} \ar[r]_{\Phi_3'''\quad} 
 & \Delta_4\Delta_3'\sigma_3'''\sigma_4''\sigma_3\sigma_2 
 \ar[d]^{\bbeta_{\sigma_2, \sigma_1}} \\
\sigma_4\sigma_1 \ar@{=}[d] \ar[rd]_{\Phi_4'} \ar[r]^{\Phi_4''\quad} 
 & \Delta_4\sigma_4''\sigma_4\sigma_1 
 \ar[d]^{b}_{\raisebox{3ex}{$\blacktriangledown$} \hspace{1em}} 
 \ar[r]^{\Phi_3'''\quad} &
 \Delta_4\Delta_3'\sigma_3'''\sigma_4''\sigma_4\sigma_1 \ar[d]^{b} 
 & \framebox{$\Delta\pi_2\rho_1$} \ar@{=}[d] \\
\framebox{$\rho_1$} &
 \Delta_4\sigma_4'\sigma_4\sigma_1 \ar[r]^{\Phi_3'''\quad} &
 \Delta_4\Delta_3'\sigma_3'''\sigma_4'\sigma_4\sigma_1 
 \ar[r]^{\boldC_{\Delta_4,\Delta_3'}} &
 \Delta\sigma_3'''\sigma_4'\sigma_4\sigma_1 
}
\end{equation}
The maps are all isomorphisms. Reduced notation has been employed 
for all the maps except those denoted by $a$ or $b$, while the $\Phi$'s
are given by the following table.
\[
\Phi_3' = \Phi_{\Delta_3 \star \sigma_3'}, \;\;\;
\Phi_3'' = \Phi_{\Delta_3 \star \sigma_3''}, \;\;\;
\Phi_3''' = \Phi_{\Delta_3'\star \sigma_3'''}^{-1}, \;\;\;
\Phi_4' = \Phi_{\Delta_4 \star \sigma_4'}^{-1}, \;\;\;
\Phi_4'' = \Phi_{\Delta_4 \star \sigma_4''}^{-1}.
\]
The ones denoted by $a$ are induced in the obvious way by the 
base-change isomorphism
$\sigma_3'{}^{\bxt}\sigma_3^{\bxt} \iso
\sigma_3''{}^{\bxt}\sigma_3^{\bxt}$, while the ones denoted
by~$b$ are induced in the obvious way by the 
base-change isomorphism
$\sigma_4''{}^{\bxt}\sigma_4^{\bxt} \iso
\sigma_4'{}^{\bxt}\sigma_4^{\bxt}$.

Let us verify that \eqref{eq:recold4} commutes. Of all its subdiagrams,
only $\blacktriangle,\blacktriangledown$ and~$\sbsq$ need be considered; 
the remaining ones commute for functorial reasons. 

In $\sbsq$, we may cancel off $\Delta_4^{\!\bxt}$
on the left and $\sigma_3^{\bxt}\sigma_2^{\bxt}$ on the right, 
from each vertex. The resulting diagram now commutes by~\ref{lem:recold1}.
In $\blacktriangle$ we first cancel 
$\sigma_2^{\bxt}$ on the right. The resulting diagram commutes 
by \ref{prop:conseq3}(i). By a similar argument
$\blacktriangledown$ commutes. Thus \eqref{eq:recold4} commutes.

Therefore, to prove that the diagram 
on the right in \eqref{eq:recold3} commutes,
it suffices to match it with the outer border of~\eqref{eq:recold4},
comprising of the framed vertices and the maps between them. 
Here is the list of the requirements wherein three edges
from~\eqref{eq:recold3} are to be compared with maps obtained
from the outer border of~\eqref{eq:recold4}, the fourth one being 
shared by both.
\begin{align}
\text{(i)}&.\;\; \Phi_{\Delta\star\pi_2}(\rho_1^{\bxt}) :
 \Delta^{\!\bxt}\pi_2^{\bxt}\rho_1^{\bxt} \xto{\quad} \rho_1^{\bxt}  
 & &= \text{ via outer border of \eqref{eq:recold4} ?} \label{eq:recold5} \\
\text{(ii)}&.\;\; \Phi_{\Delta\star\pi_1}(\rho_2^{\bxt}) :
 \Delta^{\!\bxt}\pi_1^{\bxt}\rho_2^{\bxt} \xto{\quad} \rho_2^{\bxt}
 & &= \text{ via outer border of \eqref{eq:recold4} ?} \notag \\
\text{(iii)}&.\;\; \Delta^{\!\bxt}(\bbeta_{\rho_2,\rho_1}) :
 \Delta^{\!\bxt}\pi_1^{\bxt}\rho_2^{\bxt} \xto{\quad} 
 \Delta^{\!\bxt}\pi_2^{\bxt}\rho_1^{\bxt} 
 & &= \text{ via outer border of \eqref{eq:recold4} ?} \notag
\end{align}
Out of the three, we give proofs for (i) and (ii) below;
(iii) follows essentially by functorial considerations.

In (i), $\rho_1^{\bxt}$ is common to each object
and moreover, in the corresponding portion of~\eqref{eq:recold4},
every vertex contains $\rho_1^{\bxt}=\sigma_4^{\bxt}\sigma_1^{\bxt}$. 
Upon canceling these terms, (i) reduces to 
checking that the outer border of the following diagram commutes.
\[
\xymatrix{
\framebox{$\Delta^{\bxt}\pi_2^{\bxt}$} \ar@{=}[d]^{\hspace{10em}\ddag_1} 
 \ar[rrrr]^{\Phi_{\Delta\star\pi_2}} & & & & 
 \framebox{$\mathbf{1}_{\cD_X}$} \\
\Delta^{\bxt}\sigma_3'''{}^{\bxt}\sigma_4'{}^{\bxt} 
\ar[rr]^{\text{via }\quad}_{\boldC_{\Delta_4,\Delta_3'}^{-1}\quad} & & 
 \Delta_4^{\bxt}\Delta'_3{}^{\bxt}\sigma_3'''{}^{\bxt}\sigma_4'{}^{\bxt}
 \ar[rru]^{\Phi_{...}} 
 \ar[rr]^{\quad\text{via }}_{\quad\quad\Phi_{\Delta'_3\star\sigma_3'''}} 
 & & \Delta_4^{\bxt}\sigma_4'{}^{\bxt} 
 \ar[u]_{\Phi_{\Delta_4\star\sigma'_4}}^{\ddag_2\hspace{1em}}
}
\] 
Here 
$\Phi_{...} = \Phi_{\Delta_4\star\Delta'_3\star\sigma_3'''\star\sigma_4'}$.
Now, in the bottom row, both the maps can be rewritten as
\[
\text{via }\boldC_{\Delta_4,\Delta_3'}^{-1} = 
\text{via } \Psi_{\Delta_4\star\Delta'_3, \;\Delta}, \qquad 
\text{via }\Phi_{\Delta'_3\star\sigma_3'''} =
\text{via } \Psi_{\Delta_4,\;\Delta_4\star\Delta'_3\star\sigma_3'''}, 
\]
where for the former we use \ref{prop:recC}, while for the latter we use
\ref{lem:recold2a}. With this description, both $\ddag_1$ and $\ddag_2$
are seen to commute by~\ref{lem:cocycle3}. 

In (ii), we begin as in (i) by first canceling the common factor,
which in this case is $\rho_2^{\bxt}$ or $\sigma_3^{\bxt}\sigma_2^{\bxt}$.
Thereupon, (ii) reduces to 
checking that the outer border of the following diagram 
with omitted superscripts~${}^{\bxt}$ commutes. 
\[
\xymatrix{
& \framebox{$\Delta\pi_1$} \ar@{=}[ld] 
 \ar[rd]_{\text{via } \boldC_{\Delta_4,\Delta_3'}^{-1}}^{\hspace{4em}\dag_1} 
 \ar[rr]^{\Phi} & & \framebox{$\mathbf{1}_{\cD_X}$} \\
\Delta\sigma_4'''\sigma_3' 
 \ar[rd]_{\text{via } \boldC_{\Delta_4,\Delta_3'}^{-1}} & 
 \makebox[0pt]{$\scriptstyle \dag_2$\hspace{3em}} & 
 \Delta_4\Delta'_3\pi_1 \ar@{=}[ld] \ar[ru]_{\Phi} & &
 \Delta_3\sigma_3' \ar[lu]_{\Phi} \\
& \Delta_4\Delta'_3\sigma_4'''\sigma_3' 
 \ar[rr]_{\text{via }\bbeta_{\sigma_4''',\Delta_3}} & & 
 \Delta_4\sigma_4''\Delta_3\sigma_3' 
 \ar[uu]_{\Phi\hspace{1em}\dag_4}^{\dag_3\hspace{2em}} 
 \ar[ru]_{\text{via }\Phi_{\Delta_4\star\sigma_4''}}
}
\]
In the four places where $\Phi$ occurs without a subscript, the 
missing subscript is the entire sequence given by the source of 
the corresponding map.
By \ref{lem:cocycle3}, $\dag_1$ commutes 
while~$\dag_2$ commutes for functorial reasons.
For $\dag_3$, we use \ref{lem:cocycle4} 
while~$\dag_4$ commutes by \ref{lem:strcs1b}. 
Thus \eqref{eq:recold5}(ii) is verified.

This concludes the proof of the commutativity of the 
diagram on the right in~\eqref{eq:recold3}. 
\end{proof}

\subsection{Linearity}
\label{subsec:meiosis}
We show that $\Psi_{-,-}$ is linear in the following sense.

\begin{athm}
\label{thm:meiosis1}
Let $\sigma_i \colon X \Lra Y$ for $i=1,2$ be two sequences such that
$|\sigma_1| = |\sigma_2|$. Then for any two sequences $\rho \colon W \Lra X$
and $\tau \colon Y \Lra Z$, the following two isomorphisms
are equal
\[
\rho^{\bxt}\sigma_2^{\bxt}\tau^{\bxt} 
\xto{\;\Psi_{\rho\star\sigma_1\star\tau,\;\rho\star\sigma_2\star\tau}\;} 
\rho^{\bxt}\sigma_1^{\bxt}\tau^{\bxt}, \qquad\qquad 
\rho^{\bxt}\sigma_2^{\bxt}\tau^{\bxt} 
\xto{\;\text{via }\Psi_{\sigma_1,\sigma_2}\;}
\rho^{\bxt}\sigma_1^{\bxt}\tau^{\bxt}.
\]
\end{athm}
\begin{proof}
It suffices to prove the theorem for the two special cases when 
either~$\rho$ or~$\tau$ is empty. In other words, it suffices to prove 
that the following hold
\[
\Psi_{\rho\star\sigma_1,\>\rho\star\sigma_2} = 
\rho^{\bxt}(\Psi_{\sigma_1,\sigma_2}), \qquad \qquad
\Psi_{\sigma_1\star\tau,\>\sigma_2\star\tau}
= \Psi_{\sigma_1,\sigma_2}(\tau^{\bxt}).
\]
We only show the first of these (linearity on the left). 
The second one is proved using similar methods.

Consider the following two cartesian diagrams.
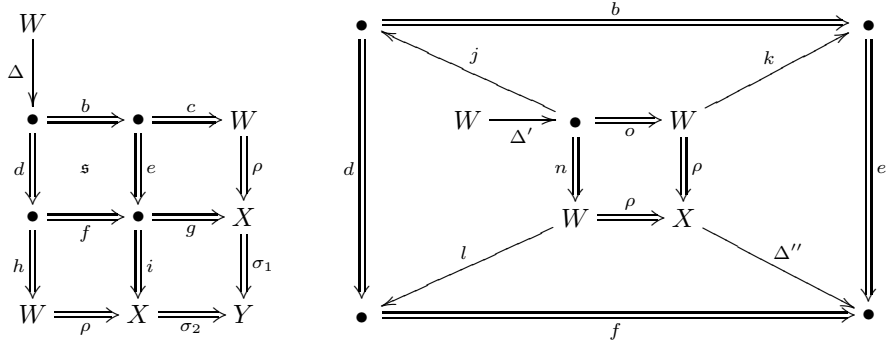
\[
\xymatrix{
W \ar[d]_{\Delta} \\
\bullet \ar@{=>}[d]_d^{\hspace{1.5em} \sfr} \ar@{=>}[r]^b 
 & \bullet \ar@{=>}[d]^e \ar@{=>}[r]^c & W \ar@{=>}[d]^{\rho} \\
\bullet \ar@{=>}[d]_h \ar@{=>}[r]_f & \bullet \ar@{=>}[d]^i \ar@{=>}[r]_g 
 & X \ar@{=>}[d]^{\sigma_1} \\
W \ar@{=>}[r]_{\rho} & X \ar@{=>}[r]_{\sigma_2} & Y 
}
\qquad
\xymatrix{
{\phantom{l}\!\!\bullet\!\!\phantom{l}} 
 \ar@{=>}[ddd]_{d}\ar@{=>}[rrrrr]^{b} & & & & &  
 {\phantom{l}\!\!\bullet\!\!\phantom{l}} \ar@{=>}[ddd]^{e} \\
& W \ar[r]_{\Delta'} & {\phantom{l}\!\!\bullet\!\!\phantom{l}} 
 \ar[llu]_{j} \ar@{=>}[d]_{n} \ar@{=>}[r]_{o} 
 & W \ar@{=>}[d]^{\rho} \ar[rru]^{k} \\
& & W \ar[lld]_{l} \ar@{=>}[r]^{\rho} 
 & X \ar[rrd]^{\Delta''} & & \\
{\phantom{l}\!\!\bullet\!\!\phantom{l}}
  \ar@{=>}[rrrrr]_{f} & & & & & \bullet 
}
\]
The diagram on the right is a cube built on $\sfr$, the northwestern square 
in the diagram on the left. Thus, the cube is obtained as a fibered product
of $\Delta''$ with~$\sfr$. (As usual, $\Delta, \Delta', \Delta''$ are 
diagonal maps with label $\sfP$.) In particular, 
$|\Delta'\star j| = |\Delta|$. 

Consider the following diagram of isomorphisms where the common superscript 
${}^{\bxt}$ has been omitted for convenience.
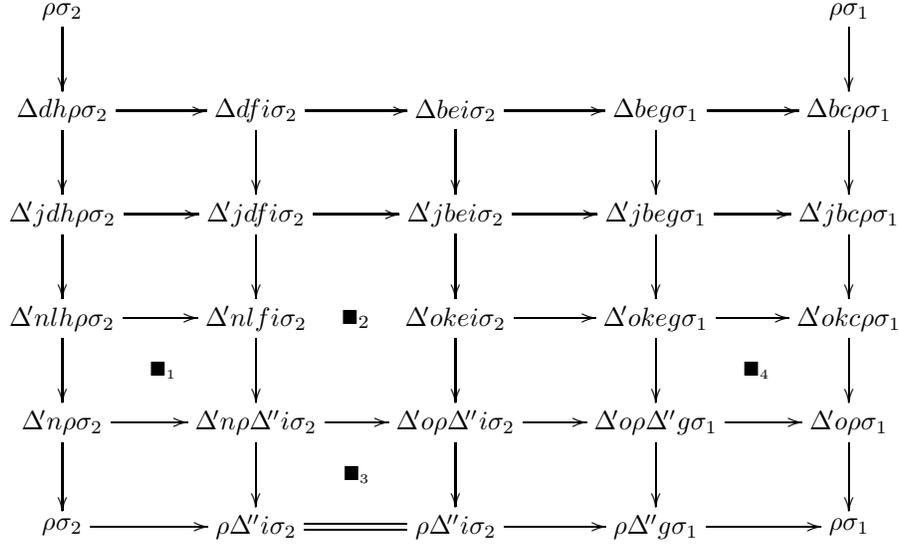
\[
\xymatrix{
\rho\sigma_2 \ar[d] & & & & \rho\sigma_1 \ar[d] \\ 
\Delta dh\rho\sigma_2 \ar[r] \ar[d] & \Delta dfi\sigma_2 \ar[r] \ar[d] & 
 \Delta bei\sigma_2 \ar[r] \ar[d] & \Delta beg\sigma_1 \ar[r] \ar[d] & 
 \Delta bc\rho\sigma_1 \ar[d] \\
\Delta\!'jdh\rho\sigma_2 \ar[r] \ar[d] & \Delta\!'jdfi\sigma_2 \ar[r] \ar[d] & 
 \Delta\!'jbei\sigma_2 \ar[r] \ar[d] & \Delta\!'jbeg\sigma_1 \ar[r] \ar[d] & 
 \Delta\!'jbc\rho\sigma_1 \ar[d] \\
\Delta\!'nlh\rho\sigma_2 \ar[r] \ar[d]^{\hspace{3em} \sbsq_1} 
 & \Delta\!'nlfi\sigma_2 \ar[d] & \Delta\!'okei\sigma_2 \ar[r] \ar[d] 
 & \Delta\!'okeg\sigma_1 \ar[r] \ar[d]^{\hspace{3em} \sbsq_4}  
 & \Delta\!'okc\rho\sigma_1 \ar[d] \\
\Delta\!'n\rho\sigma_2 \ar[r] \ar[d] & \Delta\!'n\rho\Delta\!''i\sigma_2 
 \ar[r]^{\smash{\raisebox{8ex}{$\sbsq_2$}}} \ar[d]^{\hspace{3em} \sbsq_3} 
 & \Delta\!'o\rho\Delta\!''i\sigma_2 \ar[r] \ar[d] 
 & \Delta\!'o\rho\Delta\!''g\sigma_1 \ar[r] \ar[d] & 
 \Delta\!'o\rho\sigma_1 \ar[d] \\
\rho\sigma_2 \ar[r] & \rho\Delta\!''i\sigma_2 \ar@{=}[r] & 
 \rho\Delta\!''i\sigma_2 \ar[r] & \rho\Delta\!''g\sigma_1 \ar[r] 
 & \rho\sigma_1 \\
}
\]
The morphisms are all obvious ones induced by $\boldC$,
$\bbeta$ or $\Phi$ for suitable
superscript that can be determined in an obvious way for each edge. Using this
diagram one deduces the theorem as follows.

First we verify that the preceding diagram commutes.
The unnamed rectangles commute for functorial reasons.
In $\sbsq_1$, upon canceling $\Delta\!'{}^{\bxt}n^{\bxt}$ on the left
and $\sigma_2^{\bxt}$ on the right from each vertex, we are left 
with a diagram whose commutativity follows from \ref{lem:recold1}
using the following substitutions.
\[
\sigma_1 \leadsto \Delta''\star i, \qquad \sigma_3 \leadsto l \star h,
\qquad \sigma_2 \leadsto \rho, \quad \text{etc.}
\]
By a similar argument $\sbsq_4$ commutes.
For $\sbsq_2$ we use the cube lemma. Finally, for $\sbsq_3$, 
we first cancel $\Delta''{}^{\bxt}i^{\bxt}\sigma_2^{\bxt}$ on the right
and then conclude by the reflexivity property $\Psi_{\rho,\rho} = $ identity.

Now note that that the topmost row of maps
in the preceding diagram defines 
$\Psi_{\rho\star\sigma_1,\>\rho\star\sigma_2}$ while the bottommost row 
defines $\rho^{\bxt}(\Psi_{\sigma_1,\sigma_2})$. Therefore 
to prove the theorem, it suffices to verify 
the leftmost column of maps and the rightmost one,
each composes to the corresponding identity transformation. 
Each of these verifications is carried out in the same way.
For instance, for the left column, we may first cancel 
$\rho^{\bxt}\sigma_2^{\bxt}$ from each object and then
use the outer border of the following diagram where 
each unnamed arrow is $\Phi_{-}$ with suitable subscript.
\[
\xymatrix{
\Delta^{\bxt}d^{\bxt}h^{\bxt} \ar[rr]
 & & \mathbf{1}_{\cD_X} & &
 \ar[ll] \Delta'{}^{\bxt}n^{\bxt} \\
\Delta'{}^{\bxt}j^{\bxt}d^{\bxt}h^{\bxt} 
 \ar[u]^{\text{via }\boldC_{\Delta',j}} \ar[rru]\ar[rrrr]
 & & & & \Delta'{}^{\bxt}n^{\bxt}l^{\bxt}h^{\bxt} 
 \ar[u]_{\text{via }\Phi_{l\star h}} \ar[llu]
}
\]
The triangle on the left commutes by \ref{lem:cocycle3}
(because $\boldC_{\Delta',j} = \Psi_{\Delta,\; \Delta' \star j}$),
the one in the middle commutes by 
\ref{lem:cocycle4} and the one on the right commutes by \ref{lem:strcs1b}.
\end{proof}

\section{Proofs IV (the output)}
\label{sec:proofsIV}
We give proofs for Theorem \ref{thm:output4}, 
Theorem \ref{thm:flat5}, and Theorem \ref{thm:variant3}. 
The proofs are based on the results 
of \S\ref{sec:proofsI}--\S\ref{sec:proofsIII}.
 
\subsection{Preliminaries on pseudofunctorial covers}
We discuss some basic results concerning pseudofunctorial covers 
(\ref{def:output1})
that we shall need in the proof of Theorem~\ref{thm:output4}.
Henceforth, we refer to a pseudofunctorial cover simply as a cover.

\ssss
\label{sssec:output4a}
We use the following notation.
Let $\eC = ((-)^{!}, \mu^{!}_{\times},\mu^{!}_{\ssbox})$ be a cover.
For any object $X$ in $\sfC$, let $S_X \colon X^! \to \cD_X$ be the 
associated functor.   
\begin{enumerate}
\item For any labeled map $F = (f,\lambda)\colon X \to Y$, 
we set $F^{!} \set f^!$ and for any labeled sequence $\sigma$ as follows
\[
X_1 \xto{\;F_1\;} X_2 \xto{\;F_2\;} X_3 
\xto{\quad} \cdots \xto{\quad}
X_{n} \xto{\;F_n\;} X_{n+1},
\]
we set $\sigma^{!} = F_1^! \cdots F_n^!$.
\item To every labeled map $F$ as in (i), we associate a map 
$S_XF^{!} \to F^{\bxt}S_Y$,
which, for $\lambda = \sfP$, is the canonical map 
$S_Xf^{!} \to f^{\times}S_Y$ obtained 
from $\mu^{!}_{\times}$ and for $\lambda = \sfO$, 
is the canonical map 
$S_Xf^{!} \to f^{\ssbox}S_Y$ obtained 
from~$\mu^{!}_{\ssbox}$.
To every sequence $\sigma$ as in (i), we associate the map 
$S_{X_1}\sigma^{!} \to \sigma^{\bxt}S_{X_{n+1}}$
which is the obvious one obtained by successively using 
$S_{X_i}F_i^{!} \to F_i^{\bxt}S_{X_{i+1}}$.
\end{enumerate}

It follows that if $\eC$ is a perfect cover, (\ref{def:output3}) then 
$S_{X_1}\sigma^{!} \to \sigma^{\bxt}S_{X_{n+1}}$, as defined in~(ii) 
above, is an isomorphism.

\begin{alem}
\label{lem:output5}
Let $\eC = ((-)^{!}, \mu^{!}_{\times},\mu^{!}_{\ssbox})$ be a cover.
\begin{enumerate}
\item For any labeled sequence $\sigma$ such that $|\sigma|$ is 
an identity map, say $1_X$, the following diagram of obvious
natural maps commutes.
\[
\xymatrix{
S_X\sigma^{!} \ar[rd]_{\textup{via }(-)^!} \ar[rr] 
 & & \sigma^{\bxt}S_X \ar[ld]^{\textup{via }\Phi_{\sigma}} \\
& S_X
}
\]
\item For any cartesian square $\sS$ as follows,
\[
\xymatrix{
W \ar@{=>}[d]_{\sigma_3} \ar@{=>}[r]^{\sigma_4} & X \ar@{=>}[d]^{\sigma_1} \\
Z \ar@{=>}[r]_{\sigma_2} & Y
}
\]
the following diagram of obvious natural maps commutes.
\[
\begin{CD}
S_W\sigma_4^!\sigma_1^! @>{\textup{via }(-)^!}>> S_W\sigma_3^!\sigma_2^! \\
@VVV @VVV \\ 
\sigma_4^{\bxt}S_X\sigma_1^! @. \sigma_3^{\bxt}S_Z\sigma_2^! \\
@VVV @VVV \\ 
\sigma_4^{\bxt}\sigma_1^{\bxt}S_Y @>{\textup{via }\bbeta_{\sS}}>> 
 \sigma_3^{\bxt}\sigma_2^{\bxt}S_Y  
\end{CD}
\]
\end{enumerate}
\end{alem}
\begin{proof}
(Sketch) For simplicity we shall assume that for all $X$, $X^! = \cD_X$
and that $S_X$ is the identity functor $\oneD{X}$. 
The proof in the general case is obtained by inserting $S_{-}$'s
in the appropriate places.

For (i), one looks at a staircase $\cS$ based on $\sigma$. 
Consider the following diagram whose bottom row 
spells out the definition of $\Phi_{\sigma}$ as in~\eqref{eq:psi3},
the remaining maps being the obvious ones.
\[
\begin{CD}
\sigma^! @>>> V_1^!\sigma^! @>>> \texttt{Steps}^! @>>> \oneD{X} \\ 
@VV{\hspace{2em} \ddag_1}V @VV{\hspace{2.5em} \ddag_2}V 
 @VV{\hspace{2.5em} \ddag_3}V @| \\
\sigma^{\bxt} @>>> V_1^{\bxt}\sigma^{\bxt} @>>> \texttt{Steps}^{\bxt} 
 @>>> \oneD{X}\end{CD}
\] 
To prove (i) it suffices to check that this diagram commutes.

We may expand $\scriptstyle \ddag_1$ as follows. 
\[
\xymatrix{
\sigma^! \ar[dd] \ar[rr] & & V_1^!\sigma^! \ar[d] \\
& & V_1^!\sigma^{\bxt} \ar[d] \\
\sigma^{\bxt} \ar[rr] \ar[rru] & & V_1^{\bxt}\sigma^{\bxt}
}
\]
The trapezium commutes for functorial reasons. Since $V_1$ consists 
of $\sfP$-labeled maps only, we deduce from the pseudofunctoriality
of the isomorphism $\mu^{!}_{\times}$ that the triangle in the preceding
diagram also commutes. 

In $\scriptstyle \ddag_2$, the bottom arrow is given by base-change
isomorphisms associated to each of the $n(n-1)/2$ squares occurring in the
staircase $\cS$. By using the compatibility in~\ref{def:output1}(b)  
for each of these squares one obtains that $\scriptstyle \ddag_2$ commutes.
In $\scriptstyle \ddag_3$, one uses \ref{def:output1}(c) for each of the 
$n$ steps in $\cS$ to deduce commutativity.

Finally, (ii) is proved using similar arguments.
\end{proof}

\begin{aprop}
\label{prop:output6}
Let $\eC = ((-)^{!}, \mu^{!}_{\times},\mu^{!}_{\ssbox})$ be a cover.
For any two sequences $\sigma_1,\sigma_2$ such that 
$|\sigma_1| = |\sigma_2|$, the following diagram of obvious natural 
maps commutes where $X$
is the source of $\sigma_i$ and $Y$ the target.
\[
\begin{CD}
S_X\sigma_1^! @>{\textup{via }(-)^!}>> S_X\sigma_2^! \\
@VVV @VVV \\
\sigma_1^{\bxt}S_Y @>{\textup{via }\Psi_{\sigma_2,\sigma_1}}>> 
 \sigma_2^{\bxt}S_Y 
\end{CD}
\] 
\end{aprop}
\begin{proof}
Let us assume for simplicity that for all $X$, $X^! = \cD_X$
and that $S_X$ is the identity functor $\oneD{X}$.
We use the diagram defining $\Psi_{\sigma_1,\sigma_2}$ 
(see \ref{sssec:psi5z}). It suffices to check that the 
following diagram commutes.
\[
\begin{CD}
\sigma_2^! @>>> \Delta^{!}\sigma_1'{}^{!}\sigma_2^! 
 @>>> \Delta^{!}\sigma_2'{}^{!}\sigma_1^! @>>> \sigma_1^! \\
@VVV @VVV @VVV @VVV \\
\sigma_2^{\bxt} @>>> \Delta^{\bxt}\sigma_1'{}^{\bxt}\sigma_2^{\bxt} 
 @>>> \Delta^{\bxt}\sigma_2'{}^{\bxt}\sigma_1^{\bxt} @>>> \sigma_1^{\bxt} \\
\end{CD}
\]
The square in the middle commutes by \ref{lem:output5}(ii)
while the commutativity of the squares on the either ends 
follows easily from \ref{lem:output5}(i).
\end{proof}

\subsection{Proof of Theorem \ref{thm:output4}.}
\label{subsec:Thm1}

\ssss
\label{sssec:output7} 
Proof of part (i):
We construct a perfect cover 
$\eC = ((-)^{!}, \mu^{!}_{\times},\mu^{!}_{\ssbox})$ as follows.

For any object $X$ in $\sfC$, set $X^! \set \cD_X$ and set
$S_X \set \oneD{X}$. For every map $f$ in~$\sfQ$, choose a labeled
sequence $\sigma_{\!f}$ such that $|\sigma_{\!f}| = f$ 
with the understanding that
if~$f$ is an identity map, say $1_X$, then $\sigma_{\!f}$
is a single map, viz., $(1_X, \sfP)$.
For any map~$f$ in~$\sfQ$, set $f^! \set \sigma_{\!f}^{\bxt}$. 
For any pair of maps $X \xto{f} Y \xto{g} Z$ in $\sfQ$, 
set $C^!_{f,g} = \Psi_{\sigma_{\<g\<f}, 
\;\sigma_{\<f}\> \star \>\>\sigma_{\<g}}$.  

Let us verify that $(-)^!$ is a pseudofunctor. We only check the 
associativity property of~$C^!_{-,-}$. Let $X \xto{f} Y \xto{g} Z \xto{h} W$
be maps in~$\sfQ$. Our aim is to verify that 
\[
C^!_{gf,h} \circ C^!_{f,g}(h^!) = C^!_{f,hg} \circ f^!(C^!_{g,h}).
\]
By definition, this amounts to verifying that 
\begin{equation}
\label{eq:output8}
\Psi_{\sigma_{\<h\<g\<f}, \;\sigma_{\<g\<f} \> \star \>\> \sigma_{\<h}}
\circ 
\Psi_{\sigma_{\<g\<f}, 
\;\sigma_{\<f} \> \star \>\> \sigma_{\<g}}(\sigma_{\<h}^{\bxt})
=
\Psi_{\sigma_{\<h\<g\<f}, \;\sigma_{\<f} \> \star \>\> \sigma_{\<h\<g}}
\circ 
\sigma_{\<f}^{\bxt}
(\Psi_{\sigma_{\<h\<g}, \;\sigma_{\<g} \> \star \>\>\sigma_{\<h}}). 
\end{equation}
To that end, we verify that each side equals 
$\Psi_{\sigma_{\<h\<g\<f},
\;\sigma_{\<f} \> \star \>\> \sigma_{\<g} \> \star \>\> \sigma_{\<h}}$.
For the left-hand side, we see that by~\ref{thm:meiosis1}, 
\[
\Psi_{\sigma_{\<g\<f}, 
\;\sigma_{\<f} \> \star \>\> \sigma_{\<g}}(\sigma_{\<h}^{\bxt})
=
\Psi_{\sigma_{\<g\<f} \> \star \>\> \sigma_{\<h},
\;\sigma_{\<f} \> \star \>\> \sigma_{\<g} \> \star \> \>\sigma_{\<h}}
\]
and therefore, by the cocycle rule, (\ref{thm:cocycle1}) the left-hand side 
of~\eqref{eq:output8} is the same as  
$\Psi_{\sigma_{\<h\<g\<f},
\;\sigma_{\<f} \> \star \>\> \sigma_{\<g} \> \star \>\> \sigma_{\<h}}$.
A similar argument works for the right-hand side.
Thus $C^!_{-,-}$ is associative.

Next we define on $\sfP$, a pseudofunctorial isomorphism 
$\mu^{!}_{\times} \colon (-)^!\big|_{\sfP} \iso (-)^{\times}$ as 
follows. For any $f\colon X \to Y$ in $\sfP$, we define the required map 
$f^! \iso f^{\times}$ to be $\Psi_{F,\>\sigma_{\<f}}$ where $F \set (f,\sfP)$.
For this to give a pseudofunctorial isomorphism it remains to be verified 
that for any pair of $\sfP$-morphisms $X \xto{f} Y \xto{g} Z$, 
the following diagram commutes where $F = (f,\sfP), G = (g,\sfP),
H = (gf,\sfP)$.
\[ 
\begin{CD}
f^{!}g^{!} @>{\text{via } \Psi_{F,\sigma_{\<f}}}>a> 
 f^{\times}g^{!} @>{\text{via } \Psi_{G,\sigma_{\<g}}}>b> 
 f^{\times}g^{\times} \\
@V\Psi_{\sigma_{\<g\<f}, \sigma_{\<f} \star \>\sigma_{\<g}}VV 
 @. @VcV{C^{\times}_{f,g}}V \\
(gf)^{!} @. 
 \makebox[0pt]{$\xto{\hspace{4em}\Psi_{H,\sigma_{\<g\<f}}\hspace{4em}}$} 
 @. (gf)^{\times}
\end{CD}
\]
By \ref{thm:meiosis1}, the maps $a$ and $b$ in the top row
equal $\Psi_{F\star\>\sigma_{\<g},\sigma_{\<f}\star\>\sigma_{\<g}}$
and $\Psi_{F\star\>G,F\star\>\sigma_{\<g}}$ respectively.
By \ref{prop:recC}, $c$ equals $\Psi_{H,F\star\>G}$. Thus all the maps
in the preceding diagram can be expressed as $\Psi_{-,-}$ suitably.
Now we check that each of the two paths from 
$f^!g^!$ to $(gf)^{\times}$ composes to
$\Psi_{H, \sigma_{\<f}\star\>\sigma_{\<g}}$
by repeatedly using the cocycle condition. This proves that
$\mu^{!}_{\times}$ is a pseudofunctorial map.

The remaining verifications needed to complete the construction 
of $\eC$ as a cover, such as defining $\mu^{!}_{\ssbox}$ and checking
that (b) and (c) of \ref{def:output1} hold, are carried out 
in a similar manner. For verifying commutativity of any 
diagram encountered, one uses \ref{thm:meiosis1} or the Propositions in 
\S\ref{subsec:recold} and \S\ref{subsec:recoldbc}
to first rewrite every isomorphism
encountered as $\Psi_{-,-}$ with suitable subscripts. Then, repeatedly
applying the cocycle rule gives the desired commutativity. 

\ssss
\label{sssec:output9} 
Proof of part (ii):
Let $\eC = ((-)^{!}, \mu^{!}_{\times},\mu^{!}_{\ssbox})$ and
$\eC' = ((-)^{\#}, \mu^{\#}_{\times},\mu^{\#}_{\ssbox})$. We first
show the existence of a map $\eC' \to \eC$. 

For any object $X$ in $\sfC$, let $S_X^! \colon X^! \to \cD_X$ 
and $S_X^{\#} \colon X^{\#} \to \cD_X$ 
be the functors associated to~$\eC$ and~$\eC'$ respectively.
Let $T_X \colon X^{\#} \to X^!$ be the
unique functor such that $S_X^!T_X = S_X^{\#}$. Such a functor exists 
because, by hypothesis $S_X^!$ is an isomorphism.

For any labeled sequence $\sigma \colon X \to Y$ we define 
$\tau_{\sigma} \colon T_X\sigma^{\#} \to \sigma^!T_Y$ to be the unique map 
making the following diagram commute where the remaining maps 
are the obvious ones resulting from $\eC,\eC'$ being covers. 
(see \ref{sssec:output4a})
\[
\xymatrix{
S_X^{!}T_X\sigma^{\#} \ar@{=}[d] \ar[rr]^{S_X^!(\tau_{\sigma})} 
 & & S_X^{!}\sigma^!T_Y \ar[r]^{\cong\quad} 
 & \sigma^{\bxt}S_Y^{!}T_Y \ar@{=}[d] \\
S_X^{\#}\sigma^{\#} \ar[rrr] & & & \sigma^{\bxt}S_X^{\#}
}
\]
In particular, if $\sigma$ consists of a single labeled map $(f,\sfP)$,
then the associated map
$T_Xf^{\#} \to f^{!}T_Y$ is the unique one 
determined by the commutativity of 
the diagram on the left 
in \ref{sssec:output2}; 
the analogous statement holds for $\sigma = (f,\sfO)$.

Before proceeding further we make a simplification. We shall
now assume that for any object $X$ in $\sfC$, it holds that
$X^{\#} = X^! = \cD_X$
and that each of $S^!_X,S^{\#}_X,T_X$ is the identity functor $\oneD{X}$. 
The proof in the general case is only a trivial but notationally cumbersome
modification of what we give below. 

In the simplified setup, we have so far defined, for every labeled 
sequence $\sigma$, a map $\sigma^{\#} \to \sigma^!$. 
Now consider the case of a general map $f$ in $\sfQ$. Let $\sigma$ be any 
labeled sequence such that $|\sigma|=f$. Then we obtain a map
$f^{\#} \iso \sigma^{\#} \to \sigma^{!} \iso f^!$
where the isomorphisms at either end are the canonical ones 
resulting from pseudofunctoriality. Our immediate aim is to 
verify that this map is independent of the choice of $\sigma$.

To that end, let $\sigma_1,\sigma_2$ be two sequences such that
$|\sigma_i|=f$. It suffices to verify that the outer border of the
following diagram of obvious natural maps commutes.
\[
\xymatrix{
& \sigma_1^{\#} \ar[ddd] \ar[rrrr] \ar[rrd] & & & & \sigma_1^{!} 
 \ar[lld]^{\cong} \ar[ddd] \ar[rd] \\
f^{\#} \ar[ru] \ar[rdd] & & & \sigma_1^{\bxt} 
 \ar[d]_{\Psi_{\sigma_2,\sigma_1}} & & & f^{!} \\
& & & \sigma_2^{\bxt} \\
& \sigma_2^{\#} \ar[rru] \ar[rrrr] 
 & & & & \sigma_2^{!} \ar[llu]_{\cong} \ar[ruu]
}
\]
The two triangles on the left and right sides commute by pseudofunctoriality.
The two trapeziums commute by Proposition \ref{prop:output6}. The remaining 
two triangles commute by definition.

It remains to verify that we now have a pseudofunctorial map
$(-)^{\#} \to (-)^!$ on $\sfQ$. To that end, let $X \xto{f} Y \xto{g} Z$
be maps in $\sfQ$. It suffices to verify that the following 
diagram commutes.
\[
\begin{CD}
f^{\#}g^{\#}  @>>> f^!g^! \\
@VVV @VVV \\
(gf)^{\#} @>>> (gf)^!
\end{CD}
\] 
Let $\sigma_f,\sigma_g,\sigma_{gf}$ be sequences such that 
$|\sigma_f| = f$, $|\sigma_g| = g$ and $|\sigma_{gf}| = gf$.
Then, using the definition of the horizontal maps in the preceding
diagram we expand it as follows.
\[
\begin{CD} 
f^{\#}g^{\#} @>>> \sigma_f^{\#}\sigma_g^{\#} @>>> 
 \sigma_f^{\bxt}\sigma_g^{\bxt} @>>> \sigma_f^!\sigma_g^! @>>> f^!g^! \\
@VVV @VVV @VV{\Psi_{\sigma_{gf},\sigma_f\star\sigma_g}}V @VVV @VVV \\
(gf)^{\#} @>>> \sigma_{gf}^{\#} @>>> \sigma_{gf}^{\bxt} 
 @>>> \sigma_{gf}^! @>>> (gf)^!
\end{CD}
\]
The two squares on either end commute for pseudofunctorial 
reasons. The ones in the middle commute 
by Proposition \ref{prop:output6}
for $\sigma_1 = \sigma_f\star\sigma_g$ and $\sigma_2 = \sigma_{gf}$.

Thus we have shown the existence of a map of covers $\eC' \to \eC$.
Now let us verify that there is a unique choice for such a map. 
Once again we assume, for simplicity, that for any $X$,
each of $S^!_X,S^{\#}_X,T_X$ is the identity functor $\oneD{X}$. 
Let $\tau_i \colon \eC' \to \eC$ for $i = 1,2$ be maps of covers.
Now recall that for any map $f$ in $\sfO$ or $\sfP$, the choice
for a map $f^{\#} \to f^!$ is uniquely determined by the 
commutativity of the corresponding diagram in \ref{sssec:output2}. 
Therefore, on $\sfO$ as well as on $\sfP$,
the restrictions of $\tau_i$ coincide.
Since any map in $\sfQ$ is a composite of maps in $\sfO$ or~$\sfP$,
the pseudofunctoriality of the map $\eC' \to \eC$
implies uniqueness in general.

\subsection{Proof of flat-base-change theorem}
\label{subsec:Thm2}

\ssss
\label{sssec:flat2}
We may upgrade the flat base-change isomorphisms $\sbeta^{-}_{-}$ 
of \S\ref{subsec:flat}[E3](i),(ii) to operate on the level of labeled
maps (with label $\sfO$ or $\sfP$) 
and labeled sequences in the obvious manner.
For any diagram~$\sfr$ as follows,
\[
\begin{CD}
W @>{u'}>> X \\
@V{G'}VV @VV{G}V \\
Z @>{u}>> Y
\end{CD}
\]
where $G,G'$ are labeled maps having the same label, $u,u' \in \sfF$,
and the underlying diagram is a fibered square, 
we may associate a base-change isomorphism 
$\sbeta_{\sfr} \colon u'{}^{\flat}F^{\bxt} \iso F'{}^{\bxt}u^{\flat}$
that is defined in the expected manner. More generally, 
for a diagram $\sS$ as follows, where $\sigma_1,\sigma_1'$ 
are labeled sequences, $u,u' \in \sfF$,
and the underlying diagram of ordinary maps is cartesian, 
\begin{equation}
\label{eq:flat3}
\xymatrix{
W \ar@{=>}[d]_{\sigma'} \ar[r]^{u'} & X \ar@{=>}[d]^{\sigma} \\
Z \ar[r]^u &  Y
}
\end{equation}
we associate a base-change isomorphism 
$\sbeta_{\sS} \colon u'{}^{\flat}\sigma^{\bxt} 
\iso \sigma'{}^{\bxt}u^{\flat}$, (also to be denoted
by $\sbeta_{\sigma,u}$)
that is defined via the $\sbeta_{-}$'s associated 
to each unit square in it. Note that~$\sbeta_{-}$
is transitive vertically and horizontally.
 
Any diagram, such as the preceding one shall henceforth 
be also called cartesian.

\begin{aprop}
\label{prop:flat4}
Consider two cartesian diagrams as follows where $u \in \sfF$,
$\sigma_i$ are labeled sequences and $|\sigma_1| = |\sigma_2|$,
$|\sigma_1'| = |\sigma_2'|$.
\[
\xymatrix{
W \ar@{=>}[d]_{\sigma_1'} \ar[r]^{u'} & X \ar@{=>}[d]^{\sigma_1} \\
Z \ar[r]^u &  Y
}
\qquad \qquad \qquad
\xymatrix{
W \ar@{=>}[d]_{\sigma_2'} \ar[r]^{u'} & X \ar@{=>}[d]^{\sigma_2} \\
Z \ar[r]^u &  Y
}
\]
Then the following diagram commutes.
\[
\begin{CD}
u'{}^{\flat}\sigma_1^{\bxt} @>{\qquad\sbeta_{\sigma_1,u}\qquad}>> 
 \sigma_1'{}^{\bxt}u^{\flat} \\
@V{\textup{via }\Psi_{\sigma_2,\sigma_1}}VV 
 @VV{\textup{via }\Psi_{\sigma_2',\sigma_1'}}V \\
u'{}^{\flat}\sigma_2^{\bxt} @>>{\qquad\sbeta_{\sigma_2,u}\qquad}> 
 \sigma_2'{}^{\bxt}u^{\flat} 
\end{CD}
\]
\end{aprop}
\begin{proof}
The proof uses the same kind of arguments we have been using so far
and so we only give a sketch. We proceed in three steps.

\emph{Step} 1. For the cartesian square of \eqref{eq:flat3}, we claim that 
under the further assumption that $X=Y$, $W=Z$, $|\sigma| = 1_X$,
$|\sigma'| = 1_Z$ and $u=u'$, the following diagram commutes.
\[
\xymatrix{
u^{\flat}\sigma^{\bxt} \ar[rd]_{\text{via }\Phi_{\sigma}} 
 \ar[rr]^{\sbeta_{\sigma,u}} 
 & & \sigma'{}^{\bxt}u^{\flat} \ar[ld]^{\text{via }\Phi_{\sigma'}} \\
& u^{\flat}
}
\] 
The proof of this proceeds along the same lines as \ref{lem:recold1}.

\emph{Step} 2. We claim that the cube lemma (\ref{lem:cocycle2}) 
also holds when 
the edges of the cube along a particular direction are all in $\sfF$
and for the corresponding faces we use $\sbeta$ instead of $\bbeta$.
Put differently, for the cartesian cube in \S\ref{subsec:flat}[E3](a), 
if we replace
$f,f_i,g,g_i$ by labeled sequences 
$\sigma,\sigma_i, \rho,\rho_i$ respectively, then 
the corresponding hexagon also commutes. This claim 
is proven by induction on the length of these
labeled sequences. The basis of induction is the case
when all the sequences have length one. If $\sigma$ and $\rho$ have 
different labels, then we conclude using [E3](a). If $\sigma$ and $\rho$ have
the same label, then we refer to the first half of the proof of 
the cube lemma. For the general induction argument we refer to 
the second half of the proof of the cube lemma.

\emph{Step} 3. To prove the Proposition, one now looks at
the diagram through which $\Psi_{\sigma_2,\sigma_1}$ is defined 
(see \ref{sssec:psi5z}). 
Upon taking a fibered product of this diagram with~$u$ one obtains
a 3-dimensional diagram (with a cube in it) relating the 
the definition of $\Psi_{\sigma_2,\sigma_1}$
with that of $\Psi_{\sigma_2',\sigma_1'}$.
Then one uses Steps 1 and 2 above, upon which the Proposition follows.
\end{proof}


\medskip
\ssss
Proof of Theorem \ref{thm:flat5}.
The proof proceeds along expected lines. First we 
lift the cartesian diagram $\sfr$ of the theorem to a cartesian
diagram $\sS$ such as in \eqref{eq:flat3} such that  
$|\sigma| = f$ and $|\sigma'| = f'$. (see \ref{lem:strcs1a}) 
Now we define $\sbeta^!_{\sfr}$ via 
\[
u'{}^{\flat}f^! \iso u'{}^{\flat}\sigma^! \iso u'{}^{\flat}\sigma^{\bxt}
\xto{\quad\sbeta_{\sS}\quad} 
\sigma'{}^{\bxt}u^{\flat} \iso \sigma'{}^{!}u^{\flat}
\iso f'{}^{\bxt}u^{\flat}.
\]
To check that this is independent of the choice of $\sigma$ 
it suffices to check that the following diagram obtained via
two liftings $\sS_1$ and $\sS_2$ of $\sfr$, with corresponding 
obvious notation, commutes.
\[
\xymatrix{
u'{}^{\flat}f^! \ar[rd] \ar[r] & u'{}^{\flat}\sigma_1^! \ar[r] \ar[d]
 & u'{}^{\flat}\sigma_1^{\bxt} \ar[r] \ar[d]
 & \sigma_1'{}^{\bxt}u^{\flat} \ar[r] \ar[d] & \sigma_1'{}^{!}u^{\flat} \ar[r]
 \ar[d] & f'{}^{\bxt}u^{\flat} \\
& u'{}^{\flat}\sigma_2^! \ar[r] & u'{}^{\flat}\sigma_2^{\bxt} \ar[r]
 & \sigma_2'{}^{\bxt}u^{\flat} \ar[r] & \sigma_2'{}^{!}u^{\flat} \ar[ru]
}
\]
The two triangles on either sides commute by pseudofunctoriality,
the square in the middle commutes by \ref{prop:flat4} and the remaining
squares commute for functorial reasons.

The verification that $\sbeta^!_{-}$ is transitive along both directions
is straightforward to verify. The diagrams in (ii) commute essentially
by definition, since for $f$ in $\sfP$ or $\sfO$ we may choose $\sigma$
to be $(f, \sfP)$ or $(f, \sfO)$ accordingly.

To check that $\sbeta_{-}^!$ is uniquely determined, first note that
if $f$ is in $\sfP$ or $\sfO$, then the choice for $\sbeta^!_{\sfr}$
is uniquely determined via (ii). In the general case, since
$f$ is a composite of maps in $\sfP$ and $\sfO$, vertical transitivity
of $\sbeta_{-}^!$ as required in (i) uniquely determines $\sbeta^!_{\sfr}$. 

For the last statement of the theorem, one first verifies the commutativity
in question for the case when $f$ is in $\sfO$ or $\sfP$. In the general
case, one then factors $f$ as a sequence of maps in $\sfO$ and $\sfP$.  

\subsection{Proof of the Variant Theorem}
\label{subsec:Thm3}

We give a proof of Theorem \ref{thm:variant3}.

We first show that the input conditions [A]--[D]
and [E1]--[E3] are satisfied, and then use Theorems \ref{thm:output4}
and \ref{thm:flat5}.

In view of \ref{sssec:variant1}(1) and \ref{sssec:variant1}(2),
only [D] and [E3](a),(b)
need to be verified. For [E3](a) we refer to the first half of the 
proof of the cube lemma. 

Let $X \xto{\;f\;} Y \xto{\;g\;} X$ be maps in $\sfC$ where 
$f \in \sfP$, $g \in \sfO$ and $gf = 1_X$. Then~$g$ 
is an isomrophism by~\ref{sssec:variant1}(4), and hence 
by \ref{sssec:variant1}(3), 
we have $g^{\ssbox} = g^{\times}$. We define the 
fundamental isomorphism $\phi_{g,f}\colon f^{\times}g^{\ssbox} \iso \oneD{X}$
to be the composition of the following canonical isomorphisms
\[
f^{\times}g^{\ssbox} = 
f^{\times}g^{\times} \xto{\;C^{\times}_{f,g}\;} (1_X)^{\times} = \oneD{X}.
\]

To verify [D](i)(a), we expand the diagram in question as follows.
\[
\begin{CD}
f'{}^{\times}g'{}^{\ssbox}h^{\times} @>>> 
 f'{}^{\times}h'{}^{\times}g^{\ssbox} @>>> 
 h^{\times}f^{\times}g^{\ssbox} \\
@| @| @| \\
f'{}^{\times}g'{}^{\times}h^{\times} @>>> 
 f'{}^{\times}h'{}^{\times}g^{\times} @>>> 
 h^{\times}f^{\times}g^{\times} \\
@VVV @. @VVV  \\
\oneD{X'}h^{\times} @= h^{\times} @= h^{\times}\oneD{X} 
\end{CD}
\]
The bottom portion commutes by pseudofunctoriality. The 
square on the upper right corner commutes for functorial  
reasons. In the square on the upper left corner we may first cancel 
$f'{}^{\times}$ on the left from each vertex and then conclude 
using~\ref{sssec:variant1}(5).

In case of [D](i)(b) we expand the diagram in question as follows.
\[
\begin{CD}
f'{}^{\times}g'{}^{\ssbox}h^{\ssbox} @>>> 
 f'{}^{\times}h'{}^{\ssbox}g^{\ssbox} @>>> 
 h^{\ssbox}f^{\times}g^{\ssbox} \\
@| @| @| \\
f'{}^{\times}g'{}^{\times}h^{\ssbox} @>>> 
 f'{}^{\times}h'{}^{\ssbox}g^{\times} @>>> 
 h^{\ssbox}f^{\times}g^{\times} \\
@VVV @. @VVV  \\
\oneD{X'}h^{\ssbox} @= h^{\ssbox} @= h^{\ssbox}\oneD{X} 
\end{CD}
\]
The top left square commutes by  \ref{sssec:variant1}(5),
while the top
right one commutes for functorial reasons. The bottom portion 
commutes by transitivity of base-change.

The remaining compatibilities, [D](ii) and [E3](b) are also 
easy to prove and we leave them for the reader to verify.

Now we verify that in this setup, the notion of a cover 
as in \ref{def:output1} is the same as that in \ref{sssec:variant2}.
Let $\eC = ((-)^{!}, \mu^{!}_{\times},\mu^{!}_{\ssbox})$
be a cover in the sense of \ref{def:output1} and let $g \colon Y \to X$
be an isomorphism in $\sfC$. Let $f \colon X \to Y$ be the 
inverse isomorphism. In view of our construction of the 
fundamental isomorphism $\phi_{f,g}$, by \ref{def:output1}(c), the outer
border of the following diagram commutes. 
\[
\xymatrix{ 
S_Xf^!g^! \ar[rr] \ar[d]  &  &  S_X\mathbf{1}_{X^!} \ar@{=}[dd] \\
f^{\times}S_Yg^! \ar[d] \ar[rd] \\
f^{\times}g^{\ssbox}S_X \ar@{=}[r] & f^{\times}g^{\times}S_X 
 \ar[r] & \mathbf{1}_{\cD_X}S_X
}
\]
The pentagon commutes because $(-)^! \to (-)^{\times}$ is a map 
of pseudofunctors. Thus the triangle commutes. Since $f^{\times}$
is an isomorphism, hence the maps $S_Yg^! \to g^{\times}S_X$
and $S_Yg^! \to g^{\ssbox}S_X$ are equal. Thus $\eC$ is also a cover
in the sense of~\ref{sssec:variant2}. Converse follows by reversing
the above arguments.
 
Now Theorem \ref{thm:variant3}
follows from Theorems \ref{thm:output4}
and \ref{thm:flat5}.

\section{Applications}
\label{sec:apps}
We give some applications of our abstract pasting results. 
In \S\ref{subsec:twisted}, we discuss pseudofunctorial properties
of the torsion twisted-inverse-image functor $(-)^!$
of Grothendieck duality over noetherian formal schemes.
In \S\ref{subsec:prep}, by working with larger derived 
categories we correspondingly obtain a \emph{prepseudofunctorial} 
extension of $(-)^!$. In \ref{subsec:nonnoe}, we use Lipman's
results on duality for non-noetherian formal schemes
to define $(-)^!$ over a suitable non-noetherian setup.
In \S\ref{subsec:huang}, we look at Huang's construction 
in~\cite{Hu} of a pseudofunctor $(-)^{\#}$ of 
zero-dimensional modules over the category of residually
finitely generated maps of noetherian complete local rings.
In \S\ref{subsec:Deligne}, we compare our pasting result 
with Deligne's result in~\cite{De}.

\subsection{Upper shriek for noetherian formal schemes}
\label{subsec:twisted}

Let us recall the setup of Grothendieck duality 
over noetherian formal schemes. We phrase it 
in terms of our abstract input data.

\ssss
\label{sssec:twisted-2}
Set
\begin{align}
\sfC &= \text{The category of morphisms of noetherian formal schemes;} 
 \notag \\
\sfO &= \text{The subcategory of open immersions in $\sfC$;} \notag \\
\sfP &= \text{The subcategory of pseudoproper maps in $\sfC$; 
 (\cite[1.2.2]{AJL2})} \notag \\
\sfF &= \text{The subcategory of flat maps in $\sfC$.} \notag 
\end{align}

For any object $X$ in $\sfC$, set $\cD_{X} \set \Dqct^+(X)$, the
derived category of bounded-below complexes having quasi-coherent torsion
$\cO_X$-modules as homology. (\cite[\S1.2]{AJL2})

Let $f \colon X \to Y$ be a map in $\sfC$.
\begin{align}
&\text{If $f \in \sfP$, then set } f^{\times} \set 
 \text{the right adjoint to }\R f_* \colon \cD_{X} \to \cD_{Y}  
 \;\text{(\cite[Theorem 2(a)]{AJL2})}; \notag \\
&\text{If $f \in \sfO$, then set } f^{\ssbox} \set \mathbf{L}f^* = f^*; 
 \notag \\
&\text{If $f \in \sfF$, then set } f^{\flat} \set \R\iGp{X}f^* 
 \;\text{(\cite[\S1.2]{AJL2})}. \notag 
\end{align}
The pseudofunctorial structure of $(-)^{\ssbox}$ is the canonical one.
For $(-)^{\times}$, we use the obvious one inherited, via adjointness,
from the canonical (covariant) pseudofunctor for the derived direct image.   
In case of $(-)^{\flat}$, the comparison map $C^{\flat}_{f,g}$ for 
a pair of $\sfF$-maps $X \xto{f} Y \xto{g} Z$ is given by the natural map
\[
\R\iGp{X}f^*\R\iGp{Y}g^* \to \R\iGp{X}f^*g^* \iso \R\iGp{X}(gf)^*.
\]
By \cite[Proposition 5.2.8(c)]{AJL2}, $C^{\flat}_{f,g}$ is an 
isomorphism. Associativity 
of~$C^{\flat}_{-,-}$ is easy to verify.

Consider a cartesian square $\sfr$
as follows. 
\begin{equation}
\label{eq:twisted-1}
\begin{CD}
U @>{j}>> X \\
@V{g}VV @VV{f}V \\
V @>{i}>> Y
\end{CD}
\end{equation}
If $f,g \in \sfP$ and $i,j \in \sfF$, then we choose the 
flat-base-change isomorphism
\begin{equation}
\label{eq:twisted0}
\sbeta_{\sfr} \colon \R\iGp{U}j{}^{*}f^{\times} 
\iso g{}^{\times}\R\iGp{V}i^{*}
\end{equation}
to be the natural isomorphism of \cite[Definition 7.3]{AJL2}. 
Furthermore, if $i \in \sfO$, then 
we choose $\beta_{\sfr}$ to be $\sbeta_{\sfr}$ via the canonical
identifications $\R\iGp{U}j{}^{*} = j{}^{*}$ and
$\R\iGp{V}i^{*} = i^*$. The transitivity properties of 
$\beta_{-}$ and $\sbeta_{-}$ follow from the transitivity 
properties proved in \cite[\S 7.5]{AJL2}. 


So far we have shown that conditions (1) and (2) of \ref{sssec:variant1}
are achieved. Now we move onto (3) and (4).

Let $i$ be an isomorphism in $\sfC$. 
Then~$i$ is in $\sfP$ as well as in~$\sfO$
and the functors~$i^*$ and 
$\R i_* = i_*$ are both
left-adjoint and right-adjoint to each other. We therefore 
have $i^{\ssbox} = i^* = i^{\times}$.

Let $X \xto{\;f\;} Y \xto{\;g\;} X$ be maps in $\sfC$ where 
$f \in \sfP$, $g \in \sfO$ and $gf = 1_X$. Since $g$ is a surjective
as a map of sets $Y \to X$ therefore it is an isomorphism. 
In particular,~$f$ is also an isomorphism.


Finally \ref{sssec:variant1}(5) is tackled as follows.
\begin{alem}
\label{lem:twisted1}
Let $\sfr$ denote the cartesian square of \eqref{eq:twisted-1}
with $f,g \in \sfP$ and $i,j \in \sfF$. 
If \,$i,j$ are isomorphisms, then among the following two diagrams
the one on the left commutes
while if \,$f,g$ are isomorphisms, then the one on the right commutes.
\[
\begin{CD}
j^*f^{\times} @>{\beta_{\sfr}}>> g^{\times}i^* \\
@| @| \\
j^{\times}f^{\times} @>{\textup{via } (-)^{\times}}>> g^{\times}i^{\times} 
\end{CD}
\qquad \qquad
\begin{CD}
j^{\flat}f^{\times} @>{\sbeta_{\sfr}}>> g^{\times}i^{\flat} \\
@| @| \\
j^{\flat}f^{\flat} @>{\textup{via } (-)^{\flat}}>> g^{\flat}i^{\flat} 
\end{CD}
\]
\end{alem}

\begin{proof}
(i). We start with the diagram on the left.
Set $\pi \set fj = ig$. Then~$\pi$ is in $\sfP$. 
Consider the following diagram where the morphisms are induced 
via obvious pseudofunctorial ones or via the trace/cotrace maps of the form
$\R h_*h^{\times} \to \mathbf{1}$ or $\mathbf{1} \to h^{\times}\R h_*$  
arising from adjointness. (For $h=i$ 
or $h=j$ these maps are isomorphisms.)
\[
\xymatrix{
& \framebox{$j^*f^{\times}$} \ar[ld] \ar[rd]^{\cong} \\
g^{\times}\R g_*j^*f^{\times} \ar[d]_{\cong} \ar[r]^{\cong} 
 & g^{\times}\R g_*\pi^{\times}
 \ar[d]_{\cong}^{\hspace{6em} \sbsq_1} & \pi^{\times} \ar[l] \ar[r]^{a\qquad} 
 & \pi^{\times}\R\pi_*\pi^{\times} \ar[ld]^{\cong} \ar[dd]^b \\
g^{\times}i^*i_*\R g_*j^*f^{\times} \ar[d]_{\cong} \ar[r]^{\cong} 
 & g^{\times}i^*i_*\R g_*\pi^{\times} \ar[d]_{\cong} \ar[r]^{\cong} 
 & g^{\times}i^*\R\pi_*\pi^{\times} \ar[d] \\
g^{\times}i^*\R f_*j_*j^*f^{\times} \ar[rd]_{\cong}^{\hspace{6em} \sbsq_2} 
 \ar[r]^{\cong} & g^{\times}i^*\R f_*j_*\pi^{\times} \ar[ru]_{\cong} 
 & \framebox{$g^{\times}i^*$} \ar[r]_{\cong} & \pi^{\times} \\
& g^{\times}i^*\R f_*f^{\times} \ar[ru] 
}
\]
It suffices to check that the preceding diagram commutes because,
among the two paths along the outer border between the two framed vertices,
the western path gives the base-change isomorphism defined 
in \cite[7.3]{AJL2}, while the 
eastern one gives the pseudofunctorial isomorphism in view of 
$ba = 1_{\pi^{\times}}$ (a formal consequence of 
adjointness between~$\pi^{\times}$ and~$\R\pi_*$).

Of the subdiagrams, only $\sbsq_i$ need an explanation since the others
commute for functorial reasons. But both $\sbsq_1,\sbsq_2$ are exercises
in comparing a composition of adjoints with adjoint of a composition;
for $\sbsq_1$ we first cancel $\pi^{\times}$ on the right and for 
$\sbsq_2$ we first cancel $g^{\times}i^*$ on the left.

(ii). The approach here for the other diagram of the lemma is similar to 
(i) except for the complications
introduced by the presence 
of the derived torsion functor $\R\iGp{-}$. Set $\pi \set fj = ig$. 
Then $\pi$ is in~$\sfF$. Consider the diagram in \eqref{eq:twisted2}
below consisting of obvious natural maps where we use $f^{\times} = f^*$ 
and $g^{\times} = g^*$. (Note that 
some of the functors involved operate only on the larger 
``non-torsion'' categories
$\Dqc(-)$ instead of $\Dqct^+(-)$. However the composite functor at each 
vertex still operates on $\Dqct^+(-)$.)
\begin{figure}
\rotatebox{-90}
{\begin{minipage}{8.5in}\vspace{-3mm}
\begin{equation}
\label{eq:twisted2}
\xymatrix{
\framebox{$j^{\flat}f^*$} \ar[d]_{\cong} \ar@{=}[r] & \R\iGp{U}j^*f^* 
 \ar[dd]_{\cong} & \ar[l]_{\cong} g^*g_*\R\iGp{U}j^*f^* \ar[d]_{\cong} & 
 \ar[l]_{\cong} g^*\R\iGp{V}g_*j^*f^* \ar[d]_{\cong} & 
 \ar[l] g^*\R\iGp{V}i^*i_*g_*j^*f^* 
 \ar[d]_{\cong} & \ar[l]_{\cong} g^*\R\iGp{V}i^*f_*j_*j^*f^* \ar[d]_{\cong} \\
\pi^{\flat} \ar@{=}[rd]_{\cong} & & \ar[ld]_{\cong}^{\hspace{4em} \sbsq} 
 g^*g_*\R\iGp{U}\pi^* & 
 \ar[l]_{\cong} g^*\R\iGp{V}g_*\pi^* \ar[ld]_{\cong} & 
 \ar[l] g^*\R\iGp{V}i^*i_*g_*\pi^* \ar[d]_{\cong} \ar[ld]_{\cong} & 
 \ar[l]_{\cong} g^*\R\iGp{V}i^*f_*j_*\pi^* \ar[ld]_{\cong} & 
 g^*\R\iGp{V}i^*f_*f^* \ar[lu] \\
& \R\iGp{U}\pi^* & \ar[l]_{\cong} \R\iGp{U}g^*g_*\pi^* & 
 \ar[l] \R\iGp{U}g^*i^*i_*g_*\pi^* \ar[d]_{\cong} &
 g^*\R\iGp{V}i^*\pi_*\pi^* \ar[ld]_{\cong} & \ar[l] g^*\R\iGp{V}i^* 
 \ar[ld]_{\cong} \ar[ru]_{\cong} & 
 \ar@{=}[l] \framebox{$g^*i^{\flat}$} \ar[d]_{\cong} \\
& & \ar[lu]^a \R\iGp{U}\pi^*\pi_*\pi^* & \ar[l]_{\cong}^b 
 \R\iGp{U}g^*i^*\pi_*\pi^* & \ar[l]^c \R\iGp{U}g^*i^* & 
 \ar[l]_{\cong}^d \R\iGp{U}\pi^* & \ar@{=}[l] \pi^{\flat}
}
\qquad\qquad
\end{equation}
\bigskip \bigskip
\begin{equation}
\label{eq:twisted3}
\xymatrix{
\R\iGp{U}g^*\R\iGp{V}g_* \ar@{-}[d] \ar@{-}[rrr] 
 & & & \framebox{$g^*\R\iGp{V}g_*$} 
 & & & \ar@{-}[d] \ar@{-}[lll] g^*\R\iGp{V}g_*\R\iGp{U} \\ 
\framebox{$\R\iGp{U}g^*g_*$} \ar@{-}[dd] & & & 
 \R\iGp{U}g^*\R\iGp{V}g_*\R\iGp{U} 
 \ar@{-}[d] \ar@{-}[u] \ar@{-}[lll] \ar@{-}[rrr]
 & & & \framebox{$g^*g_*\R\iGp{U}$} \ar@{-}[dd] \\ 
& & & \R\iGp{U}g^*g_*\R\iGp{U} \ar@{-}[lllu] \ar@{-}[rrru] \ar@{-}[d] \\ 
\framebox{$\R\iGp{U}$} & & & \R\iGp{U}\R\iGp{U} \ar@{-}[lll] \ar@{-}[rrr] 
 & & & \framebox{$\R\iGp{U}$}
}
\end{equation}
\end{minipage}
}
\end{figure}

It suffices to check that \eqref{eq:twisted2} commutes since the northern
route on the outer
border between the two framed vertices gives $\sbeta_{\sfr}$ while 
the southern route gives the obvious pseudofunctorial isomorphism because
$abcd$ is the identity on~$\R\iGp{U}\pi^*$. 
Commutativity of the subdiagrams are verified
easily as in (i) except that for $\sbsq$ we use commutativity
of~\eqref{eq:twisted3} which consists of isomorphisms and 
whose bottom row gives the identity on $\R\iGp{U}$.
\end{proof}

Having shown that the abstract input conditions 
\ref{sssec:variant1}(1)--(5) are achieved,
we now obtain the output given in Theorem~\ref{thm:variant3}.
For convenience, we restate it here in a 
somewhat self-contained and concise way.

\begin{athm}
\label{thm:twisted4}
On the category $\sfQ$ of composites of compactifiable maps of 
noetherian formal schemes, there is a contravariant pseudofunctor~$(-)^!$,
unique up to a unique isomorphism, such that~$(-)^!$ takes
values in~$\Dqct^+(X)$ for any
object~$X$ in~$\sfQ$ and satisfies the following conditions.
\begin{enumerate}
\item Over the subcategory $\sfP$ of pseudoproper maps, 
$(-)^!$ gives the right adjoint to the
derived direct-image pseudofunctor $(-)_*$.
\item Over the subcategory $\sfO$ of open immersions, 
$(-)^!$ is the inverse-image pseudofunctor $(-)^*$.
\item For the fibered-product diagram $\sfr$ of \eqref{eq:twisted-1},
if \, $i,j$ are in~$\sfO$, then
the associated base-change isomorphism of \eqref{eq:twisted0}
equals, via the identifications in~\textup{(i)} and~\textup{(ii)},
the obvious pseudofunctorial isomorphism given by~$(-)^!$.
\end{enumerate}
\end{athm}

\begin{athm}
\label{thm:twisted5}
Let $(-)^!$ be the uniquely determined pseudofunctor of 
\textup{Theorem~\ref{thm:twisted4}}.
Let $(-)^{\flat}$ be the natural pseudofunctor on the 
category $\sfF$ of flat maps of noetherian formal schemes, that assigns 
to any $\sfF$-map $f \colon X \to Y$, the functor 
$\R\iGp{X}f^* \colon \Dqct^+(Y) \to \Dqct^+(X)$.  
Then for any fibered-product diagram $\sfr$ such as in 
\eqref{eq:twisted-1} where $f$ is in $\sfQ$ and $i$ in $\sfF$,
there is an associated flat-base-change isomorphism
$\sbeta_{\sfr}$ such that the following hold.
\begin{enumerate}
\item If $f$ is in $\sfP$, then $\sbeta_{\sfr}$, via the identification
in \textup{\ref{thm:twisted4}(i)}, is the isomorphism of \eqref{eq:twisted0}.
\item If $f$ is in $\sfO$, then $\sbeta_{\sfr}$, via the identification
in \textup{\ref{thm:twisted4}(ii)}, is the pseudofunctorial 
isomorphism of $(-)^{\flat}$ on $\sfF$.
\item $\sbeta_{\sfr}$ is transitive vis-\`{a}-vis horizontal
extensions of $\sfr$ via $\sfF$-maps or vertical extensions via $\sfP$-maps. 
\end{enumerate}
\end{athm}

\begin{arem}
\label{rem:twisted5a}
The pseudofunctor $(-)^!$ of Theorem~\ref{thm:twisted4} also satisfies
a universal property that makes it unique up to a unique isomorphism
of pseudofunctors. This universal property is precisely of 
$(-)^!$, as a perfect cover (\ref{sssec:variant2}), 
being final in the category of all covers. 
\end{arem}

Finally let us show
that $(-)^!$ of \ref{thm:twisted4} extends to \'{e}tale maps too.
Let $\sfO'$ be the subcategory of $\sfC$ consisting of 
formally smooth pseudo-finite-type separated morphisms of relative 
dimension~0 (see \cite[Definition 2.6.2]{LNS}). 
For our purposes, we define an \'{e}tale map to be a map 
in $\sfO'$. For any \'{e}tale $f \colon X \to Y$, set 
$f^{\ssbox'} = \R\iGp{X}f^*$. The pseudofunctorial structure 
of~$(-)^{\ssbox'}$ on~$\sfO'$ is the natural one, for instance, the
restriction of $(-)^{\flat}$ to $\sfO'$. The base-change-isomorphism
associated to the fibered product of a $\sfP$-map with an $\sfO'$-map
is obtained using $\sbeta_{-}$ of \eqref{eq:twisted0}.

\begin{athm}
\label{thm:twisted6}
The pseudofunctor $(-)^!$ of \textup{Theorem~\ref{thm:twisted4}}
extends to a pseudofunctor 
on the category $\sfQ'$ of 
composites of compactifiable maps and \'{e}tale maps 
such that the conditions \textup{(i)--(iii)} of~\textup{\ref{thm:twisted4}} 
hold with $\sfO'$ in place of $\sfO$. Any two such 
extensions are isomorphic via a unique isomorphism.
Finally, the flat-base-change isomorphisms of 
Theorem~\textup{\ref{thm:twisted5}} also extend uniquely
such that \textup{(i)--(iii)} of~\textup{\ref{thm:twisted5}} 
hold with~$\sfO'$ in place of~$\sfO$. 
\end{athm}
\begin{proof}
The only hurdle in using Theorem~\ref{thm:variant3} is that condition 
(4) of \ref{sssec:variant1} does not hold in this case. The
situation is salvaged as follows.


Let $X \xto{f} Y \xto{g} X$ be $\sfC$-maps
factoring $1_X$ such that $f \in \sfP$ and $g \in \sfO'$.
We claim that~$Y$ can be written as a disconnected union $Y= Y_1 \cup Y_2$
such that~$f$ is an isomorphism onto~$Y_1$. 
Indeed, since $g$ is separated and $gf = 1_X$, 
$f$ is necessarily a closed immersion.
Let~$\I$ be the ideal in~$\cO_Y$ corresponding to~$f$. By 
\cite[Prop.~2.6.8]{LNS}, $\I = \I^2$. Therefore, using the Nakayama lemma
over the local rings of $Y$,
we see that 
the two open subsets $Y_1 \set \{ y \in Y \;|\; \I_y = 0\}$ and 
$Y_2 \set \{ y \in Y \;|\; \I_y = \cO_{Y,y} \}$ disconnect~$Y$.
Clearly $f$ maps~$X$ isomorphically to~$Y_1$.

It follows that the restriction of~$g$ to~$Y_1$ is an isomorphism
so that on~$Y_1$ we may canonically identify~$g^{\times}$ with~$g^{\ssbox'}$.
Now we may use the same arguments as in the 
proof of Theorem~\ref{thm:variant3}
(see \S\ref{subsec:Thm3}), restricting to suitable connected components 
whenever necessary, to deduce that in the case of $\sfO'$ = \'{e}tale maps,
the original input conditions [A]--[D] and [E1]--[E3] are achieved. 
The theorem then follows.
\end{proof}

\subsection{Prepseudofunctorial extension of upper shriek}
\label{subsec:prep}
We now extend the pseudofunctor $(-)^!$ of Theorem~\ref{thm:twisted4} 
above by enlarging the
category~$\cD_X$ for each~$X$. What results is no longer a pseudofunctor 
but a prespeudofunctor in the sense of Lipman, see \cite[\S1.6]{Res}. 
Let us recall this notion.

\ssss
\label{sssec:prep1}
A \emph{prespeudofunctor} $(-)^{\#}$ on a category $\sfC$ has the same data
together with conditions as a pseudofunctor except for the following 
modification:
For an object $Z$ in $\sfC$ we no longer assume that
$1_Z^{\#} = \mathbf{1}_{Z^{\#}}$ and instead assign a map
$\delta_Z^{\#} \colon 1_Z^{\#} \to \mathbf{1}_{Z^{\#}}$, 
such that for any $\sfC$-map $X \xto{f} Y$,
the following composites are identity
\[
f^{\#} \xto{(C^{\#}_{1_X,f})^{-1}} 1_X^{\#}f^{\#} 
 \xto{\text{via } \delta_X^{\#}} f^{\#}, 
\qquad \qquad
f^{\#} \xto{(C^{\#}_{f, 1_Y})^{-1}} f^{\#}1_Y^{\#} 
\xto{\text{via } \delta_Y^{\#}} f^{\#}. 
\]
In particular, $1_Z^{\#}$ is an idempotent functor. 

A morphism of prepseudofunctors $(-)^{\#} \to (-)^!$ on $\sfC$
is defined the same way as that of pseudofunctors (see beginning of 
\S\ref{sec:abs}) except that 
we also have the additional requirement that for any 
object~$X$ in~$\sfC$, the following diagram commutes.
\[
\xymatrix{
S_X1_X^{\#} \ar[rr] \ar[rd] & & 1_X^{!}S_X \ar[ld] \\
& S_X
}
\]

\ssss
\label{sssec:prep2}
We are mainly interested in prepseudofunctors arising out of 
existing pseudofunctors and idempotent functors in the following way.

Let~$(-)^{\#}$ be a pseudofunctor on a category~$\sfC$.
Suppose that for every object~$X$ in~$\sfC$ there is a 
category~$X^{\wt{\#}}$ containing~$X^{\#}$ as a full
subcategory and that the canonical 
inclusion $i_X \colon X^{\#} \to X^{\wt{\#}}$
has a right adjoint $\iG{X} \colon X^{\wt{\#}} \to X^{\#}$ such that 
the canonical map $\alpha_X \colon \mathbf{1}_{X^{\#}} \to \iG{X}i_X$ 
induced by adjointness
is an isomorphism. For convenience, we suppress all occurrences
of the term $i_X$ and also think of $\iG{X}$ as a functor
$X^{\wt{\#}} \to X^{\wt{\#}}$ that takes $X^{\#}$ to itself. 
In particular, the canonical map 
$\beta_X \colon \iG{X} \to \mathbf{1}_{X^{\wt{\#}}}$ 
resulting from adjointness, when restricted to $X^{\#}$, gives 
an isomorphism with inverse $\alpha_X$. 
Moreover, for any $A \in X^{\wt{\#}}$, the composite
$\iG{X}A \xto{\alpha_X(\iG{X}A)} \iG{X}\iG{X}A 
\xto{\iG{X}(\beta_X(A))} \iG{X}A$ is the identity. 
Thus $\iG{X}$ is idempotent. 

In this setup we obtain a prepseudofunctor~$(-)^{\wt{\#}}$ on~$\sfC$ 
satisfying the following properties.

(i) For any map 
$f \colon X \to Y$ in $\sfC$, we have $f^{\wt{\#}} = f^{\#}\iG{Y}$.

(ii) For any pair of maps $X \xto{f} Y \xto{g} Z$ in $\sfC$, 
$C^{\wt{\#}}_{f,g}$ is the composite
\[
f^{\wt{\#}}g^{\wt{\#}} = f^{\#}\iG{Y}g^{\#}\iG{Z} \iso
f^{\#}g^{\#}\iG{Z} \iso (gf)^{\#}\iG{Z} = (gf)^{\wt{\#}}.
\] 

(iii) For any object $X$ in $\sfC$, $\delta_X^{\wt{\#}}$ is given by
$1_X^{\wt{\#}} = \iG{X} \xto{\;\beta_X\;} \mathbf{1}_{X^{\wt{\#}}}$. 

Let us verify that $(-)^{\wt{\#}}$ is indeed a prepseudofunctor.
For associativity of~$C^{\wt{\#}}_{-,-}$, we consider 
maps $X \xto{f} Y \xto{g} Z \xto{h} W$ and then conclude
by looking at the outer border of
the following diagram, where the bottommost square on the
left commutes by associativity of~$C^{\#}_{-,-}$ while the remaining
squares commute for functorial reasons.
\[
\begin{CD}
f^{\#}\iG{Y}g^{\#}\iG{Z}h^{\#}\iG{W} @>>> 
f^{\#}g^{\#}\iG{Z}h^{\#}\iG{W} @>>> (gf)^{\#}\iG{Z}h^{\#}\iG{W} \\
@VVV @VVV @VVV \\ 
f^{\#}\iG{Y}g^{\#}h^{\#}\iG{W} @>>> 
f^{\#}g^{\#}h^{\#}\iG{W} @>>> (gf)^{\#}h^{\#}\iG{W} \\
@VVV @VVV @VVV \\ 
f^{\#}\iG{Y}(hg)^{\#}\iG{W} @>>> 
f^{\#}(hg)^{\#}\iG{W} @>>> (hgf)^{\#}\iG{W} \\
\end{CD}
\]
For verifying that $\delta^{\wt{\#}}$ and $C^{\wt{\#}}_{-,-}$ are 
compatible we note that for any map $f \colon X \to Y$,
the following diagrams commute and the bottom rows compose to 
identity.
\[
\begin{CD}
f^{\wt{\#}} @>>> 1_X^{\wt{\#}}f^{\wt{\#}} @>>> f^{\wt{\#}} \\
@| @| @| \\
f^{\#}\iG{Y} @>>> \iG{X}f^{\#}\iG{Y} @>>> f^{\#}\iG{Y}
\end{CD}
\qquad
\begin{CD}
f^{\wt{\#}} @>>> f^{\wt{\#}}1_Y^{\wt{\#}} @>>> f^{\wt{\#}} \\
@| @| @| \\
f^{\#}\iG{Y} @>>> f^{\#}\iG{Y}\iG{Y} @>>> f^{\#}\iG{Y}
\end{CD}
\]

\ssss
\label{sssec:prep3}
In the situation of \ref{sssec:prep2}, we say that $(-)^{\wt{\#}}$
is induced by $(-)^{\#}$ and the family of idempotents $\iG{-}$.
We also write $(-)^{\wt{\#}} \set \langle(-)^{\#}, \iG{-}\rangle$.

Let $(-)^{\wt{\#}} = \langle(-)^{\#}, \iG{-}^{\#}\rangle$
and $(-)^{\wt{!}} = \langle(-)^{!}, \iG{-}^{!}\rangle$ 
be prepseudofunctors on~$\sfC$. Suppose there
exists a map $(-)^{\#} \to (-)^{!}$ and suppose that
for any object~$X$ in~$\sfC$, the natural map $S_X \colon X^{\#} \to X^!$
extends to a map $X^{\wt{\#}} \to X^{\wt{!}}$, also to be denoted as~$S_X$.
Then we obtain a \emph{$\varGamma$-induced} 
map $(-)^{\wt{\#}} \to (-)^{\wt{!}}$, defined as follows.

For any map $f \colon X \to Y$ in $\sfC$, we define 
$S_Xf^{\wt{\#}} \to f^{\wt{!}}S_Y$ through the following natural maps 
\[
S_Xf^{\wt{\#}} = S_Xf^{\#}\iG{Y}^{\#} \to f^!S_Y\iG{Y}^{\#} 
\xto{\;\;f^!(\theta_Y)\;\;} f^!\iG{Y}^{!}S_Y = f^{\wt{!}}S_Y,
\]
where $\theta_Y$ is the unique map $S_Y\iG{Y}^{\#} \to \iG{Y}^{!}S_Y$
induced via adjointness of $\iG{Y}^{!}$ by the natural map 
$S_Y\iG{Y}^{\#} \to S_Y$ so that the following
diagram of natural maps commutes.
\[
\xymatrix{
S_Y\iG{Y}^{\#} \ar[rr]^{\theta_Y} \ar[rd] & & \iG{Y}^{!}S_Y \ar[ld] \\
& S_Y
}
\]

Since, for any object $X$, $1_X^{\wt{\#}} = \iG{X}^{\#}$ and 
$1_X^{\wt{!}} = \iG{X}^{!}$, the compatibility 
between~$\delta_{-}^{\wt{\#}}$ and~$\delta_{-}^{\wt{!}}$ follows from the 
defining property of~$\theta_{-}$ in the preceding diagram. 
It only remains to verify that the
comparison maps~$C^{\wt{\#}}_{-,-}$ and~$C^{\wt{!}}_{-,-}$ are compatible.

Let $X \xto{f} Y \xto{g} Z$ be maps in $\sfC$. Then for proving that the 
following diagram commutes,
\[
\xymatrix{
S_Xf^{\wt{\#}}g^{\wt{\#}} \ar[rr] \ar[d] & & S_X(gf)^{\wt{\#}} \ar[dd] \\
f^{\wt{!}}S_Yg^{\wt{\#}} \ar[d] \\
f^{\wt{!}}g^{\wt{!}}S_Z \ar[rr] & & (gf)^{\wt{!}}S_Z
}
\]
upon expanding the maps as per their definitions,
it suffices to prove that the following diagram of obvious
natural maps commutes.
\[
\xymatrix{
S_Xf^{\#}\iG{Y}^{\!\!\#}g^{\#}\iG{Z}^{\!\!\#} \ar[d] \ar[r] 
 & S_Xf^{\#}g^{\#}\iG{Z}^{\!\!\#} \ar[d] \ar[r] 
 & S_X(gf)^{\#}\iG{Z}^{\!\!\#} \ar[dd] \\
f^{!}S_Y\iG{Y}^{\!\!\#}g^{\#}\iG{Z}^{\!\!\#} \ar[d] \ar[r] 
 & f^{!}S_Yg^{\#}\iG{Z}^{\!\!\#} \ar[dd] \\
f^{!}\iG{Y}^{\!\!!}S_Yg^{\#}\iG{Z}^{\!\!\#} \ar[d] \ar[ru] 
 & & (gf)^{!}S_Z\iG{Z}^{\!\!\#} \ar[dd] \\ 
f^{!}\iG{Y}^{\!\!!}g^{!}S_Z\iG{Z}^{\!\!\#} \ar[d] \ar[r] 
 & f^{!}g^{!}S_Z\iG{Z}^{\!\!\#} \ar[d] \ar[ru] \\
f^{!}\iG{Y}^{\!\!!}g^{!}\iG{Z}^{\!\!!}S_Z \ar[r] 
 & f^{!}g^{!}\iG{Z}^{\!\!!}S_Z \ar[r] & (gf)^{!}\iG{Z}^{\!\!!}S_Z  
}
\]
The triangle commutes by definition of $\theta_Y$ and the pentagon
commutes because $(-)^{\#} \to (-)^!$ is a map of pseudofunctors
by assumption. The remaining subdiagrams commute for functorial reasons.

\ssss
\label{sssec:prep4}
The pseudofunctor $(-)^!$ of Theorem \ref{thm:twisted4}
is extended to a prepseudofunctor as follows.
We use the same notation as in~\ref{sssec:twisted-2}. Thus $\sfC$ 
consists of noetherian formal schemes, etc..

For any $X$ in $\sfC$,
set $\wtcD_X \set \wtDqc^+(X)$, the triangulated
subcategory of $\D(X)$ whose objects are complexes $\cFb$ 
such that $\R\iGp{X}\cFb \in \Dqc^+(X)$---or equivalently, 
$\R\iGp{X}\cFb \in \Dqct^+(X)$, see \cite[5.2.9]{AJL2}.
Set $\iG{X}^! \set \R\iGp{X}$.
Then~$\iG{X}^!$ maps~$\wtcD_X$ to~$\cD_X$ and satisfies the conditions
in~\ref{sssec:prep2} above. Thus, the pseudofunctor~$(-)^!$ 
of~\ref{thm:twisted4} induces a prepseudofunctor
$(-)^{\wt{!}} \set \langle (-)^!, \iG{-}^! \rangle$ on the 
subcategory $\sfQ = \ov{\{\sfO,\sfP\}}$ of~$\sfC$.

\begin{athm}
\label{thm:prep5}
Let notation be as in \textup{\ref{sssec:prep4}}. Let 
$\eC = ((-)^{\natural}, \mu^{\natural}_{\times},\mu^{\natural}_{\ssbox})$ 
be a pseudofunctorial cover (in the sense of \ref{sssec:variant2}) and let
$(-)^{\wt{\natural}} = \langle(-)^{\natural}, \iG{-}^{\natural}\rangle$ 
be a prepseudofunctorial extension of $(-)^{\natural}$. Suppose that
for ever $X$ in $\sfC$, the natural functor 
$S_X \colon X^{\natural} \to \cD_X$
given by $\eC$ extends to a functor $X^{\wt{\natural}} \to \wtcD_X$.
Then there exists a unique map of prepseudofunctors 
$\epsilon \colon (-)^{\wt{\natural}} \to (-)^{\wt{!}}$ 
such that the following diagrams
of obvious $\varGamma$-induced maps commutes.
\[
\xymatrix{
(-)^{\wt{\natural}}\big|_{\sfP} \ar[rd] \ar[rr]^{\epsilon|_{\sfP}} 
 & & (-)^{\wt{!}}\big|_{\sfP} \ar[ld]^{\cong} \\
& (-)^{\wt{\times}}
}
\qquad\qquad
\xymatrix{
(-)^{\wt{\natural}}\big|_{\sfO} \ar[rd] \ar[rr]^{\epsilon|_{\sfP}}  
 & & (-)^{\wt{!}}\big|_{\sfO} \ar[ld]^{\cong} \\
& (-)^{\wt{\ssbox}}
}
\]
\end{athm}
\begin{proof}
Over $\sfP$ and $\sfO$, $\epsilon$ is determined by the commutativity
of the above two diagrams. Since every map in~$\sfQ$ is a composite 
of maps in~$\sfO$ and~$\sfP$, uniqueness of~$\epsilon$ follows.
For the existence, the universal property
of~$(-)^!$ being a perfect cover gives us a map $(-)^{\natural} \to (-)^!$
and we use its $\varGamma$-induced extension for~$\epsilon$.
\end{proof}

\begin{arem}
\label{rem:prep6}
For applications, it may be desirable to work over a 
smaller subcategory of $\sfC$. Suppose $\ov{\sfC}$ is a full subcategory
of~$\sfC$ such that if~$Y$ is an object in~$\ov{\sfC}$, and 
$f \colon X \to Y$ a map in~$\sfC$, then~$f$ is also in~$\ov{\sfC}$ 
so that~$X$ is also an object in~$\ov{\sfC}$. Then~$\ov{\sfC}$ is also 
stable under fibered products. We may then replace the subcategories 
$\sfO,\sfP, \sfQ$, etc., by their respective intersections
with~$\ov{\sfC}$ to obtain $\ov{\sfO},\ov{\sfP},\ov{\sfQ}$, etc. Then we see 
that the input conditions for pasting still hold over the 
restricted categories. For the corresponding output, the prepseudofunctorial
extensions also carry over. In summary, the conclusion of 
Theorem \ref{thm:prep5} holds over any such restricted setup too.
\end{arem}

\ssss
\label{sssec:prep7}
Flat base change also holds for $(-)^{\wt{!}}$.
For a flat map $f \colon X \to Y$ we continue to use
$f^{\flat} \set \R\iGp{X}f^*$ except that~$f^{\flat}$ is now
assumed to act on the larger category~$\wtcD_Y$.
Then $(-)^{\flat}$ is a prepseudofunctor over~$\sfF$ in the obvious way.
For the fibered square in \eqref{eq:twisted-1}, with $f$ in $\sfQ$ 
and $i$ in $\sfF$, we define
the base-change isomorphism $j^{\flat}f^{\wt{!}} 
\iso g^{\wt{!}}i^{\flat}$ via 
\[
j^{\flat}f^{\wt{!}} = j^{\flat}f^{!}\iG{Y}^! \iso g^{!}i^{\flat}\iG{Y}^!
\osi g^{!}i^{\flat} \osi g^{!}\iG{V}^!i^{\flat} = g^{\wt{!}}i^{\flat}.
\]
This flat-base-change isomorphism is also vertically
and horizontally transitive. In the case of base change of a proper map, it
agrees, via canonical identifications, with the isomorphism 
in~\cite[Definition 7.3]{AJL2}.

\subsection{Upper shriek for non-noetherian ordinary schemes}
\label{subsec:nonnoe}

In \cite{NNGD}, Lipman has shown that for many of the constructions 
of Grothendieck duality over ordinary schemes, 
the noetherian hypothesis which is frequently used,
can be eliminated. What has been lacking is the construction of a
pseudofunctor over some category of not-necessarily-noetherian schemes
in which both proper maps and \'{e}tale maps are brought in.
We use our abstract pasting results to obtain such a pseudofunctor
in the non-noetherian setup.
However our results are not entirely satisfactory since we need to impose
the additional hypothesis of flatness on the maps in the working category.
The problem stems from the fact that the notion of pseudo-coherence
of maps, essential for proving many of the duality constructions,
is not preserved under base change unless flatness is also assumed.
 
\ssss
\label{sssec:nonnoe1}
Here then is the setup as per the input requirements 
in \ref{sssec:variant1}. We set
\begin{align}
\sfC = \sfF &= \text{The category of flat pseudo-coherent 
morphisms of quasi-compact} \notag \\
&\hspace{1em}\text{quasi-separated schemes; \cite[\S1]{NNGD}}  \notag \\
\sfO &= \text{The subcategory of open immersions in $\sfC$;} \notag \\
\sfP &= \text{The subcategory of proper maps in $\sfC$.} \notag 
\end{align}
Note that every map in $\sfC$ is quasi-compact \cite[I, Cor.~6.1.10]{EGAI}. 
Also, pseudo-coherence is preserved under flat base change. Hence, both,
$\sfO$ and $\sfP$ are stable under base change by maps in $\sfC$.
 
For any object $X$ in $\sfC$, set $\cD_{X} \set \Dqc^+(X)$, the
derived category of bounded-below complexes having quasi-coherent 
$\cO_X$-modules as homology. 

Let $f \colon X \to Y$ be a map in $\sfC$.
\begin{align}
&\text{If $f \in \sfP$, then set } f^{\times} \set 
 \text{the right adjoint to }\R f_* \colon \cD_{X} \to \cD_{Y} ; \notag \\
&\text{If $f \in \sfO$, then set } f^{\ssbox} \set \mathbf{L}f^* = f^*; 
 \notag \\
&\text{If $f \in \sfF$, then set } f^{\flat} \set f^*. \notag 
\end{align}
The pseudofunctorial structures for $(-)^{\times},(-)^{\ssbox},(-)^{\flat}$
are the obvious ones.

For a cartesian square $\sfr$ in $\sfC$ 
as follows, 
\begin{equation}
\label{eq:nonnoe2}
\begin{CD}
U @>{j}>> X \\
@V{g}VV @VV{f}V \\
V @>{i}>> Y
\end{CD}
\end{equation}
if $f,g \in \sfP$ and $i,j \in \sfF$, then we choose the 
flat-base-change isomorphism
\begin{equation}
\label{eq:nonnoe3}
\sbeta_{\sfr} \colon j^{*}f^{\times} \iso g{}^{\times}i^{*}
\end{equation}
to be the natural isomorphism of \cite[4.3]{NNGD}.
The transitivity properties of $\sbeta_{-}$ follow from the transitivity 
properties proved in \cite[\S1.8]{Res}.

We have therefore shown that conditions (1) and (2) 
of~\ref{sssec:variant1} hold. 
Conditions~(3) and~(4) are easily verified just as 
in~\ref{sssec:twisted-2} above.
Continuing along similar lines,~(5) also follows from 
the arguments in \ref{lem:twisted1} since the definition
of the base-change isomorphisms are formally analogous. 

As a result we now obtain the following output.

\begin{athm}
\label{thm:nonnoe4}
\textup{(1)}. On the category $\sfQ_{\textup{nn}}$ of 
composites of \'{e}tale maps and 
proper flat pseudo-coherent maps of quasi-compact quasi-separated schemes,
there exists a pseudofunctor $(-)^!$, unique up to a unique isomorphism,
taking values in $\Dqc^+(X)$ for any scheme $X$, such that 
the following conditions hold.
\begin{enumerate}
\item Over the subcategory of proper maps, 
$(-)^!$ gives the right adjoint to the
derived direct-image pseudofunctor $(-)_*$.
\item Over the subcategory of \'{e}tale maps, 
$(-)^!$ is the inverse-image pseudofunctor $(-)^*$.
\item For the fiber-product diagram $\sfr$ of \eqref{eq:nonnoe2},
if \, $i,j$ are \'{e}tale, then
the associated base-change isomorphism of \eqref{eq:nonnoe3}
equals, via the identifications in~\textup{(i)} and~\textup{(ii)},
the obvious pseudofunctorial isomorphism given by~$(-)^!$.
\end{enumerate} 

\textup{(2)}. For the diagram $\sfr$ of \eqref{eq:nonnoe2}, if $f,g$ are in
$\sfQ_{\textup{nn}}$, then there is a flat-base-change isomorphism
$\sbeta^!_{\sfr}$. Moreover $\sbeta^!_{-}$ is horizontally and 
vertically transitive and is uniquely determined by the condition
that when $f,g$ are proper, $\sbeta^!_{\sfr}$ is the isomorphism
of \eqref{eq:nonnoe3} while if $f,g$ are \'{e}tale, then $\sbeta^!_{\sfr}$
is the obvious isomorphism resulting from pseudofunctoriality of $(-)^*$. 
\end{athm}

We remark here that the replacement of open immersions by \'{e}tale maps
in the above theorem is made possible by using essentially the
same arguments as in the proof of \ref{thm:twisted6}.

\subsection{A review of Huang's $(-)_{\#}$}
\label{subsec:huang}
We review Huang's construction of a pseudofunctor $(-)_{\#}$
in \cite{Hu} in light of our pasting results.
Since Huang's result concerns \emph{covariant} pseudofunctors,
so we work with the appropriate opposite category to obtain
contravariance. We now show that in the situation of \cite[Chp.~6]{Hu},
our input conditions [A]--[D] are attained.

Let $\sfC$ be the opposite category of the category of 
residually finitely generated local homomorphisms 
of noetherian complete local rings. 
Let~$\sfP$ be the subcategory consisting of surjective homomorphisms 
and~$\sfO$ the subcategory consisting of formally smooth homomorphisms.
For any (local) ring $A$ in $\sfC$, let $\cD_A$ be the 
category of $A$-modules having zero-dimensional support.


For a map $T \to S$ in~$\sfC$, the homomorphism
of rings actually goes $S \to T$. In order to switch between these
categories, we use an ordinary letter, say $f$ for the ring homomorphism
and the dotted letter $\dot f$ for the corresponding $\sfC$-map.
For a ring~$R$ in~$\sfC$, $m_R$ stands for its maximal ideal. 

For $\dot f \colon T \to S$ in $\sfO$, we set 
$\dot f^{\ssbox} \set H^{n}_{m_T}(\omega_{T/S} \otimes_S -)$ where 
$n=\dim(T/m_ST)$ and~$\omega_{T/S}$ is top exterior power 
of the universally complete module
of relative differentials of~$T$ over~$S$ (\cite[3.10]{Hu}). For 
$\dot f \colon T \to S$ in $\sfP$ we set $\dot f^{\times} \set \Hom_S(T,-)$.

Let $T \xto{\dot f} S \xto{\dot g} R$ be maps in $\sfC$. 
If $\dot f,\dot g$ are in $\sfP$, then
the comparison map~$C_{\dot f,\dot g}^{\times}$ is the one 
given by ``evaluation at 1''
while if $\dot f,\dot g$ are in~$\sfO$, then 
there are many natural candidates for $C_{\dot f,\dot g}^{\ssbox}$. 
One possible candidate is given in \cite[2.5]{Hu}, which, for $M \in \cD_R$ 
can be described in terms of generalized fractions as follows
\begin{equation*}
\tag{\dag} 
\genfrac{[}{]}{0pt}{}{\beta \otimes \genfrac{[}{]}{0pt}{}{
 \alpha\otimes m}{\bs}}{\bt} \xto{\quad} 
\genfrac{[}{]}{0pt}{}{(\beta \wedge \alpha) \otimes m)}{\bt,\bs},
\qquad
\{m \in M, \alpha \in \omega_{S/R}, \beta \in \omega_{T/S}\}, 
\end{equation*}
where~$\bs$ is a system of parameters of~$S/m_RS$
and~$\bt$ of~$T/m_ST$.


For the fibered product of a $\sfP$-map and an $\sfO$-map, the base-change
isomorphism is given in \cite[3.6, 3.9]{Hu}. 
The transitivity of this base-change is 
readily verified from its explicit description in terms of generalized 
fractions.

For $\sfC$-maps $R \xto{\dot f} S \xto{\dot g} R$ 
where $\dot f$ is in $\sfP$ and $\dot g$ in $\sfO$ and 
$gf = 1_R$,
the fundamental isomorphism $\phi_{\dot f,\dot g}$ is obtained using residues. 
Indeed, for any $M \in \cD_R$ and any 
regular system of parameters $\bs = s_1, s_2, \ldots, s_n$ of $S/m_RS$, 
the $R$-module $\dot f^{\times}\dot g^{\ssbox}M$ can be canonically identified 
with the $R$-submodule of $\dot g^{\ssbox}M$ consisting of elements 
written in terms of generalized fraction as 
$\genfrac{[}{]}{0pt}{}{
ds_1\wedge ds_2\wedge \cdots ds_n \otimes m}{s_1, s_2, \ldots, s_n}$
and $\phi_{\dot f,\dot g}$ sends this fraction to $m$. (see \cite[Chp.~5]{Hu}.)

When it comes to verifying the compatibilities 
in \S\ref{subsec:input}[D](i), one 
runs into trouble. While [D](i)(a) and [D](ii) hold, [D](i)(b) fails.
Using the descriptions of all the isomorphisms involved in terms
of generalized fractions, the reader may verify that for any  
cartesian diagram of $\sfC$-maps as follows where $\dot f \in \sfP$,
$\dot g,\dot h \in \sfO$,
\[
\begin{CD}
R' @>{\dot f'}>> S' @>{\dot g'}>> R' \\
@VV{\dot h}V @VV{\dot h'}V @VV{\dot h}V \\
R @>{\dot f}>> S @>{\dot g}>> R
\end{CD}
\]
the commutativity of the diagram in [D](i)(b) is off by a factor 
of $(-1)^{ab}$ where $a = \dim(S/m_RS)$ and $b$ is the transcendence
degree of the residue field extension $k_R \to k_{R'}$.

The remedy lies in choosing the correct signed modification 
of $C_{-,-}^{\ssbox}$ and this is what Huang does in \cite[6.10]{Hu}. 
In the situation of~(\dag), 
let us redefine $C_{\dot f,\dot g}^{\ssbox}$ to be~$(-1)^{cd}$ times the 
isomorphism given by~(\dag) where $c = \dim(S/m_RS)$
and $d$ is the transcendence degree of 
$k_S \to k_{T}$. Then this new definition of $C_{-,-}^{\ssbox}$
still gives a pseudofunctor on~$\sfO$ and now [D](i)(b) is satisfied. 
For the other compatibilities, the calculations do not change.
Thus all the input conditions are attained.

Therefore the calculations listed above are all that are necessary
to obtain a pseudofunctor on the whole of $\sfC$.
Any such pseudofunctor 
is necessarily isomorphic to Huang's $(-)^{\#}$ of \cite[6.12]{Hu}. 
The canonicity of Huang's construction comes from the limiting process 
he employs over various factorizations of a fixed map into a $\sfO$-map
followed by a $\sfP$-map.
Here we should point out that 
the isomorphisms used in the limiting process in \cite[Chp.~6]{Hu} 
are special cases of our isomorphism $\Psi_{-,-}$ defined 
in \ref{sssec:psi5z}.
 
\subsection{Comparison with Deligne's result}
\label{subsec:Deligne}
In \cite[pp.~303--318]{De}, Deligne gave 
an abstract criterion for pasting pseudofunctors
on two subcategories into one on the whole category. Here we briefly compare 
his result with ours.

Deligne's input conditions are stated below as De-A to De-C, 
in reference to ours. For the purpose of this subsection, we assume
that in our input conditions $\sfC = \sfQ$, the smallest subcategory 
containing $\sfO$ and $\sfP$.

De-{\bf A.} Same as [A] except that in place of [A](ii), we assume that 
every map~$f$ in~$\sfC$ admits a factorization $f = ip$ where 
$i \in \sfO$ and $p \in \sfP$.

De-{\bf B.} Same as [B].

De-{\bf C.} For every \emph{commutative} square $\sfr$ as follows, 
\[
\begin{CD}
U @>{j}>> X \\
@V{g}VV @VV{f}V \\
V @>{i}>> Y
\end{CD}
\]
there is an isomorphism 
$\beta^{\text{De}}_{\sfr} \colon 
j^{\ssbox}f^{\times} \iso g^{\times}i^{\ssbox}$ that is transitive 
in the usual way.

Deligne's pasting result says that under De-A, De-B and De-C,
there is a pseudofunctor on $\sfC$ which restricts to $(-)^{\times}$
on $\sfP$, $(-)^{\ssbox}$ on $\sfO$, and that is compatible 
with~$\beta^{\text{De}}_{-}$.

We remark here that there is also a dual variant. Consider the
following condition.

De-{\bf \u A.} In De-A, assume that every $f$ factors as $f=pi$ instead. 

\noindent Then pasting also occurs for input data of De-\u A, De-B and De-C.
This is seen by interchanging the role of $\sfO$, $(-)^{\ssbox}$
with $\sfP, (-)^{\times}$ in the original criterion.

Note that Deligne's result is incomparable with ours in that,
neither De-A nor De-\u A, implies or is implied by, our condition [A],
even in the presence of the other conditions.
Of course, De-C, trivially gives us [C]. 
We now give a sketch of how De-C alone also gives us [D].
 
For a sequence $X \xto{f} Y \xto{g} X$ where $f \in \sfP, g \in \sfO$
and $gf=1_X$, the candidate for the fundamental isomorphism 
$f^{\times}g^{\ssbox} \iso \oneD{X}$ is $\beta^{\text{De}}_{\sfr}$
corresponding to the following commutative square $\sfr$. 
\[
\begin{CD}
X @= X \\
@V{f}VV @| \\
Y @>>{g}> X
\end{CD}
\]
For verifying [D](i), one takes a product of this diagram with 
$X' \to X$. Then one notes that cube lemma (\ref{lem:cocycle2}) 
also holds for a \emph{commutative} cube of unit size 
with~$\beta^{\text{De}}_{-}$ in place of~$\beta_{-}$, 
the proof being the same as the first half of the
proof of~\ref{lem:cocycle2}, (see~\ref{sssec:unitcube}).
Such a cube lemma together with an idempotence rule 
(see~\ref{lem:cafpols4}(i)) gives us~[D](i). For~[D](ii)
one completes the diagram containing~$\sfr$ there to 
a $2 \times 2$ diagram, the new maps all being identity on~$X$.
Now transitivity of~$\beta^{\text{De}}_{-}$ gives~[D](ii).

A\textsc{cknowledgements}.
I thank Joe Lipman and Pramath Sastry for constant encouragement
and the stimulating discussions which led to this work. I am also grateful
to them for sharing notes which have been helpful in 
the preparation of this paper. 


\providecommand{\bysame}{\leavevmode\hbox to3em{\hrulefill}\thinspace}
\providecommand{\MR}{\relax\ifhmode\unskip\space\fi MR }
\providecommand{\MRhref}[2]{%
  \href{http://www.ams.org/mathscinet-getitem?mr=#1}{#2}
}
\providecommand{\href}[2]{#2}

\end{document}